\documentclass[11pt,journal]{IEEEtran}

\renewcommand{\baselinestretch}{1.5}

\usepackage[mathcal]{euscript}
\usepackage[latin1]{inputenc}
\usepackage{epsfig}
\usepackage{amsmath,amstext,amsfonts,amssymb}
\usepackage{epsfig,float}
\usepackage{graphicx}
\usepackage{url}

\usepackage[french]{babel}

\usepackage{shadow,color,pifont,times,rotate}
\newcommand{\EE}{{\mathbb{E}}}
\newcommand{\tr}{{\mathrm{tr}}}
\newcommand{\PP}{{\mathbb{P}}}
\newcommand{\bs}{\boldsymbol}
\newcommand{\diag}{{\mathrm{diag}}}
\newcommand{\overcirc}{\overset{\circ}}
\newtheorem{lemma}{Lemma}

\newtheorem{proposition}{Proposition}
\newtheorem{theorem}{Theorem}

\newtheorem{remark}{Remark}
\def\adots{
  \mathinner{\mkern1mu\raise1pt\hbox{.}\mkern2mu\raise4pt\hbox{.}
  \mkern2mu\raise7pt\vbox{\kern7pt\hbox{.}}\mkern1mu}}
\def\build#1_#2^#3{\mathrel{
\mathop{\kern 0pt#1}\limits_{#2}^{#3}}}
\onecolumn

\newcommand{\mb}{\mathbf}

\newcommand{\mc}{\mathcal}

\newcommand{\tQ}{{\bf Q}}
\newcommand{\tP}{{\bf P}}
\newcommand{\tO}{\overline{\bf Q}_{*}}

\newcommand{\cerclenumero}[1]{{\ooalign {\hfil\raise .09ex\hbox{\scriptsize
(#1)}\hfil\crcr\mathhexbox20D}}}

\newtheorem{ass}{Assumption}

\usepackage[T1]{fontenc}    % Accents codés dans la fonte.
\usepackage[latin1]{inputenc}   % Accents 8 bits dans le fichier.
\usepackage{vmargin}        % Régler la taille de la feuille.
\usepackage{fancyhdr}       % Régler le titre courant et le bas de page.
\usepackage{calc}           % Faire des calculs sur les longueurs.
\usepackage{ifthen}         % Faire des tests if/then/else.
\usepackage{pifont}         % La police \ding.
\usepackage{supertabular}   % Les grands tableaux.
\usepackage{multicol}       % Plusieurs colonnes.
\usepackage{wrapfig}        % Dessins dans le texte.
\usepackage{fancybox}       % Boîtes avec une ombre pour les parties.
\usepackage{rotating}       % Pour tourner un texte.
\usepackage{xspace}         % Ajuster l'espace après des mots.

\usepackage{amsmath}    % Les symboles les plus fréquents.
\usepackage{amssymb}    % Des symboles.
\usepackage{amsfonts}   % Des fontes, eg pour \mathbb.
\usepackage{verbatim}   % Pour les codes sources en informatique.
\usepackage{mathrsfs}   % Des lettres majuscules cursives (\mathscr).
\usepackage{dsfont}

\usepackage{array}
\usepackage{curves}
\usepackage{epic}
\usepackage{eepic}
\usepackage{epsfig}
\usepackage{graphics}
\graphicspath{{PS/}}

\newcounter{ifnote}
\newif\ifremarque

\newcommand{\f}[2]{{\ensuremath{%
    \mathchoice%
    {\dfrac{#1}{#2}}
    {\dfrac{#1}{#2}}
    {\frac{#1}{#2}}
    {\frac{#1}{#2}}
}}}

\newcommand{\pa}[2][9]{%
    \ifthenelse{#1 = 0}
        {\ensuremath{#2}}{}%
    \ifthenelse{#1 = 1}
        {\ensuremath{(#2)}}{}%
    \ifthenelse{#1 = 2}
        {\ensuremath{\big(#2\big)}}{}%
    \ifthenelse{#1 = 3}
        {\ensuremath{\Big(#2\Big)}}{}%
    \ifthenelse{#1 = 4}
        {\ensuremath{\bigg(#2\bigg)}}{}%
    \ifthenelse{#1 = 5}
        {\ensuremath{\Bigg(#2\Bigg)}}{}%
    \ifthenelse{#1 = 9}
        {\ensuremath{\left(#2\right)}}{}%
}
\newcommand{\crochets}[2][9]{%
    \ifthenelse{#1 = 0}
        {\ensuremath{#2}}{}%
    \ifthenelse{#1 = 1}
        {\ensuremath{[#2]}}{}%
    \ifthenelse{#1 = 2}
        {\ensuremath{\big[#2\big]}}{}%
    \ifthenelse{#1 = 3}
        {\ensuremath{\Big[#2\Big]}}{}%
    \ifthenelse{#1 = 4}
        {\ensuremath{\bigg[#2\bigg]}}{}%
    \ifthenelse{#1 = 5}
        {\ensuremath{\Bigg[#2\Bigg]}}{}%
    \ifthenelse{#1 = 9}
        {\ensuremath{\left[#2\right]}}{}%
}
\newcommand{\accolades}[2][9]{%
    \ifthenelse{#1 = 0}
        {\ensuremath{#2}}{}%
    \ifthenelse{#1 = 1}
        {\ensuremath{\{#2\}}}{}%
    \ifthenelse{#1 = 2}
        {\ensuremath{\big\{#2\big\}}}{}%
    \ifthenelse{#1 = 3}
        {\ensuremath{\Big\{#2\Big\}}}{}%
    \ifthenelse{#1 = 4}
        {\ensuremath{\bigg\{#2\bigg\}}}{}%
    \ifthenelse{#1 = 5}
        {\ensuremath{\Bigg\{#2\Bigg\}}}{}%
    \ifthenelse{#1 = 9}
        {\ensuremath{\left\{#2\right\}}}{}%
}
\newcommand{\abs}[2][9]{%
    \ifthenelse{#1 = 0}
        {\ensuremath{#2}}{}%
    \ifthenelse{#1 = 1}
        {\ensuremath{|#2|}}{}%
    \ifthenelse{#1 = 2}
        {\ensuremath{\big|#2\big|}}{}%
    \ifthenelse{#1 = 3}
        {\ensuremath{\Big|#2\Big|}}{}%
    \ifthenelse{#1 = 4}
        {\ensuremath{\bigg|#2\bigg|}}{}%
    \ifthenelse{#1 = 5}
        {\ensuremath{\Bigg|#2\Bigg|}}{}%
    \ifthenelse{#1 = 9}
        {\ensuremath{\left|#2\right|}}{}%
}
\newcommand{\demi}{\ensuremath{\f{1}{2}}\xspace}

\newcounter{Centers}
\newcommand{\leftcentersright}[4][2]{%
    \ifthenelse{\equal{\value{Centers}}{0}}{%
        \leftcentersrightonpage[#1]{#2}{#3}{#4}%
    }{
        \leftcentersrightonitems[#1]{#2}{#3}{#4}%
    }
}

\newcounter{vcenterstest}
\newlength{\leftlength}
\newlength{\rightlength}
\newlength{\vcentersskip}
\newcommand{\leftcentersrightonitems}[4][2]{%
    \settowidth{\leftlength}{#2}%
    \settowidth{\rightlength}{#4}%
    \setcounter{vcenterstest}{#1}%
    \ifthenelse{\value{vcenterstest} = 0}
        {\setlength{\vcentersskip}{0pt}}{}%
    \ifthenelse{\value{vcenterstest} = 1}
        {\setlength{\vcentersskip}{\smallskipamount}}{}%
    \ifthenelse{\value{vcenterstest} = 2}
        {\setlength{\vcentersskip}{\medskipamount}}{}%
    \ifthenelse{\value{vcenterstest} = 3}
        {\setlength{\vcentersskip}{\bigskipamount}}{}%
    \ifthenelse{\value{vcenterstest} = 4}
        {\setlength{\vcentersskip}{1cm}}{}%
    \vskip\vcentersskip
    \noindent#2\hskip-\leftlength%
    \hfill#3\hfill%
    \mbox{}\hskip-\rightlength#4%
    \vskip\vcentersskip%
}

\newlength{\calculskip}
\newcommand{\calculvskip}[1]{%
  \ifthenelse{#1 = 0}{\setlength{\calculskip}{0pt}}{}%
  \ifthenelse{#1 = 1}{\setlength{\calculskip}{\smallskipamount}}{}%
  \ifthenelse{#1 = 2}{\setlength{\calculskip}{\medskipamount}}{}%
  \ifthenelse{#1 = 3}{\setlength{\calculskip}{\bigskipamount}}{}%
  \ifthenelse{#1 = 4}{\setlength{\calculskip}{1cm}}{}%
  \vskip\calculskip
}

\newcommand{\leftcentersrightonpage}[4][2]{%
    \settowidth{\leftlength}{#2}%
    \settowidth{\rightlength}{#4}%
    \calculvskip{#1}
    \noindent#2\hskip-\leftlength%
    \ifremarque{}\else\hskip\linewidth\hskip-\textwidth\fi%
    \hfill{#3}\hfill
    \mbox{}\hskip-\rightlength#4%
    \vskip\calculskip%
}

    {\ensuremath{\left \{ \hskip -1.5 mm \begin{array}{#1@{\ =\ }l}}}%
    {\end{array}\right.}
    {\ensuremath{\begin{array}[t]{#1@{\ =\ }l}}}%
    {\end{array}}
    {\ensuremath{\left \{ \hskip -1.5 mm%
     \begin{array}{#1@{\ }c@{\ }l}}}%
    {\end{array}\right.}
    {\ensuremath{\left \{ \hskip -1.5 mm \begin{array}{#1@{\quad}l}}}%
    {\end{array}\right.}
    {\ensuremath{\begin{array}[t]{#1@{\ \leq\ }l}}}%
    {\end{array}}
    {\ensuremath{\begin{array}[t]{#1@{\ }c@{\ }l}}}%
    {\end{array}}

\makeatletter

\newbox\bk@bxb
\newbox\bk@bxa
\newif\if@bkcont
\newif\ifbkcount
\newcount\bk@lcnt

\def\breakboxskip{2pt}
\def\breakboxparindent{1.8em}
\def\margesep{1cm}  % Écart entre la marge de gauche et le filet.
\def\intervalle{1mm}    % Écart supplémentaire entre le filet et le texte.

\def\filet{\vskip\breakboxskip\relax
\setbox\bk@bxb\vbox\bgroup \advance\linewidth -\fboxrule \advance\linewidth
-\margesep \advance\linewidth -\intervalle \advance\linewidth -\fboxsep
\hsize\linewidth\@parboxrestore
\parindent\breakboxparindent\relax}

\def\bk@split{%
\@tempdimb\ht\bk@bxb % height of original box
\advance\@tempdimb\dp\bk@bxb
\setbox\bk@bxa\vsplit\bk@bxb to\z@ % split it
\setbox\bk@bxa\vbox{\unvbox\bk@bxa}% recover height & depth of \bk@bxa
\setbox\@tempboxa\vbox{\copy\bk@bxa\copy\bk@bxb}% naive concatenation
\advance\@tempdimb-\ht\@tempboxa
\advance\@tempdimb-\dp\@tempboxa}% gap between two boxes

\def\bk@addfsepht{%
     \setbox\bk@bxa\vbox{\vskip\fboxsep\box\bk@bxa}}

\def\bk@addskipht{%
     \setbox\bk@bxa\vbox{\vskip\@tempdimb\box\bk@bxa}}

\def\bk@addfsepdp{%
     \@tempdima\dp\bk@bxa
     \dp\bk@bxa\@tempdima}

\def\bk@addskipdp{%
     \@tempdima\dp\bk@bxa
     \advance\@tempdima\@tempdimb
     \dp\bk@bxa\@tempdima}

\def\bk@line{%
    \hbox to \linewidth{\ifbkcount\smash{\llap{\the\bk@lcnt\ }}\fi
    \hskip\margesep
    \vrule \@width\fboxrule\hskip\fboxsep
    \hskip\intervalle
    \box\bk@bxa\hfil
    }}

\def\endfilet{\egroup
\ifhmode\par\fi{\noindent\bk@lcnt\@ne
\@bkconttrue\baselineskip\z@\lineskiplimit\z@ \lineskip\z@\vfuzz\maxdimen
\bk@split\bk@addfsepht\bk@addskipdp
\ifvoid\bk@bxb      % Only one line
\def\bk@fstln{\bk@addfsepdp
\vbox{\bk@line}}%
\else               % More than one line
\def\bk@fstln{\vbox{\bk@line}\hfil
\advance\bk@lcnt\@ne \loop
 \bk@split\bk@addskipdp\leavevmode
\ifvoid\bk@bxb      % The last line
 \@bkcontfalse\bk@addfsepdp
 \vtop{\bk@line}%
\else               % 2,...,(n-1)
 \bk@line
\fi
 \hfil\advance\bk@lcnt\@ne
\if@bkcont\repeat}%
\fi \leavevmode\bk@fstln\par}\vskip\breakboxskip\relax}

\bkcountfalse

\makeatother

        {\end{filet}\bigskip\remarquefalse}

\reversemarginpar \marginparwidth 3.1cm

\newcounter{Note}[page]
\newcounter{notesimple}
\newcounter{notemark}
\newcounter{notemarkref}
\newcounter{notetext}
\newcounter{notetextref}
\newcommand{\notemark}{%
    \ifthenelse{\equal{\value{notemark}}{\value{notemarkref}}}{%
            \refstepcounter{notemark}%
            \refstepcounter{notemarkref}%
            \setcounter{Note}{\value{notemark}}%
       }{%
            \refstepcounter{notemark}%
            \setcounter{Note}{\value{notemark}}%
    }%
    \setcounter{notetextref}{\value{notemarkref}}%
    \begin{picture}(0,0)%
            \put(-3,-3){\LARGE\theNote}
    \end{picture}\xspace
}

\newcommand{\notetextbase}[1]{%
        \ifthenelse{\value{ifnote} = 1}{%
        \marginpar{\notedebasenumero{#1}}%
        }{}%
}

\newcommand{\notetext}[1]{
    \ifthenelse{\equal{\value{notetext}}{\value{notetextref}}}{%
            \refstepcounter{notetext}%
            \refstepcounter{notetextref}%
            \setcounter{Note}{\value{notetext}}
        }{
            \refstepcounter{notetext}%
            \setcounter{Note}{\value{notetext}}
    }%
\setcounter{notemarkref}{\value{notetextref}}%
\notetextbase{#1}%
}

\newcommand{\notedebase}[1]{%
    \raggedright%
        \scriptsize%
        \vskip-0.5\baselineskip
        \rule[1.4mm]{\linewidth}{0.5pt} \\%
        \vskip-0.5\baselineskip
        #1 \\%
        \vskip-0.5\baselineskip
        \rule[1.4mm]{\linewidth}{0.5pt} \\%
}

\newcommand{\notedebasenumero}[1]{%
    {\LARGE\theNote}
    \notedebase{#1}%
}

\newcommand{\notesimple}[1]{%
    \ifthenelse{\value{ifnote} = 1}{%
    \marginpar{\notedebase{#1}}%
    }{}%
}

\newcommand{\notenumero}[1]{%
        \ifthenelse{\value{ifnote} = 1}{%
        \refstepcounter{Note}%
        \setcounter{notemarkref}{\value{Note}}%
        \setcounter{notetextref}{\value{Note}}%
        \setcounter{notemark}{\value{Note}}%
        \setcounter{notetext}{\value{Note}}%
        \begin{picture}(0,0)%
        \put(-3,-3){\LARGE\theNote}
        \end{picture}%
        \ifthenelse{\isodd{\value{Note}}}{%
            \protect\reversemarginpar%
            \marginpar[{\notedebasenumero{#1}}]{\notedebasenumero{#1}}}{%
            \protect\normalmarginpar%
            \marginpar[{\notedebasenumero{#1}}]{\notedebasenumero{#1}}}%
        }{}%
}

\newcommand{\note}[2][1]{%
    \ifthenelse{\equal{#1}{0} \or \equal{#1}{simple}}{%
                \setcounter{notesimple}{1}%
        }{}%
    \ifthenelse{\value{notesimple} = 0}{%
        \notenumero{#2}%
    }{%
        \notesimple{#2}%
    }
}

\begin{document}

\title{\Large {\bf On the Capacity Achieving Covariance Matrix for
    Rician MIMO Channels: An Asymptotic Approach} }

\author{J. Dumont, W. Hachem, S. Lasaulce, Ph. Loubaton and J. Najim\\
\today
\thanks{This work was partially supported by the ``Fonds National de la Science'' via the ACI program ``Nouvelles Interfaces des Math\'ematiques'', project MALCOM number 205\ and by the IST  Network of Excellence NEWCOM, 
project number 507325.}% <-this % stops a space
\thanks{J. Dumont was with France Telecom and Université Paris-Est, France. 
{\tt dumont@crans.org}\ ,}
\thanks{W. Hachem and J. Najim are with
CNRS and T\'el\'ecom Paris, Paris, France.
        {\tt \{hachem,najim\}@enst.fr}\ ,}
\thanks{S. Lasaulce is with CNRS and SUPELEC, 
France. {\tt  lasaulce@lss.supelec.fr}\ ,}  
\thanks{P. Loubaton is with Université Paris Est, IGM LabInfo, UMR CNRS 8049, France. 
        {\tt loubaton@univ-mlv.fr}\ ,}%
}

\maketitle

% pour mettre une note: \note{} pas utilisable dans tous les environnements

%-----------------------------------------------------------------------------------

\begin{abstract}
  In this contribution, the capacity-achieving input covariance
  matrices for coherent block-fading correlated MIMO Rician channels
  are determined.  In contrast with the Rayleigh and uncorrelated
  Rician cases, no closed-form expressions for the eigenvectors of the
  optimum input covariance matrix are available. Classically, both the
  eigenvectors and eigenvalues are computed by numerical techniques.
  As the corresponding optimization algorithms are not very
  attractive, an approximation of the average mutual information is
  evaluated in this paper in the asymptotic regime where the number of
  transmit and receive antennas converge to $+\infty$ at the same
  rate.  New results related to the accuracy of the corresponding
  large system approximation are provided. An attractive optimization
  algorithm of this approximation is proposed and we establish that it
  yields an effective way to compute the capacity achieving covariance
  matrix for the average mutual information. Finally, numerical simulation
  results show that, even for a moderate number of transmit and
  receive antennas, the new approach provides the same results as
  direct maximization approaches of the average mutual information,
  while being much more computationally attractive.
\end{abstract}
%-----------------------------------------------------------------------------------

%-----------------------------------------------------------------------------------

%-----------------------------------------------------------------------------------

\IEEEpeerreviewmaketitle
%%%%%%%%%%%%%%%%%%%%%%%%%%%%%%%%%%%%%%%%%%%%%%%%%%%%%%%%%%%%%%%%%%%%%%%%%%

%%%%%%%%%%%%%%%%%%%%%%%%%%%%%%%%%%%%%%%%%%%%%%%%%%%%%%%%%%%%%%%%%%%%%%%%%%

%%%%%%%%%%%%%%%%%%%%%%%%%%%%%%%%%%%%%%%%%%%%%%%%%%%%%%%%%%%%%%%%%%%%%%%%%%

%%%%%%%%%%%%%%%%%%%%%%%%%%%%%%%%%%%%%%%%%%%%%%%%%%%%%%%%%%%%%%%%%%%%%%%%%%

%++++++++++++++++++++++++++++++++++++++++++++++++++++++++++++++++++++

%++++++++++++++++++++++++++++++++++++++++++++++++++++++++++++++++++++

%++++++++++++++++++++++++++++++++++++++++++++++++++++++++++++++++++++

\section{Introduction}
%++++++++++++++++++++++++++++++++++++++++++++++++++++++++++++++++++++

%++++++++++++++++++++++++++++++++++++++++++++++++++++++++++++++++++++

%++++++++++++++++++++++++++++++++++++++++++++++++++++++++++++++++++++

Since the seminal work of Telatar \cite{Telatar-ETT-99}, the advantage
of considering multiple antennas at the transmitter and the receiver in
terms of capacity, for Gaussian and fast Rayleigh fading single-user
channels, is well understood. In that paper, the figure of merit chosen
for characterizing the performance of a
coherent\footnote{Instantaneous channel state information is assumed
  at the receiver but not necessarily at the transmitter.}
communication over a fading Multiple Input Multiple Output (MIMO)
channel is the Ergodic Mutual Information (EMI). This choice will be
justified in section \ref{sec:asymptotic}. Assuming the knowledge of
the channel statistics at the transmitter, one important issue is then to
maximize the EMI with respect to the channel input distribution.
Without loss of optimality, the search for the optimal input
distribution can be restricted to circularly Gaussian inputs. The
problem then amounts to finding the optimum covariance matrix.

This optimization problem has been addressed extensively in the case
of certain Rayleigh channels. In the context of the so-called
Kronecker model, it has been shown by various authors (see e.g.
\cite{Goldsmith-Jafar-etal-03} for a review) that the eigenvectors of
the optimal input covariance matrix must coincide with the eigenvectors of
the transmit correlation matrix. It is therefore sufficient to
evaluate the eigenvalues of the optimal matrix, a problem which can be
solved by using standard optimization algorithms. Note that
\cite{Tulino-Verdu-03} extended this result to more general (non
Kronecker) Rayleigh channels.

Rician channels have been comparatively less studied from this point
of view. Let us mention the work \cite{Hoesli-Kim-Lapidoth-05} devoted
to the case of uncorrelated Rician channels, where the authors proved
that the eigenvectors of the optimal input covariance matrix are the
right-singular vectors of the line of sight component of the channel.
As in the Rayleigh case, the eigenvalues can then be evaluated by
standard routines. The case of correlated Rician channels is undoubtedly
more complicated because the eigenvectors of the optimum matrix have
no closed form expressions. Moreover, the exact expression of the EMI
being complicated (see e.g. \cite{Kang-Alouini-Rice-06}), both the
eigenvalues and the eigenvectors have to be evaluated numerically. In
\cite{Vu-Paulraj-05}, a barrier interior-point method is proposed and
implemented to directly evaluate the EMI as an expectation.  The
corresponding algorithms are however not very attractive because they
rely on computationally-intensive Monte-Carlo simulations.

In this paper, we address the optimization of the input covariance of
Rician channels with a two-sided (Kronecker) correlation. 
As the exact expression of the EMI is
very complicated, we propose to evaluate an approximation
of the EMI, valid when the number of transmit and receive antennas converge to
$+\infty$ at the same rate, and then to optimize this
asymptotic approximation. This will turn out to be a simpler problem.
The results of the present contribution have been presented in part in the 
short conference paper \cite{Dumont-Loubaton-Lasaulce-2006}. 

The asymptotic approximation of the mutual information has been
obtained by various authors in the case of MIMO Rayleigh channels, and
has shown to be quite reliable even for a moderate number of
antennas. The general case of a Rician correlated channel has recently
been established in \cite{HLN07} using large random matrix theory and
completes a number of previous works among which \cite{Chuah-Tse-02},
\cite{Tulino-Verdu-Book} and \cite{Moustakas-Simon-Sengupta-03}
(Rayleigh channels), \cite{Cottatelluci-Debbah-04} and
\cite{Moustakas-Simon-05} (Rician uncorrelated channels),
\cite{Dumont-Loubaton-Lasaulce-Debbah-05} (Rician receive correlated
channel) and \cite{Tarrico-06} (Rician correlated channels). Notice
that the latest work (together with \cite{Moustakas-Simon-Sengupta-03}
and \cite{Moustakas-Simon-05}) relies on the powerful but non-rigorous
replica method. It also gives an expression for the variance of the
mutual information. We finally mention the recent paper \cite{Tarrico-Riegler-2007}
in which the authors generalize our approach sketched in \cite{Dumont-Loubaton-Lasaulce-2006}
to the MIMO Rician channel with interference. The optimization algorithm of the large system
approximant of the EMI proposed in  \cite{Tarrico-Riegler-2007} is however different
from our proposal.   

In this paper, we rely on the results of \cite{HLN07} in which a
closed-form asymptotic approximation for the mutual information is
provided, and present new results concerning its accuracy. We then
address the optimization of the large system approximation w.r.t. the
input covariance matrix and propose a simple iterative maximization
algorithm which, in some sense, can be seen as a generalization to the
Rician case of \cite{Wen-Com-06} devoted to the Rayleigh context: Each
iteration will be devoted to solve a system of two nonlinear equations
as well as a standard waterfilling problem. Among the convergence
results that we provide (and in contrast with \cite{Wen-Com-06}): We
prove that the algorithm converges towards the optimum input
covariance matrix as long as it converges. We also prove that the
matrix which optimizes the large system approximation asymptotically
achieves the capacity. This result has an important practical range as
it asserts that the optimization algorithm yields a procedure that
asymptotically achieves the {\em true} capacity. Finally, simulation
results confirm the relevance of our approach.

The paper is organized as follows. Section \ref{sec:model} is devoted
to the presentation of the channel model and the underlying
assumptions. The asymptotic approximation of the ergodic mutual
information is given in section \ref{sec:asymptotic-emi}.  In section
\ref{sec:concavity}, the strict concavity of the asymptotic
approximation as a function of the covariance matrix of the input
signal is established; it is also proved that the resulting optimal
argument asymptotically achieves the true capacity. The maximization
problem of the EMI approximation is studied in section \ref{sec:algo}.
Validations, interpretations and numerical results are provided in
section \ref{sec:simulations}.

%++++++++++++++++++++++++++++++++++++++++++++++++++++++++++++++++++++

%++++++++++++++++++++++++++++++++++++++++++++++++++++++++++++++++++++

%++++++++++++++++++++++++++++++++++++++++++++++++++++++++++++++++++++

\section{Problem statement}
\label{sec:model}
%++++++++++++++++++++++++++++++++++++++++++++++++++++++++++++++++++++

%++++++++++++++++++++++++++++++++++++++++++++++++++++++++++++++++++++

%++++++++++++++++++++++++++++++++++++++++++++++++++++++++++++++++++++

\subsection{General Notations} In this paper, the notations $s$, ${\bf
  x}$, $\mb{M}$ stand for scalars, vectors and matrices, respectively.
As usual, $\|{\bf x} \|$ represents the Euclidian norm of vector ${\bf
  x}$ and $\| {\bf M} \|$ stands for the spectral norm of matrix ${\bf
  M}$. The superscripts $(.)^T$ and $(.)^H$ represent respectively the
transpose and transpose conjugate. The trace of $\mb{M}$ is denoted by
$\mathrm{Tr}(\mb{M})$. The mathematical expectation operator is
denoted by $\mathbb{E}(\cdot)$ and the symbols $\Re$ and $\Im$ denote
respectively the real and imaginary parts of a given complex number.
If $x$ is a possibly complex-valued random variable, $\mathrm{Var}(x)
= \mathbb{E}|x|^{2} - \left|\mathbb{E}(x)\right|^{2}$ represents the
variance of $x$.

All along this paper, $r$ and $t$ stand for the number of transmit and
receive antennas.  Certain quantities will be studied in the
asymptotic regime $t \rightarrow \infty$, $r \rightarrow \infty$ in
such a way that $\f{t}{r} \rightarrow c \in (0,+\infty)$. In order to
simplify the notations, $t \rightarrow +\infty$ should be understood
from now on as $t \rightarrow \infty$, $r \rightarrow \infty$ and
$\f{t}{r} \rightarrow c \in (0,+\infty)$. A matrix ${\bf M}_t$ 
whose size depends on $t$ is said to be uniformly
bounded if $\sup_{t} \| {\bf M}_t \| < +\infty$.

Several variables used throughout this paper depend on various parameters, e.g. 
the number of antennas, the noise level, the covariance matrix of the transmitter, etc.
In order to simplify the notations, we may not always mention all these dependencies.

%--------------------------------------------------------------------

\subsection{Channel model}
%--------------------------------------------------------------------

We consider a wireless MIMO link with $t$ transmit and $r$ receive
antennas. In our analysis, the channel matrix can possibly vary from
symbol vector (or space-time codeword) to symbol vector. The channel
matrix is assumed to be perfectly known at the receiver whereas the
transmitter has only access to the statistics of the channel. The
received signal can be written as
\begin{equation}
{\bf y}(\tau) = \mb{H}(\tau) {\bf x}(\tau) + {\bf z}(\tau)
\end{equation}
where ${\bf x}(\tau)$ is the $t\times 1$ vector of transmitted symbols
at time $\tau$, $\mb{H}(\tau)$ is the $r\times t$ channel matrix
(stationary and ergodic process) and ${\bf z}(\tau)$ is a complex
white Gaussian noise distributed as $N(0, \sigma^2
\mathbf{I}_r)$. For the sake of simplicity, we omit the time
index $\tau$ from our notations. The channel input is subject to a
power constraint $\mathrm{Tr}\left[\mathbb{E}({\bf x} {\bf x}^H)
\right] \leq t$. Matrix $\mb{H}$ has the following structure:
\begin{equation}
\label{eq:modeleH}
\mb{H} = \sqrt{\frac{K}{K+1}} {\bf A} + \frac{1}{\sqrt{K+1}} {\bf V}\ ,
\end{equation}
where matrix $\bf{A}$ is deterministic, ${\bf V}$ is a random matrix and constant 
$K \geq 0$ is the so-called Rician factor which expresses the relative strength of the direct and
scattered components of the received signal.
Matrix $\bf{A}$ satisfies $\frac{1}{r} \mathrm{Tr}({\bf A} {\bf A}^{H}) = 1$ while 
${\bf V}$ is given by 
\begin{equation}
\label{eq:modeleV}
{\bf V} = \f{1}{\sqrt{t}} \mb{C}^{\demi} \mb{W} \mb{\tilde{C}}^{\demi}\ ,
\end{equation}
where $\mb{W}=(W_{ij})$ is a $r\times t$ matrix whose entries are independent and identically distributed (i.i.d.)
complex circular Gaussian random variables ${\mathcal CN}(0,1)$, i.e. $W_{ij}= \Re W_{ij} +\mathbf{i} \Im W_{ij}$
where $\Re W_{ij}$ and $\Im W_{ij}$ are independent centered real Gaussian random variables with variance $\frac 12$.
The matrices $\mb{\tilde{C}} > 0$ and $\mb{C} > 0$ account for the transmit
and receive antenna correlation effects respectively and satisfy 
$\frac{1}{t} \mathrm{Tr}(\mb{\tilde{C}}) = 1$ and $\frac{1}{r} \mathrm{Tr}(\mb{C}) = 1$.
This correlation structure is often referred to as a separable or Kronecker correlation model. 
\begin{remark}
Note that no extra assumption related to the rank of the 
deterministic component ${\bf A}$ of the
channel is done. Generally, it is often assumed that $\bf{A}$ has rank one
(\cite{Goldsmith-Jafar-etal-03}, \cite{Lozano-Tulino-Verdu-03},
\cite{hansen-isit-2004}, \cite{Lebrun-etal-2006}, etc..) because of the relatively small path loss
exponent of the direct path. Although the rank-one assumption is often relevant,
it becomes questionable if 
one wants to address, for instance, a multi-user setup and determine the
sum-capacity of a cooperative multiple access or broadcast channel in the high
cooperation regime. Consider for example 
a macro-diversity situation in the downlink: Several base stations
interconnected \footnote{For example in a cellular system the base stations are
connected with one another via a radio network controller.} through ideal wireline
channels cooperate to maximize the performance of a given multi-antenna
receiver. Here the matrix ${\bf A}$ is likely to have a rank higher than one or
even to be of full rank: Assume that the receive array of antennas is linear 
and uniform. Then a typical structure for ${\bf A}$ is
\begin{equation}
\label{eq:exempleA} {\bf A} =  \f{1}{\sqrt{t}} \left[{\bf
a}(\theta_1), \ldots, {\bf a}(\theta_t) \right] {\bs \Lambda}\ ,
\end{equation}
where ${\bf a}(\theta) = (1, e^{i \theta}, \ldots, e^{ i (r-1) \theta})^{T}$
and ${\bf \Lambda}$ is a diagonal matrix whose entries represent the complex
amplitudes of the $t$ line of sight (LOS) components. 
\end{remark}

%--------------------------------------------------------------------

\subsection{Maximum ergodic mutual information} \label{sec:asymptotic}
%--------------------------------------------------------------------

We denote by ${\cal C}$ the cone of nonnegative Hermitian $t \times t$
matrices and by ${\cal C}_1$ the subset of all matrices $\mb{Q}$ of ${\cal C}$
for which $\f{1}{t} \mbox{Tr}({\bf Q}) = 1$. Let $\mb{Q}$ be an element of
${\cal C}_1$ and denote by $I(\mb{Q})$ the ergodic mutual information (EMI)
defined by:
\begin{equation}
\label{eq:defmutual}
 I(\mb{Q}) =
\mathbb{E}_{\bs{H}} \left[ \log \det \left( {\bf I}_r + \f{1}{\sigma^2} \mb{H}
\mb{Q} \mb{H}^H \right) \right].
\end{equation}
Maximizing the EMI with respect to the input covariance matrix $\mb{Q}
= \mathbb{E}({\bf x} {\bf x}^H)$ leads to the channel Shannon capacity
for \textit{fast} fading MIMO channels i.e. when the channel vary from
symbol to symbol. This capacity is achieved by averaging over channel
variations over time. 

%For slow fading MIMO channels, i.e. when the
%channel matrix remains constant over a certain block duration much
%smaller than the channel coherence time, such an averaging is not possible
%and one has to communicate at rates smaller than the ergodic capacity.
%The maximum EMI is therefore a rate upper bound for \textit{slow}
%fading MIMO channels and only a fraction of it can be
%achieved\footnote{This fraction is called the multiplexing gain in
%  \cite{tse-dmt-it-2003} where the authors introduced the famous
%  diversity multiplexing trade-off.}. A more suited performance metric
%to study slow-fading channels is the outage capacity whose computation
%would require the knowledge of the variance of the mutual information.
%This is beyond the aim of this paper where we limit ourselves to the
%calculation of an asymptotic approximation of the mean of the mutual
%information. The computations performed in this article can be seen as a first step
%toward the evaluation of the variance of the EMI and its outage probability. 

We will denote by $C_E$ the maximum  value of the EMI over the set
$\mc{C}_1$:
\begin{equation}
\label{eq:defcapa} C_E = \sup_{\mb{Q} \in {\cal C}_1} I(\mb{Q}).
\end{equation}
The optimal input covariance matrix thus coincides with the argument
of the above maximization problem. Note that $I: \mb{Q} \mapsto
I(\mb{Q})$ is a strictly concave function on the convex set ${\cal
  C}_1$, which guarantees the existence of a unique maximum ${\bf Q}_*$
(see \cite{Luenberger}). When $\mb{\tilde{C}} = \mb{I}_t$, $\mb{C} =
\mb{I}_r $, \cite{Hoesli-Kim-Lapidoth-05} shows that the eigenvectors
of the optimal input covariance matrix coincide with the
right-singular vectors of ${\bf A}$. By adapting the proof of
\cite{Hoesli-Kim-Lapidoth-05}, one can easily check that this result
also holds when $\mb{\tilde{C}} = \mb{I}_t$ and ${\bf C}$ and ${\bf A}
{\bf A}^{H}$ share a common eigenvector basis. Apart from these two
simple cases, it seems difficult to find a closed-form expression for
the eigenvectors of the optimal covariance matrix.  Therefore the
evaluation of $C_E$ requires the use of numerical techniques (see e.g.
\cite{Vu-Paulraj-05}) which are very demanding since they rely on
computationally-intensive Monte-Carlo simulations. This problem can be
circumvented as the EMI $I(\mb{Q})$ can be approximated by a simple
expression denoted by $\bar{I}({\bf Q})$ (see section
\ref{sec:asymptotic-emi}) as $t\rightarrow \infty$ which in turn will
be optimized with respect to ${\bf Q}$ (see section \ref{sec:algo}).

\begin{remark} Finding the optimum covariance matrix is useful in
  practice, in particular if the channel input is assumed to be
  Gaussian. In fact, there exist many practical space-time encoders
  that produce near-Gaussian outputs (these outputs are used as inputs
  for the linear precoder $\mb{Q}^{1/2}$). See for instance
  \cite{rekaya-isita-2004}.
\end{remark}
\subsection{Summary of the main results.}
The main contributions of this paper can be summarized as follows:
\begin{enumerate}
\item We derive an accurate approximation of $I({\bf Q})$ as $t
  \rightarrow +\infty$: $I({\bf Q}) \simeq \bar{I}({\bf Q})$
  where
\begin{equation}
\label{eq:expre-preliminaire-Ibarre}
\bar{I}({\bf Q}) = \log \mathrm{det} \left[ {\bf I}_t + {\bf G}(\delta({\bf Q}, \tilde{\delta}({\bf Q})) {\bf Q} \right]
+ i(\delta({\bf Q}), \tilde{\delta}({\bf Q}))
\end{equation}
where $\delta({\bf Q})$ and $\tilde{\delta}({\bf Q})$ are two positive
terms defined as the solutions of a system of 2 equations (see Eq.
(\ref{eq:canonique})). The functions ${\bf G}$ and $i$ depend on
$(\delta({\bf Q}), \tilde{\delta}({\bf Q}))$, $K$, ${\bf A}$, ${\bf C}$, 
$\tilde{{\bf C}}$, and on the noise variance $\sigma^{2}$. They
are given in closed form.

The derivation of $\bar{I}({\bf Q})$ is based on the observation
that the eigenvalue distribution of random matrix ${\bf H} {\bf Q}
{\bf H}^{H}$ becomes close to a deterministic distribution as $t \rightarrow
+\infty$. This in particular implies that if $(\lambda_i)_{1\le i\le r}$
represent the eigenvalues of ${\bf H} {\bf Q} {\bf
  H}^{H}$, then:
\[
\frac{1}{r} \log \mathrm{det} \left[  {\bf I}_r + \f{1}{\sigma^2} \mb{H}
\mb{Q} \mb{H}^H \right] = \frac{1}{r} \sum_{i=1}^{r} \log \left( 1 + \frac{\lambda_i}{\sigma^{2}} \right) 
\]
has the same behaviour as a deterministic term, which turns out to be
equal to $\frac{\bar{I}({\bf Q})}{r}$. Taking the mathematical
expectation w.r.t. the distribution of the channel, and multiplying by
$r$ gives $I({\bf Q}) \simeq \bar{I}({\bf Q})$.

The error term $I({\bf Q}) - \bar{I}({\bf Q})$ is shown to be of order
$O(\frac{1}{t})$.  As $I({\bf
  Q})$ is known to increase linearly with $t$, the relative error
$\frac{I({\bf Q}) - \bar{I}({\bf Q})}{I({\bf Q})}$ is of order
$O(\frac{1}{t^{2}})$.  This supports the fact that $\bar{I}({\bf Q})$ is
an accurate approximation of $I({\bf Q})$, and that it is relevant to
study $\bar{I}({\bf Q})$ in order to obtain some insight on $I({\bf
  Q})$.

\item We prove that the function ${\bf Q} \mapsto \bar{I}({\bf Q})$ is
  strictly concave on ${\cal C}_1$. As a consequence, the maximum of
  $\bar{I}$ over ${\cal C}_1$ is reached for a unique matrix $\tO$.
  We also show that $I(\tO) - I({\bf Q}_*) = O(1/t)$ where we recall that 
  ${\bf Q}_*$ is the capacity achieving covariance matrix. Otherwise stated,
  the computation of $\tO$ (see below) allows one to (asymptotically)
  achieve the capacity $I({\bf Q}_*)$.

\item
  We study the structure of $\tO$ and establish that $\tO$ is
  solution of the standard waterfilling problem:
\[
\max_{{\bf Q} \in {\cal C}_1} \log \mathrm{det} \left( {\bf I} + {\bf G}(\delta_{*}, \tilde{\delta}_{*}) {\bf Q} \right)\ ,
\]
where $\delta_{*} = \delta(\tO)$, $\tilde{\delta}_{*} = \tilde{\delta}(\tO)$ and 
\[ {\bf G}(\delta_{*}, \tilde{\delta}_{*}) = \frac{\delta_{*}}{K+1}
\tilde{{\bf C}} + \frac{1}{\sigma^{2}} \frac{K}{K+1} {\bf A}^H \left(
  {\bf I}_r+\frac{\tilde{\delta}_{*}}{K+1} \,{\bf C} \right)^{-1}{\bf
  A}\ .
\]

This result provides insights on the structure of the approximating
capacity achieving covariance matrix, but cannot be used to evaluate
$\tO$ since the parameters $\delta_{*}$ and
$\tilde{\delta}_{*}$ depend on the optimum matrix $\tO$. We
therefore propose an attractive iterative maximization algorithm of
$\bar{I}({\bf Q})$ where each iteration consists in solving a
standard waterfilling problem and a $2\times 2$ system characterizing
the parameters $(\delta, \tilde{\delta})$.
\end{enumerate}

%++++++++++++++++++++++++++++++++++++++++++++++++++++++++++++++++++++

%++++++++++++++++++++++++++++++++++++++++++++++++++++++++++++++++++++

%++++++++++++++++++++++++++++++++++++++++++++++++++++++++++++++++++++

\section{Asymptotic behavior of the ergodic mutual information}
\label{sec:asymptotic-emi}

In this section, the input covariance matrix $\mb{Q} \in {\cal C}_1$
is fixed and the purpose is to evaluate the asymptotic behaviour of
the ergodic mutual information $I(\mb{Q})$ as $t \rightarrow \infty$
(recall that $t \rightarrow +\infty$ means $t \rightarrow \infty$, $r
\rightarrow \infty$ and $t/r \rightarrow c \in (0,+\infty)$).

As we shall see, it is possible to evaluate in closed form an accurate
approximation $\bar{I}({\bf Q})$ of $I(\mb{Q})$.  The corresponding
result is partly based on the results of \cite{HLN07} devoted to the
study of the asymptotic behaviour of the eigenvalue distribution of
matrix $\bs{\Sigma} \bs{\Sigma}^{H}$ where $\bs{\Sigma}$ is given by
\begin{equation}
\label{eq:modele-equivalent-iid}
\bs{\Sigma} = {\bf B} + {\bf Y}\ ,
\end{equation}
matrix ${\bf B}$ being a deterministic $r \times t$ matrix, and ${\bf
  Y}$ being a $r \times t$ zero mean (possibly complex circular
Gaussian) random matrix with independent entries whose variance
is given by $\mathrm{E}|Y_{ij}|^{2} = \frac{\sigma^{2}_{ij}}{t}$.
Notice in particular that the variables $(Y_{ij};\ 1\le i \le r,\ 1\le
j \le t)$ are not necessarily identically distributed. We shall refer
to the triangular array $(\sigma^{2}_{ij};\ 1\le i \le r,\ 1\le j \le
t)$ as the variance profile of $\bs{\Sigma}$; we shall say that it is
separable if $\sigma^{2}_{ij} = d_i \tilde{d}_j$ where $d_i \geq 0$
for $1\le i\le r$ and $\tilde{d}_j \geq 0$ for $1\le j\le t$. Due to
the unitary invariance of the EMI of Gaussian channels, the study of
$I(\mb{Q})$ will turn out to be equivalent to the study of the EMI of
model (\ref{eq:modele-equivalent-iid}) in the complex circular
Gaussian case with a separable variance profile.

\subsection{Study of the EMI of the equivalent model (\ref{eq:modele-equivalent-iid}).}
\label{subsec:modele-equivalent}

We first introduce the resolvent and the Stieltjes transform
associated with ${\bf \Sigma \Sigma}^H$ (Section
\ref{subsub-resolvent}); we then introduce auxiliary quantities
(Section \ref{subsub-aux}) and their main properties; we finally
introduce the approximation of the EMI in this case
(Section \ref{subsub-approx}).

\subsubsection{The resolvent, the Stieltjes transform}\label{subsub-resolvent}

Denote by ${\bf S}(\sigma^{2})$ and 
$\tilde{{\bf S}}(\sigma^{2})$ the resolvents of matrices $\bs{\Sigma} \bs{\Sigma}^{H}$ and 
$ \bs{\Sigma}^{H} \bs{\Sigma}$ defined by: 
\begin{equation}
\label{eq:def-resolvente}
{\bf S}(\sigma^{2}) =  \left[\bs{\Sigma} \bs{\Sigma}^{H} + \sigma^{2} {\bf I}_{r} \right]^{-1},\qquad
\tilde{{\bf S}}(\sigma^{2}) =  \left[\bs{\Sigma}^{H} \bs{\Sigma} + \sigma^{2} {\bf I}_{t} \right]^{-1}\ .
\end{equation}
These resolvents satisfy the obvious, but useful property:
\begin{equation}
  \label{eq:borne-resolvente}
  {\bf S}(\sigma^{2}) \leq \frac{{\bf I}_r}{\sigma^{2}}, \qquad  \tilde{{\bf S}}(\sigma^{2}) \leq \frac{{\bf I}_t}{\sigma^{2}}\ .
\end{equation}
Recall that the Stieltjes transform of a nonnegative measure $\mu$ is
defined by $\int \frac{\mu(d\lambda)}{\lambda-z}$. The quantity
$s(\sigma^{2}) = \frac{1}{r} \mathrm{Tr}({\bf S}(\sigma^{2}))$
coincides with the Stieltjes transform of the eigenvalue distribution
of matrix $\bs{\Sigma} \bs{\Sigma}^{H}$ evaluated at point
$z=-\sigma^2$. In fact, denote by $(\lambda_i)_{1\le i\le r}$ its
eigenvalues , then:
\[
s(\sigma^{2}) = \frac{1}{r} \sum_{i=1}^{r} \frac{1}{\lambda_i + \sigma^{2}} = \int_{\mathrm{R}^{+}} \frac{\nu(d\lambda)}{\lambda + \sigma^{2}}\ ,
\]
where $\nu$ represents the eigenvalue distribution of $\bs{\Sigma} \bs{\Sigma}^{H}$ defined as the probability 
distribution:
\[
\nu = \frac{1}{r} \sum_{i=1}^{r} \delta_{\lambda_i}
\]
where $\delta_x$ represents the Dirac distribution at point $x$. The
Stieltjes transform $s(\sigma^{2})$ is important as the
characterization of the asymptotic behaviour of the eigenvalue
distribution of $\bs{\Sigma} \bs{\Sigma}^{H}$ is equivalent to the
study of $s(\sigma^{2})$ when $t \rightarrow +\infty$ for each
$\sigma^{2}$. This observation is the starting point of the approaches
developed by Pastur \cite{Mar-Pas-67}, Girko 
\cite{Girko-01}, Bai and Silverstein \cite{BaiSil98}, etc.

We finally recall that a positive $p \times p$ matrix-valued measure
$\bs{\mu}$ is a function defined on the Borel subsets of $\mathbb{R}$
onto the set of all complex-valued $p \times p$ matrices satisfying:
\begin{itemize}
\item[(i)] For each Borel set $B$, $\bs{\mu}(B)$ is a Hermitian
  nonnegative definite $p \times p$ matrix with complex entries;
\item[(ii)] $\bs{\mu}(0) = {\bf 0}$;
\item[(iii)] For each countable family $(B_n)_{n \in \mathrm{N}}$ of
  disjoint Borel subsets of $\mathrm{R}$,
\[
\bs{\mu}(\cup_{n} B_n) = \sum_{n} \bs{\mu}(B_n)\ .
\]
\end{itemize}
Note that for any nonnegative Hermitian $p\times p$ matrix ${\bf M}$,
then $\mathrm{Tr}({\bf M} \bs{\mu})$ is a (scalar) positive measure.
The matrix-valued measure $\bs{\mu}$ is said to be finite if
$\mathrm{Tr}(\bs{\mu}(\mathbb{R})) < + \infty$.

\subsubsection{ The auxiliary quantities $\beta, \tilde \beta$, ${\bf T}$ and ${\bf \tilde T}$}\label{subsub-aux}
We gather in this section many results of \cite{HLN07} that will be of help in the sequel.

\begin{ass}
\label{ass:normes2}
Let $({\bf B}_t)$ be a family of $r\times t$ deterministic matrices such that: 
$\sup_{t,i} \sum_{j=1}^{t} |B_{ij}|^{2} < +\infty, \; \sup_{t,j} \sum_{i=1}^{r} |B_{ij}|^{2} < +\infty$\ .
\end{ass}

\begin{theorem}
\label{theo:canonique-preparatoire}
Recall that ${\bf \Sigma}= {\bf B} +{\bf Y}$ and assume that 
${\bf Y} = \frac 1{\sqrt{t}}{\bf D}^{\frac 12} {\bf X}\,
\tilde{{\bf D}}^{\frac 12}$, where ${\bf D}$ and $\tilde{{\bf D}}$
represent the diagonal matrices ${\bf D} = \mbox{diag}(d_i,\, 1\le
i\le r)$ and $\tilde{{\bf D}} = \mbox{diag}(\tilde{d}_j,\ 1\le j\le
t)$ respectively, and where ${\bf X}$ is a matrix whose entries are
i.i.d. complex centered with variance one. The following facts hold true: 
\begin{itemize}
\item[(i)] ({\em Existence and uniqueness of auxiliary quantities}) For $\sigma^{2}$ fixed, consider the
system of equations:
\begin{equation}
  \label{eq:canonique-preparatoire}
\left\{
\begin{array}{l}
  \beta 
  = \displaystyle{ \frac{1}{t} \textrm{Tr} \left[ {\bf D} \left( \sigma^{2}({\bf I}_r+{\bf D} \tilde{\beta}) + {\bf B} ({\bf I}_t+\tilde{{\bf D}}\beta)^{-1} {\bf B}^{H} \right)^{-1} \right] }
  \\
  \tilde{\beta} = \displaystyle{\frac{1}{t} \textrm{Tr} \left[ \tilde{{\bf D}} \left( \sigma^{2}({\bf I}_t+\tilde{{\bf D}} \beta) + {\bf B}^{H} 
        ({\bf I}_r+ {\bf D} \tilde{\beta})^{-1} {\bf B} \right)^{-1} \right]}
\end{array}\right. .
\end{equation}
Then, the system (\ref{eq:canonique-preparatoire}) admits a unique couple of positive solutions 
$(\beta(\sigma^{2}), \tilde{\beta}(\sigma^{2}))$. 
Denote by ${\bf T}(\sigma^{2})$ and $\tilde{{\bf T}}(\sigma^{2})$ the following matrix-valued functions:
\begin{equation}
\label{eq:defT}
\left\{
  \begin{array}{ccc}
    {\bf T}(\sigma^{2}) & = & \left[ \sigma^{2}({\bf I} +  \tilde{\beta}(\sigma^{2}) {\bf D}) + {\bf B}({\bf I} + \beta(\sigma^{2}) \tilde{{\bf D}})^{-1} {\bf B}^{H} \right]^{-1} \\
    \tilde{{\bf T}}(\sigma^{2}) & = & \left[ \sigma^{2}({\bf I} + \beta(\sigma^{2}) \tilde{{\bf D}}) + {\bf B}^{H}({\bf I} + \tilde{\beta}(\sigma^{2}) {\bf D} )^{-1} {\bf B} \right]^{-1}
\end{array}\right. \ .
\end{equation}
Matrices ${\bf T}(\sigma^{2})$ and $\tilde{{\bf T}}(\sigma^{2})$ satisfy 
\begin{equation}
\label{eq:borne-T}
{\bf T}(\sigma^{2}) \leq \frac{{\bf I}_r}{\sigma^{2}}, \qquad \tilde{{\bf T}}(\sigma^{2}) \leq \frac{{\bf I}_t}{\sigma^{2}}\ .
\end{equation}

\item[(ii)] ({\em Representation of the auxiliary quantities}) 
There exist two uniquely defined positive matrix-valued measures $\bs{\mu}$ and $\tilde{\bs{\mu}}$ such that 
$\bs{\mu}(\mathbb{R}^{+}) = {\bf I}_r$, $\tilde{\bs{\mu}}(\mathbb{R}^{+}) = {\bf I}_t$ and 
\begin{equation}
\label{eq:representationT}
{\bf T}(\sigma^{2}) =  \int_{\mathrm{R}^{+}} \frac{\bs{\mu}(d\lambda)}{\lambda + \sigma^{2}},\qquad 
\tilde{{\bf T}}(\sigma^{2}) =  \int_{\mathrm{R}^{+}} \frac{ \tilde{\bs{\mu}}(d \lambda)}{\lambda + \sigma^{2}} \ .
\end{equation}
The solutions $\beta(\sigma^{2})$ and  $\tilde{\beta}(\sigma^{2})$ of system \eqref{eq:canonique-preparatoire} 
are given by:
\begin{equation}
\label{eq:exprebeta}
\beta(\sigma^{2}) =  \frac{1}{t} \mathrm{Tr} {\bf D} {\bf T}(\sigma^{2})\ , \qquad 
\tilde{\beta}(\sigma^{2}) =  \frac{1}{t} \mathrm{Tr} \tilde{{\bf D}} \tilde{{\bf T}}(\sigma^{2})\ ,
\end{equation}
and can thus be written as 
\begin{equation}
\label{eq:represention-beta}
\beta(\sigma^{2}) =  \int_{\mathrm{R}^{+}} \frac{\mu_b(d \lambda)}{\lambda + \sigma^{2}} \ ,\qquad 
\tilde{\beta}(\sigma^{2}) =  \int_{\mathrm{R}^{+}} \frac{\tilde{\mu}_b(d \lambda) }{\lambda + \sigma^{2}} 
\end{equation}
where $ \mu_b$ and $\tilde{\mu}_b$ are nonnegative scalar measures
defined by 
$$ 
\mu_b(d \lambda) = \frac{1}{t} \mathrm{Tr}({\bf D} 
\bs{\mu}(d \lambda))\quad \textrm{and}\quad  \tilde{\mu}_b(d\lambda) = \frac{1}{t}
\mathrm{Tr}(\tilde{{\bf D}} \tilde{\bs{\mu}}(d\lambda)).
$$

\item[(iii)] ({\em Asymptotic approximation})
Assume that Assumption \ref{ass:normes2} holds and that
$$
\sup_{t} \|{\bf D}\|< d_{\max} < +\infty\quad \textrm{and}\quad \sup_{t}
\|\tilde{{\bf D}}\| < \tilde{d}_{\max} < +\infty\ .
$$ 
For every deterministic matrices ${\bf M}$ and $\tilde{{\bf M}}$ satisfying $\sup_{t} \|{\bf M}\| < +\infty$ 
and $\sup_{t} \|\tilde{{\bf M}}\| < +\infty$, the following limits hold true almost surely:
\begin{equation}
\label{eq:convergence-normalise}
\left\{
\begin{array}{ccc}
\lim_{t \rightarrow +\infty} \frac{1}{r} \mathrm{Tr} \left[ ({\bf S}(\sigma^{2}) - {\bf T}(\sigma^{2})) {\bf M} \right] &  = &  0\ \\
\lim_{t \rightarrow +\infty} \frac{1}{t} \mathrm{Tr} \left[ (\tilde{{\bf S}}(\sigma^{2}) - \tilde{{\bf T}}(\sigma^{2})) \tilde{{\bf M}} \right] &  = &  0
\end{array}\right. \ .
\end{equation}
Denote by $\mu$ and $\tilde{\mu}$ the (scalar) probability measures
$\mu = \frac{1}{r} \mathrm{Tr} \bs{\mu}$ and $\tilde{\mu} =
\frac{1}{t} \mathrm{Tr} \tilde{\bs{\mu}}$, by $(\lambda_i)$ (resp.
$(\tilde{\lambda}_j)$) the eigenvalues of $\bs{\Sigma}\bs{\Sigma}^{H}$
(resp. of $\bs{\Sigma}^{H}\bs{\Sigma}$). The following limits hold
true almost surely:
\begin{equation}
\label{eq:convergence-distribution}
\left\{
\begin{array}{ccc}
  \lim_{t \rightarrow +\infty} \frac{1}{r} \sum_{i=1}^{r} \phi(\lambda_i) - \int_{0}^{+\infty} \phi(\lambda) \; \mu(d \lambda) & = & 0 \\
  \lim_{t \rightarrow +\infty} \frac{1}{t} \sum_{j=1}^{t}\tilde \phi(\lambda_j) - \int_{0}^{+\infty} \tilde{\phi}(\lambda) \; \tilde{\mu}(d\lambda) & = & 0
\end{array} \right. \ ,
\end{equation}
for continuous bounded functions $\phi$ and $\tilde{\phi}$ defined on $\mathbb{R}^{+}$.

\end{itemize}
\end{theorem} 

The proof of $(i)$ is provided in Appendix
\ref{proof-existence-unicite} (note that in \cite{HLN07}, the
existence and uniqueness of solutions to the system
(\ref{eq:canonique-preparatoire}) is proved in a certain class of
analytic functions depending on $\sigma^{2}$ but this does not imply
the existence of a unique solution $(\beta,\tilde \beta)$ when
$\sigma^2$ is fixed). The rest of the statements of Theorem
\ref{theo:canonique-preparatoire} have been established in
\cite{HLN07}, and their proof is omitted here.

\begin{remark} As shown in \cite{HLN07}, the results in Theorem \ref{theo:canonique-preparatoire}
do not require any Gaussian assumption for  
$\bs{\Sigma}$. Remark that
(\ref{eq:convergence-normalise}) implies in some sense that the
entries of ${\bf S}(\sigma^{2})$ and $\tilde{{\bf S}}(\sigma^{2})$
have the same behaviour as the entries of the deterministic matrices
${\bf T}(\sigma^{2})$ and $\tilde{{\bf T}}(\sigma^{2})$  (which can be
evaluated by solving the system (\ref{eq:canonique-preparatoire})). In
particular, using (\ref{eq:convergence-normalise}) for ${\bf M} = {\bf
  I}$, it follows that the Stieltjes transform $s(\sigma^{2})$ of the
eigenvalue distribution of $\bs{\Sigma} \bs{\Sigma}^{H}$ behaves like
$\frac{1}{r} \mathrm{Tr} {\bf T}(\sigma^{2})$, which is itself the
Stieltjes transform of measure $\mu = \frac{1}{r} \mathrm{Tr}
\bs{\mu}$. The convergence statement
(\ref{eq:convergence-distribution}) which states that the eigenvalue
distribution of $\bs{\Sigma} \bs{\Sigma}^{H}$ (resp. $\bs{\Sigma}^{H}
\bs{\Sigma}$) has the same behavior as
$\mu$ (resp. $\tilde{\mu}$) directly follows from this observation. 
\end{remark}  
\subsubsection{The asymptotic approximation of the EMI}\label{subsub-approx}
Denote by $J(\sigma^{2}) = \mathbb{E} \log \det \left({\bf I}_r + \sigma^{-2}
  \bs{\Sigma} \bs{\Sigma}^{H} \right) $ the EMI associated with matrix ${\bf \Sigma}$.
First notice that 
$$
\log \det \left( {\bf I} + \frac{\bs{\Sigma}\bs{\Sigma}^{H}}{\sigma^{2}}\right)
= \sum_{i=1}^r \log \left( 1+\frac{\lambda_i}{\sigma^2}\right)\ ,
$$
where the $\lambda_i$'s stand for the eigenvalues of
$\bs{\Sigma}\bs{\Sigma}^{H}$. Applying
(\ref{eq:convergence-distribution}) to function $\phi(\lambda) = \log
(\lambda+\sigma^{2})$ (plus some extra work since $\phi$ is not
bounded), we obtain:
\begin{equation}
\label{eq:convergence-logdet-etape0}
\lim_{t \rightarrow +\infty} \left( \frac{1}{r} \log \det \left( {\bf I} + \frac{\bs{\Sigma}\bs{\Sigma}^{H}}{\sigma^{2}}\right)  -  \int_{0}^{+\infty} \log(\lambda + \sigma^{2}) \; d\mu(\lambda) \right) = 0\ .
\end{equation}
% It is possible to obtain a rather explicit expression of
% $\int_{0}^{+\infty} \log(\lambda + \sigma^{2}) \; d\mu(\lambda)$.
Using the well known relation: 
\begin{eqnarray}
\label{eq:lien-trace-log}
\frac{1}{r} \log \det \left( {\bf I} + \frac{\bs{\Sigma}\bs{\Sigma}^{H}}{\sigma^{2}}\right) &=&
\int_{\sigma^{2}}^{+\infty} \left( \frac{1}{\omega} - \frac{1}{r} \mathrm{Tr} (\bs{\Sigma} \bs{\Sigma}^{H} + \omega {\bf I})^{-1} \right) 
\; d\omega\nonumber \\
&=& \int_{\sigma^{2}}^{+\infty} \left( \frac{1}{\omega} - \frac{1}{r} \mathrm{Tr}\, {\bf S}(\omega)\right) 
\; d\omega\ ,
\end{eqnarray}
together with the fact that ${\bs S}(\omega) \approx {\bs T}(\omega)$
(which follows from Theorem \ref{theo:canonique-preparatoire}), it is
proved in \cite{HLN07} that:
\begin{equation}
\label{eq:convergence-logdet}
\lim_{t \rightarrow +\infty} \left[\frac{1}{r} \log \det \left( {\bf I} + \frac{\bs{\Sigma}\bs{\Sigma}^{H}}{\sigma^{2}}\right)  - 
  \int_{\sigma^{2}}^{+\infty} \left( \frac{1}{\omega} - \frac{1}{r} \mathrm{Tr} {\bf T}(\omega) \right) \; d\omega \right]=0 
\end{equation}
almost surely. Define by $\bar{J}(\sigma^{2})$ the quantity:
\begin{equation}
\label{eq:defJbarre}
\bar{J}(\sigma^{2}) 
= r \int_{\sigma^{2}}^{+\infty} \left( \frac{1}{\omega} - \frac{1}{r} \mathrm{Tr} {\bf T}(\omega) \right) \; d\omega\ .
\end{equation}
Then, $\bar{J}(\sigma^{2})$ can be expressed more explicitely as:
\begin{multline}
\label{eq:expreJbarre}
\bar{J}(\sigma^{2})
= \log \det \left[{\bf I}_r + \tilde{\beta}(\sigma^{2}) {\bf D} + \frac{1}{\sigma^{2}} {\bf B} ({\bf I}_t + \beta(\sigma^{2}) \tilde{{\bf D}})^{-1} {\bf B}^{H} \right] \\
+ \log \det \left[{\bf I}_t + \beta(\sigma^{2}) \tilde{{\bf D}}
\right] - \sigma^{2} t \beta(\sigma^{2}) \tilde{\beta}(\sigma^{2})\ ,
\end{multline}
or equivalently as
\begin{multline}
\label{eq:expreJbarrebis}
\bar{J}(\sigma^{2}) = \log \det \left[{\bf I}_t + \beta(\sigma^{2}) \tilde{{\bf D}} + \frac{1}{\sigma^{2}} {\bf B}^{H} ({\bf I}_r + \tilde{\beta}(\sigma^{2}) {\bf D})^{-1} {\bf B} \right] \\
+  \log \det \left[{\bf I}_r + \tilde{\beta}(\sigma^{2}) {\bf D} \right] - \sigma^{2} t \beta(\sigma^{2}) \tilde{\beta}(\sigma^{2})\ .
\end{multline}
Taking the expectation with respect to the channel $\bs{\Sigma}$ in
(\ref{eq:convergence-logdet}), the EMI $J(\sigma^{2}) = \mathbb{E} \log \det \left({\bf I}_r + \sigma^{-2}
  \bs{\Sigma} \bs{\Sigma}^{H} \right) $ can 
be approximated by $\bar{J}(\sigma^{2})$:
\begin{equation}
\label{eq:vitesse-minable}
J(\sigma^{2}) = \bar{J}(\sigma^{2}) + o(t)
\end{equation}
as $t \rightarrow +\infty$. This result is fully proved in
\cite{HLN07} and is of potential interest since the numerical
evaluation of $\bar{J}(\sigma^{2})$ only requires to solve the
$2\times 2$ system (\ref{eq:canonique-preparatoire}) while the
calculation of $J(\sigma^{2})$ either rely on Monte-Carlo simulations
or on the implementation of rather complicated explicit formulas (see
for instance \cite{Kang-Alouini-Rice-06}). 

In order to evaluate the precision of the asymptotic approximation
$\bar{J}$, we shall improve \eqref{eq:vitesse-minable} and get the
speed $ J(\sigma^{2}) = \bar{J}(\sigma^{2}) + O(t^{-1}) $ in the next
theorem. This result completes those in \cite{HLN07} and on the
contrary of the rest of Theorem \ref{theo:canonique-preparatoire}
heavily relies on the Gaussian structure of ${\bf \Sigma}$. We
first introduce very mild extra assumptions:
\begin{ass}
\label{ass:norme-spectrale} Let $({\bf B}_t)$ be a family of $r\times t$ deterministic matrices such that 
$$
\sup_{t} \| {\bf B} \| < b_{\max} < + \infty\ .
$$
\end{ass}

\begin{ass}
\label{ass-borneinf}
Let ${\bf D}$ and $\tilde{{\bf D}}$ be respectively $r\times r$ and $t\times t$ diagonal matrices such that 
$$
\sup_{t} \|{\bf D}\|< d_{\max} < +\infty\quad \textrm{and} \quad 
\sup_{t} \|\tilde{{\bf D}}\| < \tilde{d}_{\max} < +\infty\ .
$$
\indent Assume moreover that 
$$
\inf_{t} \frac{1}{t} \mathrm{Tr} {\bf D} > 0 \quad \textrm{and}\quad
\inf_{t} \frac{1}{t} \mathrm{Tr} \tilde{{\bf D}} > 0\ .
$$
\end{ass}
\begin{theorem}
\label{theo:convergence-rate}
Recall that ${\bf \Sigma}= {\bf B} +{\bf Y}$ and assume that ${\bf Y}
= \frac 1{\sqrt{t}}{\bf D}^{\frac 12} {\bf X}\, \tilde{{\bf D}}^{\frac
  12}$, where ${\bf D}=\mbox{diag}(d_i)$ and $\tilde{{\bf
    D}}=\mbox{diag}(\tilde{d}_j)$ are $r \times r$ and $t\times t$
diagonal matrices and where ${\bf X}$ is a matrix whose entries are
i.i.d. complex circular Gaussian variables ${\mathcal CN}(0,1)$.
Assume moreover that Assumptions \ref{ass:norme-spectrale} and \ref{ass-borneinf} hold true.
Then, for every deterministic matrices ${\bf M}$ and $\tilde{{\bf M}}$ satisfying 
$\sup_{t} \|{\bf M}\| < +\infty$ and $\sup_{t} \|\tilde{{\bf M}}\| < +\infty$, 
the following facts hold true: 
\begin{equation}
\label{eq:vitesse-convergence-normalise}
\mathrm{Var} \, \left( \frac{1}{r} \mathrm{Tr} \left[ {\bf S}(\sigma^{2}) {\bf M} \right] \right)   =  O\left(\frac{1}{t^{2}}\right) \quad \textrm{and}\quad 
\mathrm{Var} \, \left(  \frac{1}{t} \mathrm{Tr} \left[ \tilde{{\bf S}}(\sigma^{2}) \tilde{{\bf M}} \right]  \right)  =  O\left(\frac{1}{t^{2}}\right)
\end{equation}
where $\mathrm{Var}(.)$ stands for the variance. Moreover, 
\begin{equation}
\label{eq:vitesse-convergence-biais}
\begin{array}{ccc}
\frac{1}{r} \mathrm{Tr} \left[ (\mathbb{E}({\bf S}(\sigma^{2})) - {\bf T}(\sigma^{2})) {\bf M} \right]  
&=&  O\left(\frac{1}{t^{2}}\right)\\
\frac{1}{t} \mathrm{Tr} \left[ (\mathbb{E}(\tilde{{\bf S}}(\sigma^{2})) - \tilde{{\bf T}}(\sigma^{2})) \tilde{{\bf M}} \right]   &=& O\left(\frac{1}{t^{2}}\right)
\end{array}
\end{equation}
and  
\begin{equation}
\label{eq:vitesse-rapide}
J(\sigma^{2}) = \bar{J}(\sigma^{2}) + O\left(\frac{1}{t}\right)\ .
\end{equation}
\end{theorem}
The proof is given in Appendix \ref{proof-theo-convergence-rate}. We provide here some comments.
\begin{remark}
The proof of Theorem \ref{theo:convergence-rate} takes full advantage of the Gaussian
structure of matrix $\bs{\Sigma}$
and relies on two simple ingredients:
\begin{itemize}
\item[(i)] An integration by parts formula that provides an expression
  for the expectation of certain functionals of Gaussian vectors,
  already well-known and widely used in Random Matrix Theory
  \cite{Kho-Pas-93,Pas-Kho-Vas-95}.
\item[(ii)] An inequality known as Poincar\'e-Nash inequality that
  bounds the variance of functionals of Gaussian vectors. Although
  well known, its application to random matrices is fairly recent
  (\cite{Chatterjee-Bose-2004}, \cite{Pastur-06}, see also \cite{HKLNP06pre}).
\end{itemize}
\end{remark}

\begin{remark}
Equations (\ref{eq:vitesse-convergence-normalise}) also hold
in the non Gaussian case and can be established by using the so-called REFORM 
(Resolvent FORmula Martingale) method introduced by Girko (\cite{Girko-01}).

Equations (\ref{eq:vitesse-convergence-biais}) and
(\ref{eq:vitesse-rapide}) are specific to the complex Gaussian
structure of the channel matrix $\bs{\Sigma}$. In particular, in the
non Gaussian case, or in the real Gaussian case, one would get $J(\sigma^{2}) = \bar{J}(\sigma^{2}) +
O(1)$. These two facts are in accordance with:
\begin{itemize}
\item[(i)] The work of \cite{BaiSil04} in which a weaker
  result ($o(1)$ instead of $O(t^{-1})$) is proved in the simpler case
  where ${\bf B} = {\bf 0}$;
\item[(ii)] The predictions of the replica method in
  \cite{Moustakas-Simon-Sengupta-03} (resp. \cite{Moustakas-Simon-05})
  in the case where ${\bf B} = \mb{0}$ (resp. in the case where $\mb{\tilde{D}} = \mb{I}_t$ and
${\bf D} = {\bf I}_r$); 
\end{itemize}
\end{remark}

\begin{remark} [Standard deviation and bias] Eq.
  (\ref{eq:vitesse-convergence-normalise}) implies that the standard
  deviation of $ \frac{1}{r} \mathrm{Tr} \left[ ({\bf S}(\sigma^{2}) -
    {\bf T}(\sigma^{2})) {\bf M} \right]$ and $\frac{1}{t} \mathrm{Tr}
  \left[ (\tilde{{\bf S}}(\sigma^{2}) - \tilde{{\bf T}}(\sigma^{2}))
    \tilde{{\bf M}} \right]$ are of order $O(t^{-1})$ terms. However,
  their mathematical expectations (which correspond to the bias)
  converge much faster towards $0$ as
  (\ref{eq:vitesse-convergence-biais}) shows (the order is $O(t^{-2})$).
\end{remark}

\begin{remark}
  By adapting the techniques developed in the course of the proof of
  Theorem \ref{theo:convergence-rate}, one may establish that 
${\bf u}^{H} \mathbb{E}{\bf S}(\sigma^{2}) {\bf v} - {\bf u}^{H} {\bf T}(\sigma^{2}) 
{\bf v} = O\left(\frac{1}{t}\right)\ ,$
where ${\bf u}$ and ${\bf v}$ are uniformly bounded $r$-dimensional vectors.
%
% ${\bf u}^{H} {\bf S}(\sigma^{2}) {\bf v} - {\bf u}^{H} {\bf
%   T}(\sigma^{2}) {\bf v}$ is a $O_p(\frac{1}{\sqrt{t}})$ term, and
% that its mathematical expectation is a $O(\frac{1}{t})$ term.
\end{remark}

\begin{remark}
  Both $J(\sigma^{2})$ and $ \bar{J}(\sigma^{2})$ increase
  linearly with $t$. Equation (\ref{eq:vitesse-rapide}) thus
  implies that the relative error $\frac{J(\sigma^{2}) -
    \bar{J}(\sigma^{2})}{J(\sigma^{2})}$ is of order $O(t^{-2})$.
  This remarkable convergence rate strongly supports the observed fact
  that approximations of the EMI remain reliable even for small numbers
  of antennas (see also the numerical results in section
  \ref{sec:simulations}).  Note that similar observations have been
  done in other contexts where random matrices are used, see e.g. 
  \cite{Biglieri-Taricco-Tulino-02}, \cite{Moustakas-Simon-Sengupta-03}. 
\end{remark}

\subsection{Introduction of the virtual channel ${\bf H Q}^{\frac 12}$}\label{sec:virtual}
The purpose of this section is to establish a link between the
simplified model \eqref{eq:modele-equivalent-iid}: ${\bf \Sigma} ={\bf
  B} + {\bf Y}$ where ${\bf Y} = \frac 1{\sqrt{t}} {\bf D}^{\frac 12}
{\bf X} \tilde {\bf D}^{\frac 12}$, ${\bf X}$ being a matrix with
i.i.d ${\mathcal C} N(0,1)$ entries, ${\bf D}$ and $\tilde{\bf D}$
being diagonal matrices, and the Rician model \eqref{eq:modeleH} under
investigation: $\mb{H} = \sqrt{\frac{K}{K+1}} {\bf A} +
\frac{1}{\sqrt{K+1}} {\bf V}$ where ${\bf V} = \frac 1{\sqrt{t}} {\bf
  C}^{\frac 12} {\bf W} \tilde {\bf C}^{\frac 12}$. As we shall see,
the key point is the unitary invariance of the EMI of Gaussian
channels together with a well-chosen eingenvalue/eigenvector
decomposition.

We introduce the virtual channel ${\bf H} {\bf Q}^{\demi}$ which can be written as:
\begin{equation}
  \label{eq:virtual-representation}
  {\bf H} {\bf Q}^{\frac 12} = \sqrt{\frac{K}{K+1}} {\bf A}  {\bf Q}^{\frac 12} + \frac{1}{\sqrt{K+1}} {\bf C}^{\frac 12} \frac{{\bf W}}{\sqrt{t}} \bs{\Theta}
  ({\bf Q}^{\frac 12} \tilde{{\bf C}} {\bf Q}^{\frac 12})^{\frac 12}\ ,
\end{equation}
where $\bs{\Theta}$ is the deterministic unitary $t \times t$ matrix defined by 
\begin{equation}\label{eq:theta}
\bs{\Theta} = \tilde{{\bf C}}^{\frac 12} {\bf Q}^{\frac 12} ({\bf
  Q}^{\frac 12} \tilde{{\bf C}} {\bf Q}^{\frac 12})^{-\frac 12}\ .
\end{equation}
The virtual channel ${\bf H} {\bf Q}^{\frac 12}$ has thus a structure similar to ${\bf H}$, where 
$({\bf A}, {\bf C}, \tilde{{\bf C}}, {\bf W})$ are respectively replaced with 
$({\bf A} {\bf Q}^{\frac 12}, {\bf C}, {\bf Q}^{\frac 12} \tilde{{\bf C}} {\bf Q}^{\frac 12}, {\bf W} \bs{\Theta})$.\\

Consider now the eigenvalue/eigenvector decompositions of matrices
$\frac{{\bf C}}{\sqrt{K+1}}$ and $\frac{{\bf Q}^{\frac 12} \tilde{{\bf
      C}} {\bf Q}^{\frac 12}}{\sqrt{K+1}}$:
\begin{equation}\label{eq:eigen}
\frac{{\bf    C}}{\sqrt{K+1}} = {\bf U} {\bf D} {\bf U}^{H}\qquad \textrm{and}\qquad \frac{{\bf
    Q}^{\frac 12} \tilde{{\bf C}} {\bf Q}^{\frac 12}}{\sqrt{K+1}} =
\tilde{{\bf U}} \tilde{{\bf D}} \tilde{{\bf U}}^{H}\ .
\end{equation}
Matrices ${\bf U}$ and $\tilde{{\bf U}}$ are the eigenvectors matrices
while ${\bf D}$ and $\tilde{{\bf D}}$ are the eigenvalues diagonal
matrices.  It is then clear that the ergodic mutual information of
channel ${\bf H} {\bf Q}^{\frac 12}$ coincides with the EMI of
$\bs{\Sigma} = {\bf U}^{H} {\bf H} {\bf Q}^{1/2} \tilde{{\bf U}}$.
Matrix ${\bf \Sigma}$ can be written as $\bs{\Sigma} = {\bf B} + {\bf
  Y}$ where 
\begin{equation}\label{eq:definition-B}
{\bf B} = \sqrt{\frac K{K+1}}{\bf U}^{H} {\bf A} {\bf
  Q}^{\frac 12} \tilde{{\bf U}}\qquad \textrm{and}\qquad  {\bf Y} = \frac 1{\sqrt{t}} {\bf
  D}^{\frac 12} {\bf X} \tilde{{\bf D}}^{\frac 12}\quad \textrm{with}\quad  {\bf X}={\bf U}^{H}
{\bf W} \bs{\Theta} \tilde{{\bf U}}\ .
\end{equation}
As matrix ${\bf W}$ has i.i.d. ${\mathcal C}N(0,1)$ entries, so has
matrix ${\bf X}= {\bf U}^{H} {\bf W} \bs{\Theta} \tilde{{\bf U}}$ due
to the unitary invariance. Note that the entries of ${\bf Y}$ are
independent since ${\bf D}$ and $\tilde{{\bf D}}$ are diagonal.
We sum up the previous discussion in the following proposition.

\begin{proposition}\label{prop:virtual}
  Let ${\bf W}$ be a $r\times t$ matrix whose individual entries are
  i.i.d. ${\mathcal C}N(0,1)$ random variables. The two ergodic
  mutual informations
$$
I({\bf Q}) =\mathbb{E} \log \det \left( {\bf I} +\frac{ {\bf H Q H}^H}{\sigma^2} \right)
\qquad \textrm{and}\qquad 
J(\sigma^2) = \mathbb{E} \log \det \left( {\bf I} +\frac{ {\bf \Sigma \Sigma}^H}{\sigma^2} \right)
$$
are equal provided that channel ${\bf H}$ is given by:
 $$\mb{H} = \sqrt{\frac{K}{K+1}} {\bf A} +
\frac{1}{\sqrt{K+1}} {\bf V}
$$ with ${\bf V}  =  \frac 1{\sqrt{t}} {\bf
  C}^{\frac 12} {\bf W} \tilde {\bf C}^{\frac 12}$; channel ${\bf \Sigma}$ by  
$
{\bf \Sigma} = {\bf B} + {\bf Y}
$ with ${\bf Y}  =  \frac 1{\sqrt{t}} {\bf  D}^{\frac 12} {\bf X} \tilde{{\bf D}}^{\frac 12}$ 
and that \eqref{eq:theta}, \eqref{eq:eigen} and \eqref{eq:definition-B} hold true.  
\end{proposition}

\subsection{Study of the EMI $I({\bf Q})$.}

We now apply the previous results to the study of the EMI of channel
${\bf H}$. We first state the corresponding result.
\begin{theorem}
\label{theo:canonique}
For ${\bf Q} \in {\cal C}_1$, consider the system of equations
\begin{equation}
\label{eq:canonique}
\left\{
\begin{array}{ccc}
\delta & = & f(\delta, \tilde{\delta}, \mb{Q})  \\
\tilde{\delta} & = & \tilde{f}(\delta, \tilde{\delta}, \mb{Q})
\end{array}\right.\ ,
\end{equation}
where $f(\delta, \tilde{\delta}, \mb{Q})$ and
$\tilde{f}(\delta, \tilde{\delta}, \mb{Q})$ are given by: 
\begin{multline}
\label{eq:expref}
f(\delta, \tilde{\delta}, \mb{Q}) =   \f{1}{t}\,  \mathrm{Tr} \bigg\{ {\bf C} \Big[ \sigma^2\,\big({\bf
I}_r+  \frac{\tilde\delta}{K+1} \, {\bf C}\,\big)\\
+ \frac{K}{K+1} {\bf A}{\bf Q}^{\demi}\pa{{\bf
I}_t+ \frac{\delta}{K+1} \,{\bf Q}^{\demi} {\bf \tilde{C}}{\bf Q}^{\demi}}^{-1}{\bf Q}^{\demi}{\bf A}^H \Big]^{-1}\bigg\} \ , 
\end{multline}
\begin{multline}
\tilde{f}(\delta, \tilde{\delta}, \mb{Q}) =  \f{1}{t}\, \mathrm{Tr}
\bigg\{ {\bf Q}^{\demi} {\bf \tilde{C}} {\bf Q}^{\demi} \Big[\sigma^2\,\big({\bf I}_t+\frac{\delta}{K+1}\, {\bf Q}^{\demi} \tilde{{\bf
C}}{\bf Q}^{\demi} \,\big)\\
+ \frac{K}{K+1} {\bf Q}^{\demi}{\bf A}^H\pa{{\bf
I}_r+\frac{\tilde\delta}{K+1}\,{\bf C}}^{-1}{\bf A}{\bf Q}^{\demi}\Big]^{-1}\bigg\} \ .\label{eq:expreftildef}
\end{multline}
Then the system of equations (\ref{eq:canonique}) has a unique
strictly positive solution $(\delta(\mb{Q}), \tilde{\delta}(\mb{Q}))$. \\
Furthermore, assume that $\sup_{t} \| \mb{Q} \| < +\infty$, 
$\sup_{t} \|{\bf A}\| < +\infty$,  $\sup_{t} \|{\bf C}\| < +\infty$, 
and $\sup_{t} \|\tilde{{\bf C}}\| < +\infty$. 
Assume also that $\inf_{t} \lambda_{\min}(\tilde{{\bf C}}) > 0$ where 
$\lambda_{\min}(\tilde{{\bf C}})$ represents the smallest eigenvalue of $\tilde{{\bf C}}$. 
Then, as $t \rightarrow +\infty$,
\begin{equation}
\label{eq:equivalent1} I(\mb{Q}) = \bar{I}(\mb{Q}) + O\left(\frac{1}{t}\right)
\end{equation}
where the asymptotic approximation $\bar{I}(\mb{Q})$ is given by
\begin{multline}
\label{eq:expreCbarre}
\bar{I}(\mb{Q}) = \log{}\det \left( \, {\bf I}_t+
  \frac{\delta(\mb{Q})}{K+1} \, {\bf Q}^{\demi} \tilde{{\bf C}} {\bf
    Q}^{\demi} + \f{1}{\sigma^2} \, \frac{K}{K+1} \,
  {\bf Q}^{\demi} {\bf A}^H \left( {\bf I}_r+\frac{\tilde{\delta}(\mb{Q})}{K+1} \,{\bf C} \right)^{-1}{\bf A}{\bf Q}^{\demi}\,\right)\\
+ \log \det \left({\bf I}_r+ \frac{\tilde \delta(\mb{Q})}{K+1} {\bf
    C}\right) - \frac{ t \sigma^2}{K+1}
\delta(\mb{Q})\,\tilde{\delta}(\mb{Q})\ ,
\end{multline}
%%%%%%%%%%%%%%%%%%%%%%%%%%%%%%%%%%%%%%%%%%%%%%%%%%%%%%%%%%%%%%%%%%%%%%%%%%%%%%%%%%%%%%%%%%%%%
or equivalently by
\begin{multline}
\label{eq:expreCbarrebis}
\bar{I}(\mb{Q}) =  \log \det\left( {\bf I}_r+ \frac{\tilde\delta(\mb{Q})}{K+1}
{\bf C}+\f{1}{\sigma^2}\,\frac{K}{K+1} \, {\bf A}{\bf Q}^{\demi} \left( {\bf
I}_t+\frac{\delta(\mb{Q})}{K+1} \,{\bf Q}^{\demi} {\bf \tilde{C}} {\bf Q}^{\demi} \right)^{-1} {\bf Q}^{\demi}{\bf A}^H\,\right)\\
+ \log \det \left({\bf I}_t+ \frac{\delta(\mb{Q})}{K+1} {\bf Q}^{1/2} {\bf \tilde{C}} {\bf Q}^{1/2} \,\right)
-\frac{t \sigma^2}{K+1}  \delta(\mb{Q})\,\tilde \delta(\mb{Q}).
\end{multline}
\end{theorem}
%%%%%%%%%%%%%

\medskip
%%%%%%%%%%%%%%%%%%%%%%%%%%%%%%%%%%%%%%%%%%%%%%%%%%%%%%%%%%%%%%%%%%%%%%%%%%%%%%%%%%%%%%%%%%%%%

\begin{proof} We rely on the virtual channel introduced in Section
  \ref{sec:virtual} and on the eigenvalue/eigenvector decomposition
  performed there. 

  Matrices ${\bf B}$, ${\bf D}$, $\tilde{{\bf D}}$ as introduced in
  Proposition \ref{prop:virtual} are clearly uniformly bounded, while
  $\inf_{t} \frac{1}{t} \mathrm{Tr} {\bf D} = \inf_{t} \frac{1}{t}
  \mathrm{Tr} {\bf C} =1$ due to the model specifications and $\inf_t
  \frac{1}{t} \mathrm{Tr} {\bf Q}^{\frac 12} \tilde{{\bf C}} {\bf
    Q}^{\frac 12} \geq \inf_t \lambda_{\min}(\tilde{{\bf C}})
  \frac{1}{t} \mathrm{Tr} {\bf Q} > 0$ as $\frac{1}{t} \mathrm{Tr}
  {\bf Q} = 1$. Therefore, matrices ${\bf B}$, ${\bf D}$ and
  $\tilde{{\bf D}}$ clearly satisfy the assumptions of Theorems
  \ref{theo:canonique-preparatoire} and \ref{theo:convergence-rate}.

We first apply the results of Theorem
\ref{theo:canonique-preparatoire} to matrix $\bs{\Sigma}$, and use the
same notations as in the statement of Theorem
\ref{theo:canonique-preparatoire}. Using the unitary invariance of the
trace of a matrix, it is straightforward to check that:
\begin{eqnarray*}
  \frac{f(\delta, \tilde{\delta}, {\bf Q})}{\sqrt{K+1}} 
& = & \frac{1}{t} \textrm{Tr} \left[ {\bf D} \left( \sigma^{2}\left({\bf I}+ {\bf D} \frac{\tilde{\delta}}{\sqrt{K+1}}\right) 
+ {\bf B} \left({\bf I}+\tilde{{\bf D}}\frac{\delta}{\sqrt{K+1}}\right)^{-1} {\bf B}^{H} \right)^{-1} \right]\ , \\
  \frac{\tilde{f}(\delta, \tilde{\delta}, {\bf Q})}{\sqrt{K+1}} 
& = & \frac{1}{t} \textrm{Tr} \left[ \tilde{{\bf D}} \left( \sigma^{2}\left({\bf I}+ \tilde{{\bf D}} \frac{\delta}{\sqrt{K+1}}\right) + {\bf B}^{H} \left({\bf I}+{\bf D} \frac{\tilde{\delta}}{\sqrt{K+1}}\right)^{-1} {\bf B} \right)^{-1} \right] \ .
\end{eqnarray*}
Therefore, $(\delta, \tilde{\delta})$ is solution of
(\ref{eq:canonique}) if and only if $(\frac{\delta}{\sqrt{K+1}},
\frac{\tilde{\delta}}{\sqrt{K+1}})$ is solution of
(\ref{eq:canonique-preparatoire}). As the system (\ref{eq:canonique-preparatoire})
admits a unique solution, say $(\beta, \tilde{\beta})$, the solution $(\delta, \tilde{\delta})$
to (\ref{eq:canonique}) exists, is unique and is related to $(\beta, \tilde{\beta})$ by the relations:
\begin{equation}
\label{eq:lien-delta-beta}
\beta = \frac{\delta}{\sqrt{K+1}}, \; \; \tilde{\beta} = \frac{\tilde{\delta}}{\sqrt{K+1}}\ .
\end{equation}
In order to justify (\ref{eq:expreCbarre}) and
(\ref{eq:expreCbarrebis}), we note that $J(\sigma^{2})$ coincides with
the EMI $I({\bf Q})$. Moreover, the unitary invariance of the
determinant of a matrix together with (\ref{eq:lien-delta-beta}) imply
that $\bar{I}({\bf Q})$ defined by (\ref{eq:expreCbarre}) and
(\ref{eq:expreCbarrebis}) coincide with the approximation $\bar{J}$
given by (\ref{eq:expreJbarre}) and (\ref{eq:expreJbarrebis}). This
proves (\ref{eq:equivalent1}) as well. 
\end{proof}

In the following, we denote by ${\bf T}_K(\sigma^{2})$ and
$\tilde{{\bf T}}_K(\sigma^{2})$ the following matrix-valued functions:
\begin{equation}
\label{eq:defTK}
\left\{
\begin{array}{ccc}
{\bf T}_K(\sigma^{2}) & = & \left[ \sigma^{2}({\bf I} +  \frac{\tilde{\delta}}{K+1} {\bf C}) + \frac{K}{K+1} {\bf A} {\bf Q}^{\frac 12} ({\bf I} + \frac{\delta}{K+1} 
{\bf Q}^{\frac 12} \tilde{{\bf C}} {\bf Q}^{\frac 12})^{-1}  {\bf Q}^{\frac 12} {\bf A}^{H} \right]^{-1} \\
\tilde{{\bf T}}_K(\sigma^{2}) & = & \left[ \sigma^{2}({\bf I} + \frac{\delta}{K+1} {\bf Q}^{\frac 12} \tilde{{\bf C}} {\bf Q}^{\frac 12} ) + \frac{K}{K+1} {\bf Q}^{\frac 12} {\bf A}^{H}({\bf I} + \frac{\tilde{\delta}}{K+1} {\bf C} )^{-1} {\bf A} {\bf Q}^{\frac 12} \right]^{-1} 
\end{array}\right. \ .
\end{equation}
They are related to matrices ${\bf T}$ and $\tilde{{\bf T}}$ defined by (\ref{eq:defT}) by the relations:
\begin{equation}\label{equiv-tk}
\left\{
\begin{array}{ccc}
{\bf T}_K(\sigma^{2}) & = & {\bf U} {\bf T}(\sigma^{2}) {\bf U}^{H} \\
\tilde{{\bf T}}_K(\sigma^{2}) & = & \tilde{{\bf U}} \tilde{{\bf T}}(\sigma^{2}) \tilde{{\bf U}}^{H}
\end{array} \right.\ ,
\end{equation}
and their entries represent deterministic approximations of $({\bf H} {\bf Q} {\bf H}^{H} +
\sigma^{2} {\bf I}_r)^{-1}$ and $({\bf Q}^{\frac 12}{\bf H}^{H} {\bf
  H}{\bf Q}^{\frac 12} + \sigma^{2} {\bf I}_t)^{-1}$ (in the sense of 
  Theorem \ref{theo:canonique-preparatoire}).

  As $\frac{1}{r} \mathrm{Tr} {\bf T}_K = \frac{1}{r} \mathrm{Tr} {\bf
    T}$ and $\frac{1}{t} \mathrm{Tr} \tilde{{\bf T}}_K = \frac{1}{t}
  \mathrm{Tr} \tilde{{\bf T}}$, the quantities $\frac{1}{r}
  \mathrm{Tr} {\bf T}_K$ and $\frac{1}{t} \mathrm{Tr} \tilde{{\bf
      T}}_K$ are the Stieltjes transforms of probability measures
  $\mu$ and $\tilde{\mu}$ introduced in Theorem
  \ref{theo:canonique-preparatoire}. As matrices ${\bf H} {\bf Q} {\bf
    H}^{H}$ and $\bs{\Sigma} \bs{\Sigma}^{H}$ (resp. ${\bf Q}^{\frac
    12}{\bf H}^{H} {\bf H}{\bf Q}^{\frac 12}$ and $\bs{\Sigma}^{H}
  \bs{\Sigma}$) have the same eigenvalues,
  (\ref{eq:convergence-distribution}) implies that the eigenvalue
  distribution of ${\bf H} {\bf Q} {\bf H}^{H}$ (resp. ${\bf Q}^{\frac
    12}{\bf H}^{H} {\bf H}{\bf Q}^{\frac 12}$) behaves like $\mu$
  (resp. $\tilde{\mu}$).

We finally mention that $\delta(\sigma^{2})$ and  $\tilde{\delta}(\sigma^{2})$ are given by
\begin{equation}
\label{eq:expredelta}
\delta(\sigma^{2})  =  \frac{1}{t} \mathrm{Tr} {\bf C} {\bf T}_K(\sigma^{2}) \qquad \textrm{and} \qquad
\tilde{\delta}(\sigma^{2}) =  \frac{1}{t} \mathrm{Tr} {\bf Q}^{\frac 12} \tilde{{\bf C}} {\bf Q}^{1/2} \tilde{{\bf T}}_K(\sigma^{2})\ ,
\end{equation}
and that the following representations hold true:
\begin{equation}
\label{eq:representation-delta}
\delta(\sigma^{2}) =  \int_{\mathrm{R}^{+}} \frac{\mu_d(d\,\lambda)}{\lambda + \sigma^{2}} \qquad \textrm{and}\qquad 
\tilde{\delta}(\sigma^{2}) =  \int_{\mathrm{R}^{+}} \frac{\tilde{\mu}_d(d\,\lambda) }{\lambda + \sigma^{2}} 
\ ,
\end{equation}
where $\mu_d$ and $\tilde{\mu}_d$ are positive measures on
$\mathbb{R}^{+}$ satisfying $\mu_d(\mathbb{R}^{+}) = \frac{1}{t}
\mathrm{Tr} {\bf C}$ and $\tilde{\mu}_d(\mathbb{R}^{+}) = \frac{1}{t}
\mathrm{Tr} {\bf Q}^{1/2} \tilde{{\bf C}} {\bf Q}^{1/2}$.

\section{Strict concavity of $\bar{I}({\bf Q})$ and approximation of the capacity $I({\bf Q}_*)$}
\label{sec:concavity}
\subsection{Strict concavity of $\bar{I}({\bf Q})$} The strict concavity of $\bar{I}({\bf Q})$ is an important issue for
optimization purposes (see Section \ref{sec:algo}). The main result of the section is the following: 
\begin{theorem}\label{theo:con}
The function ${\bf Q}\mapsto \bar{I}({\bf Q})$ is strictly concave on ${\mathcal C}_1$.
\end{theorem}

As we shall see, the concavity of $\bar{I}$ can be
established quite easily by relying on the concavity of the EMI
$I({\bf Q})=\mathbb{E} \log \det \left( {\bf I} + \frac{{\bf H Q
    H}^H}{\sigma^2}\right)$. The strict concavity is more demanding and
its proof is mainly postponed to Appendix \ref{app:constric}.  

Recall that we denote by ${\mathcal C_1}$ the set of nonnegative Hermitian 
$t\times t$ matrices whose normalized trace is equal to one (i.e. $t^{-1} \mathrm{Tr}\, {\bf Q}=1$). 
In the sequel, we shall rely on the following straightforward but
useful result:
\begin{proposition}
\label{prop:prop-concavite} Let $f:{\cal C}_1\rightarrow \mathbb{R}$ be a real function.
Then $f$ is strictly concave if and only if for every matrices ${\bf
  Q}_1, {\bf Q}_2$ (${\bf Q}_1\neq {\bf Q}_2$) of ${\cal C}_1$,
the function $\phi(\lambda)$ defined on $[0,1]$ by
\[
\phi(\lambda) = f \left( \lambda {\bf Q}_1 + (1 - \lambda) {\bf Q}_2 \right)
\]
is strictly concave.
\end{proposition}

\subsubsection{Concavity of the EMI}\label{concavity-EMI} We first recall that $ I({\bf Q})=
\mathbb{E} \log \det \left( {\bf I} +\frac{{\bf H} {\bf Q}
    {\bf H}^H}{\sigma^2}\right) $ is concave on ${\mathcal C}_1$, 
and provide a proof for the sake of completeness. 
Denote by
${\bf Q}=\lambda {\bf Q}_1 + (1-\lambda) {\bf Q}_2$ and let 
$\phi(\lambda)=I(\lambda {\bf Q}_1 + (1-\lambda) {\bf Q}_2)$.
Following Proposition \ref{prop:prop-concavite}, it is sufficient to
prove that $\phi$ is concave. As $\log \det \left( {\bf I} +\frac{{\bf HQ
    H}^H}{\sigma^2}\right) = \log \det \left( {\bf I} +\frac{{\bf H}^H  {\bf HQ}
    }{\sigma^2}\right)$, we have:
\begin{eqnarray*}
  \phi(\lambda) & = &\mathbb{E} \log \det \left( {\bf I} +\frac{{\bf HQ
      H}^H}{\sigma^2}\right)\ , \\
  \phi'(\lambda) & = & \mathbb{E}\, \mathrm{Tr} \left( {\bf I} +\frac{{\bf H}^H {\bf H Q}}{\sigma^2} \right)^{-1} 
\frac{{\bf H}^H {\bf H}}{\sigma^2} ({\bf Q_1 -Q_2})\ ,\\
\phi''(\lambda) & = & - \mathbb{E}\, \mathrm{Tr} \left[ \left( {\bf I} +\frac{{\bf H}^H {\bf H Q}}{\sigma^2} \right)^{-1}
 \frac{{\bf H}^H {\bf H}}{\sigma^2} ({\bf Q_1 -Q_2})\left( {\bf I} +\frac{{\bf H}^H {\bf H Q}}{\sigma^2} \right)^{-1}
 \frac{{\bf H}^H {\bf H}}{\sigma^2} ({\bf Q_1 -Q_2})
\right]  \ .
\end{eqnarray*} 
In order to conclude that $\phi''(\lambda)\le 0$, we notice
that $\left( {\bf I} +\frac{{\bf H}^H {\bf H Q}}{\sigma^2}
\right)^{-1} \frac{{\bf H}^H {\bf H}}{\sigma^2}$ coincides with 
\[
{\bf H}^{H} \left( {\bf I} +\frac{{\bf H} {\bf Q} {\bf H}^{H}}{\sigma^{2}} \right)^{-1} \frac{{\bf H}}{\sigma^{2}}
\]
(use the
well-known inequality $({\bf I} +{\bf UV})^{-1} {\bf U} = {\bf U}({\bf
  I} +{\bf VU})^{-1}$ for ${\bf U} = {\bf H}^{H}$ and ${\bf V} = \frac{{\bf H} {\bf Q}}{\sigma^{2}}$ ). 
We denote by ${\bf M}$ the non negative matrix 
\[
{\bf M} = {\bf H}^{H} \left( {\bf I} +\frac{{\bf H} {\bf Q} {\bf H}^{H}}{\sigma^{2}} \right)^{-1} \frac{{\bf H}}{\sigma^{2}}
\]
and remark that
\begin{equation}
\label{eq:autre-expression-phi''}
 \phi''(\lambda) = - \mathbb{E} \mathrm{Tr} \left[ {\bf M} ({\bf Q}_1 - {\bf Q}_2) {\bf M} ({\bf Q}_1 - {\bf Q}_2) \right]
\end{equation}
or equivalently that
\[
 \phi''(\lambda) = - \mathbb{E} \mathrm{Tr} \left[ {\bf M}^{1/2}  ({\bf Q}_1 - {\bf Q}_2) {\bf M}^{1/2}  {\bf M}^{1/2}  ({\bf Q}_1 - {\bf Q}_2) {\bf M}^{1/2} \right]\ .
\]
As matrix  ${\bf M}^{1/2}  ({\bf Q}_1 - {\bf Q}_2) {\bf M}^{1/2}$ is Hermitian, this of course implies that 
$\phi''(\lambda) \le 0$. The concavity of $\phi$ and of $I$ are established.

% and to write:
%\begin{multline*}
%  \phi''(\lambda) = - \mathbb{E} \mathrm{Tr}\Bigg[ ({\bf Q_1 -Q_2})^{\frac
%    12} \left( {\bf I} +\frac{{\bf H}^H {\bf H Q}}{\sigma^2}
%  \right)^{-1} \frac{{\bf H}^H {\bf H}}{\sigma^2} ({\bf Q_1
%    -Q_2})^{\frac 12} \\
%\times ({\bf Q_1 -Q_2})^{\frac
%    12} \left( {\bf I} +\frac{{\bf H}^H {\bf H Q}}{\sigma^2}
%  \right)^{-1} \frac{{\bf H}^H {\bf H}}{\sigma^2} ({\bf Q_1
%    -Q_2})^{\frac 12} \Bigg]\ .
%\end{multline*}
%Now $\phi''(\lambda)=-\mathbb{E}\, \mathrm{Tr} {\bf M M}^H \le 0$ and the concavity of $\phi$ is established.
\subsubsection{Using an auxiliary channel to establish concavity of $\bar{I}({\bf Q})$}
Denote by $\otimes$ the Kronecker product of matrices. We introduce the following matrices:
$$
{\bf \Delta}={\bf I}_m \otimes {\bf C},\quad \tilde {\bf \Delta}={\bf I}_m \otimes \tilde {\bf C},\quad 
\check{\bf A} = {\bf I}_m \otimes {\bf A},\quad \check {\bf Q} = {\bf I}_m \otimes {\bf Q}\ .
$$
Matrix ${\bf \Delta}$ is of size $rm\times rm$, matrices $\tilde {\bf
  \Delta}$ and $\check {\bf Q}$ are of size $tm \times tm$, and $\check {\bf
  A}$ is of size $rm \times tm$. Let us now introduce:
$$
\check {\bf V} =\frac 1{\sqrt{mt}} {\bf \Delta}^{\frac 12} \check{\bf W} \tilde {\bf \Delta}^{\frac 12} \quad
\textrm{and} \quad \check {\bf H} = \sqrt{\frac K{K+1}} \check{\bf A} + \frac1{\sqrt{K+1}} \check{\bf V}\ ,
$$
where $\check {\bf W}$ is a $rm \times tm$ matrix whose entries are
i.i.d ${\mathcal CN}(0,1)$-distributed random variables. Denote by $I_m(\check {\bf Q})$ the 
EMI associated with channel $\check {\bf H}$:
$$
I_m(\check {\bf Q})= \mathbb{E} \log \det \left( {\bf I}
  +\frac{\check {\bf H} \check {\bf Q} \check{\bf H}^H}{\sigma^2} \right).
$$
Applying Theorem \ref{theo:canonique} to the channel $\check {\bf H}$,
we conclude that $I_m(\check {\bf Q})$ admits an asymptotic
approximation $\bar{I}_m(\check {\bf Q})$ defined by the system
\eqref{eq:expref}-\eqref{eq:expreftildef} and formula
\eqref{eq:expreCbarre}, where one will substitute the quantities related to channel ${\bf H}$ by those related to 
channel $\check {\bf H}$, i.e.:
$$
t \leftrightarrow mt,\quad
r \leftrightarrow mr,\quad
{\bf A} \leftrightarrow \check{\bf A},\quad
{\bf Q} \leftrightarrow \check {\bf Q},\quad
{\bf C} \leftrightarrow {\bf \Delta},\quad
\tilde {\bf C}  \leftrightarrow \tilde {\bf \Delta}\ .
$$
Due to the block-diagonal nature of matrices $\check{\bf A}$, $\check
{\bf Q}$, ${\bf \Delta}$ and $\tilde {\bf \Delta}$, the system
associated with channel $\check {\bf H}$ is exactly the same as the one
associated with channel ${\bf H}$. Moreover, a straightforward computation 
yields:
$$
\frac 1m \bar{I}_m(\check {\bf Q}) =\bar{I}({\bf Q}),\qquad \forall\, m\ge 1\ .
$$
It remains to apply the convergence result \eqref{eq:equivalent1} to conclude that
$$
\lim_{m\rightarrow \infty} \frac 1m I_m(\check {\bf Q}) = \bar{I}({\bf Q})\ .
$$
Since ${\bf Q} \mapsto I_m(\check {\bf Q}) = I_m ({\bf I}_m \otimes
{\bf Q})$ is concave, $\bar{I}$ is
concave as a pointwise limit of concave functions.

\subsubsection{Uniform strict concavity of the EMI of the
  auxiliary channel - Strict concavity of $\bar{I}({\bf Q})$}
In order to establish the strict concavity of $\bar{I}({\bf Q})$, we
shall rely on the following lemma: 
\begin{lemma}\label{unif-strict-con}
Let $\bar{\phi}:[0,1] \rightarrow \mathbb{R}$ be a real function such that there exists a family 
$(\phi_m)_{m\ge 1}$ of real functions satisfying:
\begin{itemize}
\item[(i)] The functions $\phi_m$ are twice differentiable and there exists $\kappa<0$ such that
\begin{equation}\label{strict-con}
\forall m\ge 1,\quad \forall \lambda \in [0,1],\qquad \phi_m''(\lambda) \le \kappa <0\ .
\end{equation}
\item[(ii)] For every $\lambda\in [0,1]$, $\phi_m(\lambda)\xrightarrow[m\rightarrow \infty]{} \bar{\phi}(\lambda)$.
\end{itemize} 
Then $\bar{\phi}$ is a strictly concave real function. 
\end{lemma}
Proof of Lemma \ref{unif-strict-con} 
is postponed to Appendix \ref{app:constric}.

Let ${\bf Q}_1$, ${\bf Q}_2$ in ${\cal C}_1$; denote by ${\bf Q}=\lambda {\bf Q}_1 +(1-\lambda) {\bf Q}_2$, 
$\check {\bf Q}_1= I_m \otimes {\bf Q}_1$, $\check {\bf Q}_2= I_m \otimes {\bf Q}_2$, $\check {\bf Q}
= I_m \otimes {\bf Q}$. Let $\check {\bf H}$ be the matrix associated with the auxiliary channel and denote by: 
$$
\phi_m(\lambda)= \frac 1m \mathbb{E} \log \det \left( {\bf I}
  +\frac{\check {\bf H} \check {\bf Q} \check {\bf H}^H}{\sigma^2} \right)\ .
$$  
We have already proved that $\phi_m(\lambda)
\xrightarrow[m\rightarrow \infty]{} \bar{\phi}(\lambda)\stackrel{\triangle}{=}\bar{I} (\lambda {\bf Q}_1
+(1-\lambda) {\bf Q}_2 )$. In order to fulfill assumptions of
Lemma \ref{unif-strict-con}, it is sufficient to prove that there exists $\kappa<0$ such that
for every $\lambda \in [0,1]$, 
\begin{equation}\label{unif-strict-phi}
\limsup_{m\rightarrow \infty} \phi_m''(\lambda)\le \kappa <0\ .
\end{equation}
(\ref{unif-strict-phi}) is proved in the Appendix \ref{app:constric}.

%++++++++++++++++++++++++++++++++++++++++++++++++++++++++++++++++++++

%++++++++++++++++++++++++++++++++++++++++++++++++++++++++++++++++++++

%++++++++++++++++++++++++++++++++++++++++++++++++++++++++++++++++++++
\subsection{Approximation of the capacity $I({\bf Q}_*)$}

Since $\bar{I}$ is strictly concave over the compact set ${\mathcal
  C}_1$, it admits a unique argmax we shall denote by $\tO$, i.e.:
$$
\bar{I}(\tO)=\max_{{\bf Q} \in {\mathcal C}_1} \bar{I}({\bf Q})\ .
$$
As we shall see in Section \ref{sec:algo}, matrix ${\tO}$ can be
obtained by a rather simple algorithm. Provided that $\sup_t \|\tO\|$
is bounded, Eq. \eqref{eq:equivalent1} in Theorem \ref{theo:canonique}
yields $I(\tO) -\bar{I}(\tO)\rightarrow 0$ as $t\rightarrow \infty$.
It remains to check that $I({\bf Q}_*) -I(\tO)$ goes asymptotically to
zero to be able to approximate the capacity. This is the purpose of
the next proposition.
   
\begin{proposition}\label{approx-capacity}
Assume that $\sup_t \| {\bf A} \| < \infty$,  
$\sup_t \| \tilde{\bf C} \| < \infty$, $\sup_t \| {\bf C} \| < \infty$, 
$\inf_t \lambda_{\mathrm{min}} ( \tilde{\bf C} ) > 0$, and 
$\inf_t \lambda_{\mathrm{min}} ( {\bf C} ) > 0$. 
Let $\tO$ and ${\bf Q}_*$ be the maximizers over ${\cal C}_1$ of $\bar{I}$ and
$I$ respectively. Then the following facts hold true:
\begin{itemize}
\item[(i)] $\sup_t \| \tO \| < \infty$. 
\item[(ii)] $\sup_t \| \mb{Q}_* \| < \infty$. 
\item[(iii)] ${I}(\tO) = I({\bf Q}_*) + O(t^{-1})$. 
\end{itemize} 
\end{proposition}

\begin{proof} The proof of items (i) and (ii) is postponed to Appendix
  \ref{proof-approx-capacity}. Let us prove (iii). As
\begin{multline} 
\begin{array}[t]{c}
\underbrace{\left( I({\bf Q}_*) - I(\tO) \right)} \\
\geq 0 
\end{array}
+ 
\begin{array}[t]{c}
\underbrace{\left( \bar I(\tO) - \bar I({\bf Q}_*) \right)} \\
\geq 0 
\end{array} \\ 
= 
\begin{array}[t]{c}
\underbrace{\left( I({\bf Q}_*) - \bar I({\bf Q}_*) \right)} \\
= O(t^{-1}) \\
\text{by (ii) and Th. \ref{theo:canonique} Eq. \eqref{eq:equivalent1}} 
\end{array} 
+
\begin{array}[t]{c}
\underbrace{\left(\bar I(\tO) - I(\tO) \right)} \\
= O(t^{-1}) \\
\text{by (i) and Th. \ref{theo:canonique} Eq. \eqref{eq:equivalent1}} 
\end{array} 
\label{eq-I(Q)-I(bar Q)} 
\end{multline} 
where the two terms of the lefthand side are nonnegative due to the
fact that ${\bf Q}_*$ and $\tO$ are the maximizers of $I$ and $\bar I$
respectively.  As a direct consequence of \eqref{eq-I(Q)-I(bar Q)}, we
have $I({\bf Q}_*) - I(\tO) = O(t^{-1})$ and the proof is completed.
\end{proof}

\section{Optimization of the input covariance matrix} \label{sec:algo}
%++++++++++++++++++++++++++++++++++++++++++++++++++++++++++++++++++++

%++++++++++++++++++++++++++++++++++++++++++++++++++++++++++++++++++++ 

%++++++++++++++++++++++++++++++++++++++++++++++++++++++++++++++++++++

In the previous section, we have proved that matrix $\tO$
asymptotically achieves the capacity. The purpose of this section is
to propose an efficient way of maximizing the asymptotic approximation
$\bar{I}(\mb{Q})$ without using complicated numerical optimization
algorithms. In fact, we will show that our problem boils down to
simple waterfilling algorithms.

\subsection{Properties of the maximum of  $\bar{I}(\mb{Q})$.}

In this section, we shall establish some of $\tO$'s properties. We
first introduce a few notations.  Let $V(\kappa, \tilde{\kappa},
\mb{Q})$ be the function defined by:
\begin{multline}
\label{eq:defv-1}
V(\kappa, \tilde{\kappa}, \mb{Q}) = \log \det \left( {\bf I}_t+
  \frac{\kappa}{K+1} \, {\bf Q}^{\demi} \tilde{{\bf C}} {\bf
    Q}^{\demi} + \frac{K}{\sigma^2(K+1)} \,
  {\bf Q}^{\demi} {\bf A}^H \left( {\bf I}_r+\frac{\tilde{\kappa}}{K+1} \,{\bf C} 
\right)^{-1}{\bf A}{\bf Q}^{\demi}\right)\\
+ \log\det \left({\bf I}_r+ \frac{\tilde{\kappa}}{K+1} {\bf
    C}\right)-  \frac{t\sigma^2 \kappa\tilde \kappa}{K+1} \, .
\end{multline}
or equivalently by
\begin{multline}
\label{eq:defv-2}
V(\kappa, \tilde{\kappa}, \mb{Q}) =  \log\det \left( {\bf I}_r+ \frac{\tilde{\kappa}}{K+1}{\bf C}
+\frac{K}{\sigma^2(K+1)} {\bf A}{\bf Q}^{\demi} \left( {\bf
I}_t+\frac{\kappa}{K+1} {\bf Q}^{\demi} {\bf \tilde{C}} {\bf Q}^{\demi} \right)^{-1} {\bf Q}^{\demi}{\bf A}^H\right) \\
+ \log\det \left({\bf I}_t+ \frac{\kappa}{K+1} {\bf Q}^{1/2} {\bf \tilde{C}} {\bf Q}^{1/2} \right)
-\frac{t\sigma^2 \kappa \tilde{\kappa}}{K+1}\ .
\end{multline}
Note that if $(\delta(\mb{Q}), \tilde{\delta}(\mb{Q}))$ is the solution of system (\ref{eq:canonique}), then:
$$
\bar{I}({\bf Q}) = V (\delta(\mb{Q}), \tilde{\delta}(\mb{Q}), {\bf Q})\ .
$$
Denote by $(\delta_{*},\tilde{\delta}_{*})$ the solution
$(\delta(\tO), \tilde{\delta}(\tO))$ of (\ref{eq:canonique})
associated with $\tO$. The aim of the section is to prove that $\tO$
is the solution of the following standard waterfilling problem:
$$
\bar{I}(\tO) =\max_{{\bf Q} \in {\cal C}_1} V(\delta_*, \tilde \delta_*, {\bf Q})\ .
$$

Denote by ${\bf G}(\kappa, \tilde{\kappa})$ the $t \times t$ matrix given by:
\begin{equation}
\label{eq:defG}
{\bf G}(\kappa, \tilde{\kappa}) =  \frac{\kappa}{K+1} \tilde{{\bf C}} +  \frac{K}{\sigma^2(K+1)}  {\bf A}^H \left( {\bf I}_r+\frac{\tilde{\kappa}}{K+1} \,{\bf C} \right)^{-1}{\bf A}\ .
\end{equation}
Then,  $V(\kappa, \tilde{\kappa}, \mb{Q})$ 
also writes
\begin{equation}
\label{eq:expreVter}
V(\kappa, \tilde{\kappa}, \mb{Q})  =  \log \mathrm{det} \left( {\bf I} + {\bf Q} {\bf G}(\kappa, \tilde{\kappa}) \right)
+ \log\det \left({\bf I}_r+ \frac{\tilde{\kappa}}{K+1} {\bf C}\right) -\frac{t \sigma^2 \kappa \tilde{\kappa}}{K+1}\ ,
\end{equation}
which readily implies the differentiability
of $(\kappa, \tilde \kappa, \tQ) \mapsto V(\kappa, \tilde \kappa,
\tQ)$ and the strict concavity of  $\tQ \mapsto V(\kappa, \tilde \kappa,
\tQ)$ ($\kappa$ and $\tilde \kappa$ being frozen).

In the sequel, we will denote by $\nabla F(x)$ the derivative of the
differentiable function $F$ at point $x$ ($x$ taking its values in
some finite-dimensional space) and by $\langle \nabla F(x), y\rangle$
the value of this derivative at point $y$. Sometimes, a function is
not differentiable but still admits {\em directional derivatives}: The
directional derivative of a function $F$ at $x$ in direction $y$ is
$$
F'(x;y) =\lim_{t\downarrow  0} \frac {F(x+ty) -F(x)}t
$$
when the limit exists. Of course, if $F$ is differentiable at $x$, then 
$F'(x;y) = \langle \nabla F(x), y \rangle$. The following proposition captures the main features 
needed in the sequel.
\begin{proposition}
\label{prop:caracterisation}
Let $F:{\cal C}_1 \rightarrow \mathbb{R}$ be a concave function. Then:
\begin{itemize}
\item[(i)] The directional derivative $F'(\tQ; \tP-\tQ)$ exists in $(-\infty,\infty]$ for all $\tQ,\tP$ in ${\cal C}_1$.
\item[(ii)] ({\em necessary condition}) If $F$ attains
  its maximum for $\tO \in {\cal C}_1$, then:
\begin{equation}\label{opti-ness}
\forall \tQ\in {\cal C}_1,\quad F'(\tO; \tQ -\tO) \le 0\ .
\end{equation}
\item[(iii)] ({\em sufficient condition}) Assume that there exists $\tO\in {\cal C}_1$ such that:
\begin{equation}\label{opti-suff}
\forall \tQ \in {\cal C}_1,\quad  F'(\tO;\tQ -\tO) \le 0.
\end{equation}
Then $F$ admits its maximum at $\tO$ (i.e. $\tO$ is an argmax of $F$ over ${\cal C}_1$).
\end{itemize} 
If $F$ is differentiable then both conditions \eqref{opti-ness} and \eqref{opti-suff} write:
$$
\forall \tQ \in {\cal C}_1,\quad \langle \nabla  F(\tO),\tQ -\tO  \rangle \le 0.
$$
\end{proposition}
Although this is standard material (see for
instance \cite[Chapter 2]{BorLew00}), we provide some elements of
proof for the reader's convenience.

\begin{proof}
  Let us first prove item (i). As $\tQ+t(\tP-\tQ)= (1-t)\tQ +t\tP \in {\cal C}_1$, 
$\Delta(t)\stackrel{\triangle}{=}t^{-1}\left(
    F(\tQ+t(\tP-\tQ)) - F(\tQ)\right)$ is well-defined. Let $0\le s\le t\le 1$ and
  consider 
\begin{eqnarray*}
\Delta(t) -\Delta(s) &=&\frac 1s \left\{ \frac st F\left( (1-t)\tQ + t\tP \right) +\frac{t-s}t F({\bf Q})
- F\left( (1-s)\tQ + s\tP \right) \right\}\ ,\\
&\stackrel{(a)}{\le}& \frac 1s \left\{ F\left(s\frac{(1-t)\tQ + t\tP }t +\frac{t-s}t {\bf Q} \right)
- F\left( (1-s)\tQ + s\tP \right) \right\}\ ,\\
&=& \frac 1s \left\{ F\left( (1-s)\tQ + s\tP \right) -F\left( (1-s)\tQ + s\tP \right) \right\}\quad =\quad 0 \,,
\end{eqnarray*}
where $(a)$ follows from the concavity of $F$. This shows that
$\Delta(t)$ increases as $t \downarrow 0$, and in particular always
admits a limit in $(-\infty,\infty]$.

Item (ii) readily follows from the fact that $F((1-t)\tO +t\tP) \le
F(\tO)$ due to the mere definition of $\tO$. This implies that
$\Delta(t)\le 0$ which in turn yields \eqref{opti-ness}.

We now prove (iii). The concavity of $F$ yields:
$$
\Delta(t)=\frac{F(\tO+t(\tP-\tO)) - F(\tO) }t \ge F(\tP) -F(\tO). 
$$
As $\lim_{t\downarrow 0} \Delta(t) \le 0$ by \eqref{opti-suff}, one gets: $\forall \tP \in
{\cal C}_1$, $F(\tP) -F(\tO)\le 0$. Otherwise stated, $F$ attains its
maximum at $\tO$ and Proposition \ref{prop:caracterisation} is proved.
%
%  As $W$ is differentiable, it is in particular continuous over
%  ${\mathcal C}_1$.  The existence and uniqueness of $\tO$ follows
%  from $W$'s continuity and strict concavity together with the fact
%  that ${\cal C}_1$ is compact.
%
%
%Let $\tP \in {\cal C}_1$, and consider the function
%$\phi$ defined for $\lambda\in[0,1]$ by:
%$$
%\phi(\lambda) = \frac{W \left( \tO + \lambda(\tP - \tO) \right) - W(\tO)}{\lambda}\ .
%$$
%The mere definition of $\tO$ yields $\phi(\lambda) \leq 0$ for every
%$\lambda \in [0,1]$. This in particular implies that $ \langle \nabla 
%W({\bf Q_*}),{\bf P} -{\bf Q}_*\rangle =\lim_{\lambda \rightarrow
%  0^{+}}\phi(\lambda)$ is lower than 0, which is exactly condition
%(\ref{eq:condition-opt}). Conversely, let $\tO'$ be an element of
%${\cal C}_1$ satisfying (\ref{eq:condition-opt}). As $W$ is
%concave, $\phi(\lambda)\ge W(\tP) - W(\tO')$. Taking
%$\lambda \rightarrow 0^{+}$ and applying (\ref{eq:condition-opt})
%yields 
%$$ 
%0\ge \langle \nabla  W({{\bf Q}'_*}),{\bf P} -{\bf Q}'_*\rangle \geq W(\tP) -
%W({\bf Q}_*')\ ,\qquad \forall \tP \in {\cal C}_1\ ,
%$$ 
%in particular $W(\tO')\ge W(\tP)$ for every $\tP \in {\cal C}_1$. 
%Finally, the argmax being unique, $\tO' = \tO$
\end{proof}

In the following proposition, we gather various properties related
to $\bar{I}$.

\begin{proposition}\label{properties-ibar}
Consider the functions $\delta({\bf Q}), \tilde \delta({\bf Q})$ and $\bar{I}({\bf Q})$ from ${\cal C}_1$ to $\mathbb{R}$.
The following properties hold true:
\begin{itemize}
\item[(i)] Functions $\delta({\bf Q}), \tilde \delta({\bf Q})$ and
  $\bar{I}({\bf Q})$ are differentiable (and in particular continuous)
  over ${\cal C}_1$.
\item[(ii)] Recall that $\tO$ is the argmax of $\bar{I}$ over 
${\cal C}_1$, i.e.
$
\forall {\bf Q} \in {\cal C}_1,\ \bar{I}({\bf Q}) \le \bar{I}(\tO)\ .
$
Let ${\bf Q}\in {\cal C}_1$. The following property:
$$
\forall {\bf P} \in {\cal C}_1,\quad \langle \nabla  \bar{I}({\bf Q}),{\bf P} -{\bf Q}\rangle \le 0
$$
holds true if and only if ${\bf Q}= \tO$.
\item[(iii)] \label{prop:waterfilling} Denote by $\delta_{*}$ and
  $\tilde{\delta}_{*}$ the quantities $\delta(\tO)$ and
  $\tilde{\delta}(\tO)$. Matrix $\tO$ is the solution
  of the standard waterfilling problem: Maximize over ${\bf Q} \in
  {\cal C}_1$ the function $V(\delta_{*}, \tilde{\delta}_{*}, \tQ)$ or
  equivalently the function $\log \mathrm{det}({\bf I} + {\bf Q} {\bf
    G}(\delta_{*}, \tilde{\delta}_{*}))$.
\end{itemize}
\end{proposition}

\begin{proof}
  (i) is established in the Appendix. 
Let us establish (ii). Recall that $\bar{I}({\bf Q})$ is strictly
concave by Theorem \ref{theo:con} (and therefore its
maximum is attained at at most one point). On the other hand,
$\bar{I}({\bf Q})$ is continuous by (i) over ${\cal C}_1$ which is compact. Therefore, the maximum of 
$\bar{I}({\bf Q})$ is uniquely attained at a point $\tO$. Item (ii) follows then from
Proposition \ref{prop:caracterisation}.

Proof of item (iii) is based on the following identity, to be proved below:
\begin{equation}
\label{eq:egalite}
\langle \nabla \bar{I}({\tO}),\tQ -  \tO \rangle = \langle \nabla_{\tQ} V \left( \delta_{*}, \tilde{\delta}_{*},\tO \right) ,\tQ -  \tO\rangle\ ,
\end{equation}
where $\nabla_{\tQ}$ denote the derivative of $V(\kappa, \tilde
\kappa, \tQ)$ with respect to $V$'s third component, i.e.
$\nabla_{\tQ} V(\kappa, \tilde \kappa, \tQ) = \nabla \Gamma (\tQ)$
with $\Gamma: \tQ \mapsto V(\kappa, \tilde \kappa, \tQ)$.  Assume that
(\ref{eq:egalite}) holds true. Then item (ii) implies that $\langle
\nabla_{\tO} V \left( \delta_{*}, \tilde{\delta}_{*},\tO \right) ,\tQ
- \tO\rangle \le 0$ for every $\tQ \in {\cal C}_1$. As ${\bf Q}
\mapsto V(\delta_{*}, \tilde{{\delta}}_{*}, {\bf Q})$ is strictly
concave on ${\cal C}_1$, $\tO$ is the argmax of $V(\delta_{*},
\tilde{{\delta}}_{*}, \cdot )$ by Proposition \ref{prop:caracterisation} and we are done.

It remains to prove \eqref{eq:egalite}. Consider $\tQ$ and $\tP$ in ${\cal C}_1$, and use the identity
\begin{multline*}
  \langle \nabla \bar{I}(\tP), \tQ - \tP \rangle = \langle \nabla_{\bf Q} V(\delta(\tP),
  \tilde{\delta}(\tP), \tP), \tQ -\tP)\rangle  \\
  + \left( \frac{\partial V}{\partial \kappa} \right)(\delta(\tP), \tilde{\delta}(\tP), \tP) 
  \ \langle \nabla \delta(\tP), \tQ -  \tP \rangle \\ 
  + \left( \frac{\partial V}{\partial \tilde{\kappa}}
  \right)(\delta(\tP), \tilde{\delta}(\tP), \tP) \ \langle
  \nabla \tilde{\delta}(\tP), \tQ - \tP\rangle\ .
\end{multline*}
We now compute the partial derivatives of $V$ and obtain:
\begin{equation}\label{eq:derivees}
\left\{
\begin{array}{ccc}
\displaystyle{\f{\partial V}{\partial \kappa}} & = & \displaystyle{- \frac{t
\sigma^{2}}{K+1} \left( \tilde{\kappa} - \tilde{f}(\kappa, \tilde{\kappa}, \mb{Q})
\right)}\\
\displaystyle{\f{\partial V}{\partial \tilde{\kappa}}} & = & \displaystyle{ - \frac{t
\sigma^{2}}{K+1} \left( \kappa - f(\kappa, \tilde{\kappa}, \mb{Q}) \right)}
\end{array}
\right.\ ,
\end{equation}
where $f$ and $\tilde f$ are defined by \eqref{eq:expref} and
\eqref{eq:expreftildef}.  The first relation follows from 
(\ref{eq:defv-1}) and the second relation from (\ref{eq:defv-2}).  
As
$(\delta({\bf Q}), \tilde{\delta}({\bf Q}))$ is the solution of 
system (\ref{eq:canonique}), equations (\ref{eq:derivees}) imply that:
\begin{equation}
\label{eq:crucial}
\frac{\partial V}{\partial \kappa} (\delta({\bf Q}), \tilde{\delta}({\bf Q}), {\bf Q}) =  
\frac{\partial V}{\partial \tilde{\kappa}} (\delta({\bf Q}), \tilde{\delta}({\bf Q}), {\bf Q})  =  0\ .
\end{equation}
Letting $\tP = \tO$ and taking into account (\ref{eq:crucial}) yields:
\[
\langle \nabla \bar{I}(\tO), \tQ - \tO \rangle = \langle \nabla_{\tQ}
V(\delta(\tO), \tilde{\delta}(\tO), \tO), \tQ -\tO \rangle\ ,
\]
and (iii) is established.
\end{proof}

\begin{remark} The quantities $\delta_{*}$ and $\tilde{\delta}_*$
  depend on matrix $\tO$. Therefore, Proposition
  \ref{prop:waterfilling} does not provide by itself any optimization
  algorithm. However, it gives valuable insights on the structure of
  $\tO$. Consider first the case ${\bf C} = {\bf I}$ and
  $\tilde{{\bf C}} = {\bf I}$. Then, $ {\bf G}(\delta_{*},
  \tilde{\delta}_{*})$ is a linear combination of ${\bf I}$ and matrix
  ${\bf A}^{H} {\bf A}$.  The eigenvectors of $\tO$ thus
  coincide with the right singular vectors of matrix ${\bf A}$, a
  result consistent with the work \cite{Hoesli-Kim-Lapidoth-05}
  devoted to the maximization of the EMI $I({\bf \tQ})$. If ${\bf C} =
  {\bf I}$ and $\tilde{{\bf C}} \neq {\bf I}$, $ {\bf G}(\delta_{*},
  \tilde{\delta}_{*})$ can be interpreted as a linear combination of
  matrices $\tilde{{\bf C}}$ and ${\bf A}^{H} {\bf A}$. Therefore, if
  the transmit antennas are correlated, the eigenvectors of the
  optimum matrix $\tO$ coincide with the eigenvectors of some
  weighted sum of $\tilde{{\bf C}}$ and ${\bf A}^{H} {\bf A}$. This
  result provides a simple explanation of the impact of correlated
  transmit antennas on the structure of the optimal input covariance
  matrix. The impact of correlated receive antennas on $\tO$
  is however less intuitive because matrix ${\bf A}^{H} {\bf A}$ has
  to be replaced with ${\bf A}^{H}({\bf I} + \tilde{\delta}_{*} {\bf
    C})^{-1} {\bf A}$.
\end{remark}

\subsection{The optimization algorithm.}

We are now in position to introduce our maximization algorithm of $\bar{I}$. It is mainly motivated
by the simple observation that for each fixed $(\kappa, \tilde{\kappa})$, the maximization w.r.t.
$\tQ$ of function $V(\kappa, \tilde{\kappa}, \tQ)$ defined by (\ref{eq:expreVter}) can be 
achieved by a standard waterfilling procedure, which,
of course, does not need the use of numerical techniques. On the other hand,
for $\tQ$ fixed, the equations (\ref{eq:canonique}) have unique solutions that, in practice,
can be obtained using a standard fixed-point algorithm. Our algorithm thus consists in adapting
parameters $\tQ$ and $\delta, \tilde{\delta}$ separately by the following iterative scheme:
\begin{itemize}
\item Initialization: $\tQ_{0} = {\bf I}$, $(\delta_1, \tilde{\delta}_1)$ are defined as the unique
solutions of system (\ref{eq:canonique}) in which $\tQ = \tQ_{0} = {\bf I}$.
Then, define $\tQ_1$ are the maximum of function $\tQ \rightarrow V(\delta_1, \tilde{\delta}_1, \tQ)$
on ${\cal C}_1$, which is obtained through a standard waterfilling procedure.
\item Iteration $k$: assume $\tQ_{k-1}$, $(\delta_{k-1}, \tilde{\delta}_{k-1})$ available. Then,
 $(\delta_k, \tilde{\delta}_k)$ is defined as the unique solution of (\ref{eq:canonique}) in which $\tQ = \tQ_{k-1}$.
Then, define $\tQ_k$ are the maximum of function $\tQ \rightarrow V(\delta_{k}, \tilde{\delta}_{k}, \tQ)$
on ${\cal C}_1$.
\end{itemize}
One can notice that this algorithm is the generalization of the
procedure used by \cite{Wen-Com-06} for
optimizing the input covariance matrix for correlated Rayleigh
MIMO channels. \\ 

We now study the convergence properties of this algorithm, and state a result which 
implies that, if the algorithm converges, then it converges to the unique 
argmax $\tO$ of $\bar{I}$.
\begin{proposition}
\label{prop:convergence-algo}
Assume that the two sequences $(\delta_k)_{k \geq 0}$ and  $(\tilde{\delta}_k)_{k \geq 0}$ verify
\begin{equation}
\label{eq:condition}
\lim_{k \rightarrow +\infty} \delta_k - \delta_{k-1} \rightarrow 0, \lim_{k \rightarrow +\infty} \tilde{\delta}_k - \tilde{\delta}_{k-1} \rightarrow 0
\end{equation}
Then, the sequence $(\tQ_{k})_{k \geq 0}$ converges toward the maximum $\tO$ of $\bar{I}$ on ${\cal C}_1$.
\end{proposition}
The proof is given in the appendix. 

\begin{remark}
If the algorithm is convergent,
i.e. if sequence $({\bf Q}_k)_{k \geq 0}$ converges towards a matrix
${\bf P}_{*}$, Proposition \ref{prop:convergence-algo} implies that
${\bf P}_{*} = \tO$. In fact, functions ${\bf Q} \mapsto
\delta({\bf Q})$ and ${\bf Q} \mapsto \tilde{\delta}({\bf Q})$ are
continuous by Proposition \ref{properties-ibar}. As $\delta_k = \delta({\bf Q}_{k-1})$ and
$\tilde{\delta}_k = \tilde{\delta}({\bf Q}_{k-1})$, the convergence of
$({\bf Q}_k)$ thus implies the convergence of
$(\delta_k)$ and $(\tilde{\delta}_k)$, and (\ref{eq:condition}) is fulfilled.
Proposition \ref{prop:convergence-algo} immediately yields 
${\bf P}_{*} = \tO$. Although we have not been able to prove the convergence of
the algorithm, the above result is encouraging, and tends to indicate
the algorithm is reliable. In particular, all the numerical
experiments we have conducted indicates that the algorithm converges
towards a certain matrix which must coincide by Proposition
\ref{prop:convergence-algo} with $\tO$.
\end{remark}

%--------------------------------------------------------------------

\section{Numerical experiments.} \label{sec:simulations}

%++++++++++++++++++++++++++++++++++++++++++++++++++++++++++++++++++++

%++++++++++++++++++++++++++++++++++++++++++++++++++++++++++++++++++++

%++++++++++++++++++++++++++++++++++++++++++++++++++++++++++++++++++++

%--------------------------------------------------------------------

\subsection{When is the number of antennas large enough to reach the asymptotic regime?}
\label{validation}
 All our analysis is based on the approximation
of the ergodic mutual information. This approximation consists in
assuming the channel matrix to be large. Here we provide typical
simulation results showing that the asymptotic regime is reached
for relatively small number of antennas. For the simulations
provided here we assume:
\begin{itemize}
    \item  $\mb{Q} = \mb{I}_t$.
    \item The chosen line-of-sight (LOS) component $\mb{A}$ is based on
 equation (\ref{eq:exempleA}). The angle of arrivals are chosen
 randomly according to a uniform distribution. 
    \item Antenna correlation is assumed to decrease exponentially
    with the inter-antenna distance i.e. $\tilde{{\bf C}}_{ij} \sim
    \rho_{T}^{|i-j|}$, ${\bf C}_{ij} \sim \rho_{R}^{|i-j|}$ with $0 \leq \rho_T \leq
    1$ and $0 \leq \rho_R \leq 1$.
    \item $K$ is equal to $1$.
    \end{itemize}

\begin{figure}[ht]
  \centerline{\includegraphics[width=14cm]{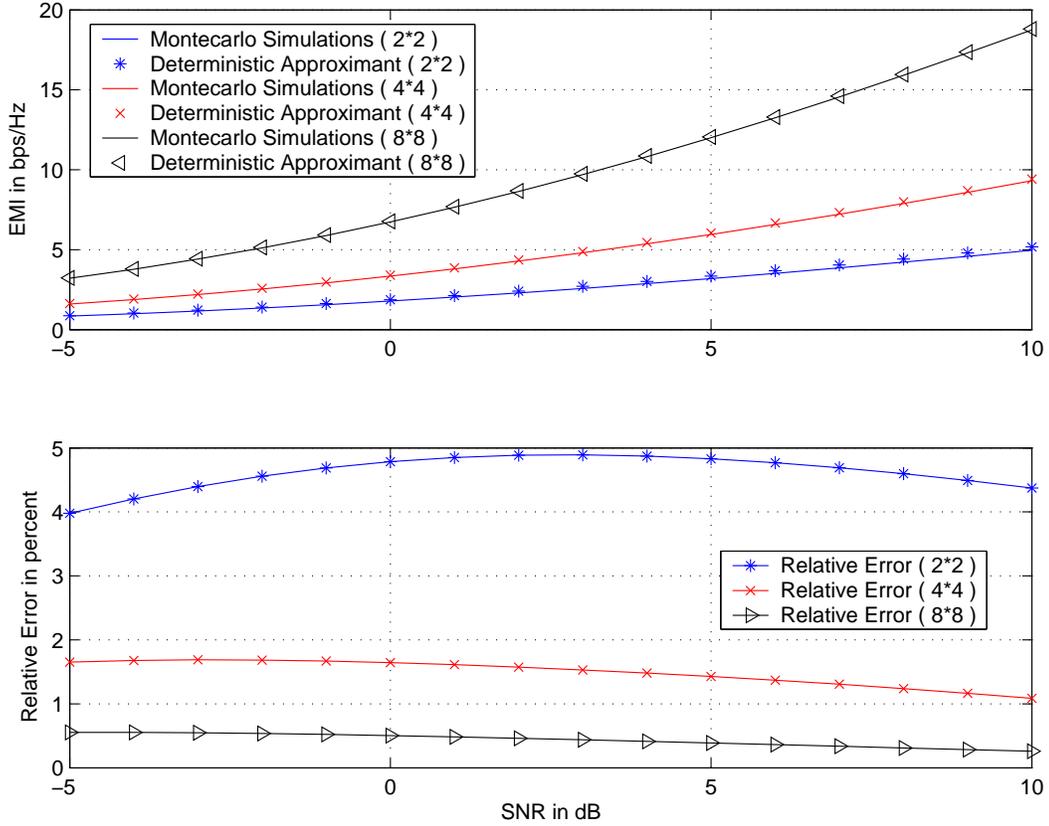}}
\caption{The large system approximation is accurate for correlated
Rician MIMO channels. The relative difference between the EMI
approximation and that obtained by Monte-Carlo simulations is less
than $5 \ \%$ for a $2 \times 2$ system and less than $1 \ \%$ for
a $8 \times 8$ system.}\label{fig:capa}
\end{figure}

Figure \ref{fig:capa} represents the EMI $I({\bf Q})$ evaluated 
by Monte Carlo simulations and its approximation 
$\bar{I}({\bf Q})$ as well as their relative difference (in
percentage). Here, the correlation coefficients are equal to $(\rho_T, \rho_R)= (0.8,
0.3)$ and three different pairs of numbers of antenna are considered: $(t,r) \in
\{(2,2), (4,4), (8,8)\}$. Figure \ref{fig:capa} shows that the approximation 
is reliable even for $r=t=2$ in a wide range of SNR. 

%--------------------------------------------------------------------

%--------------------------------------------------------------------

\subsection{Comparison with the Vu-Paulraj method.}
\begin{figure}[ht]
\centerline{\includegraphics[scale=0.6]{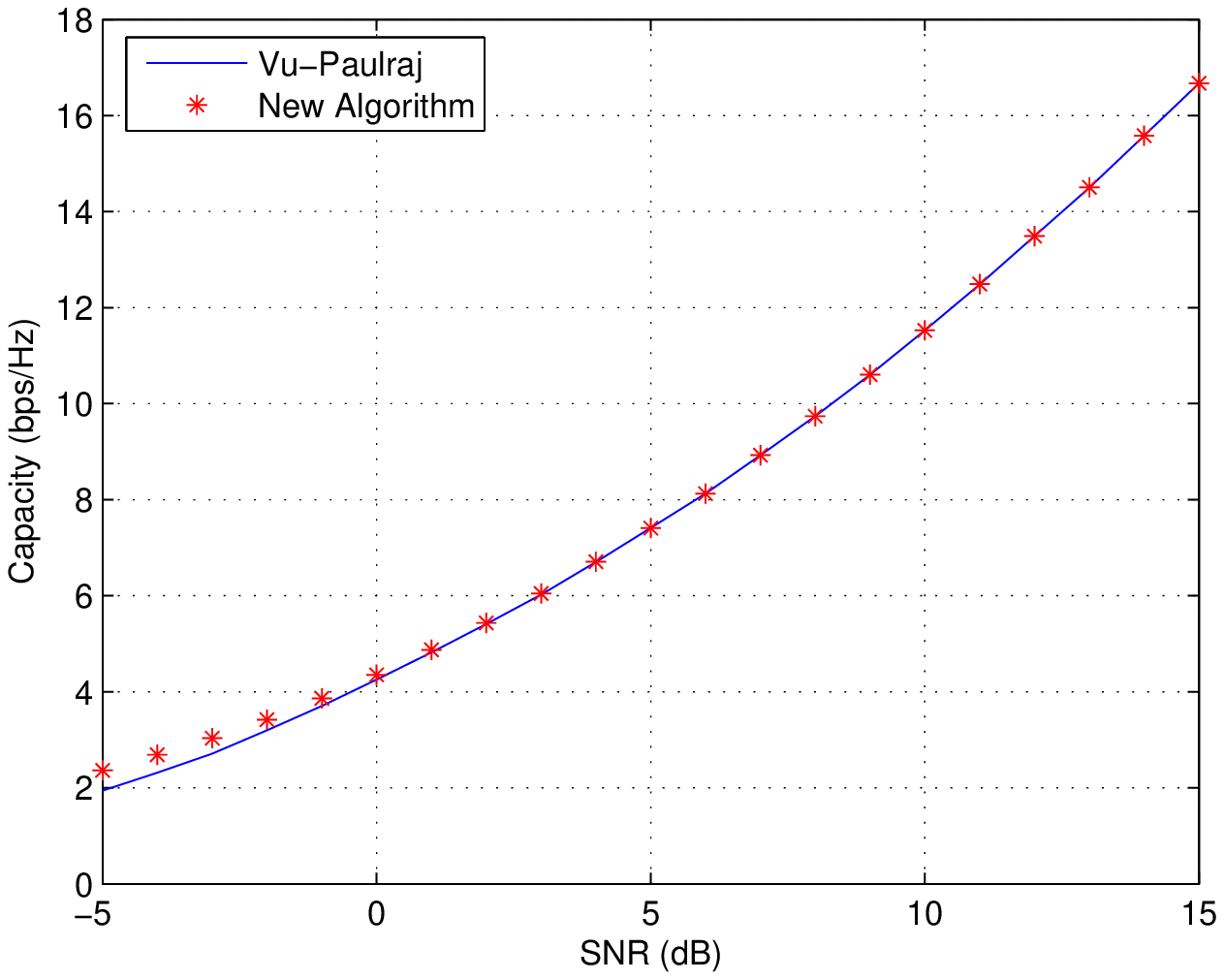}}\caption{Comparison with
the Vu-Paulraj algorithm I} \label{fig:canal-paulraj}
\end{figure}

In this paragraph, we compare our algorithm with the method presented
in \cite{Vu-Paulraj-05} based on the maximization of $I({\bf Q})$.  We
recall that Vu-Paulraj's algorithm is based on a Newton method and a
barrier interior point method. Moreover, the average mutual
informations and their first and second derivatives are evaluated by
Monte-Carlo simulations. In fig.  \ref{fig:canal-paulraj}, we have
evaluated $C_E = \max_{ {\bf Q} \in {\cal C}_1} I({\bf Q})$ versus the
SNR for $r=t=4$. Matrix ${\bf H}$ coincides with the example
considered in \cite{Vu-Paulraj-05}. The solid line corresponds to the
results provided by the Vu-Paulraj's algorithm; the number of trials
used to evaluate the mutual informations and its first and second
derivatives is equal to $30.000$, and the maximum number of iterations
of the algorithm in \cite{Vu-Paulraj-05} is fixed to 10. The dashed
line corresponds to the results provided by our algorithm: Each point
represents $I(\tO)$ at the corresponding SNR, where $\tO$ is the
argmax of $\bar{I}$; the average mutual information at point $\tO$ is
evaluted by Monte-Carlo simulation (30.000 trials are used). The
number of iterations is also limited to 10. Figure
\ref{fig:canal-paulraj} shows that our asymptotic approach provides
the same results than the Vu-Paulraj's algorithm. However, our
algorithm is computationally much more efficient as the above table
shows.  The table gives the average executation time (in sec.) of one
iteration for both
algorithms for $r=t=2, r=t=4, r=t=8$. \\

\begin{figure}[ht]
\centerline{\includegraphics[scale=0.6]{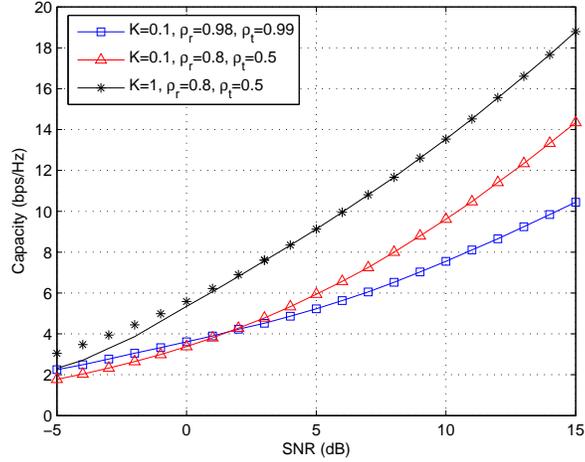}} \caption{Comparison with
the Vu-Paulraj algorithm II} \label{fig:canal-aleatoire}
\end{figure}

In fig. \ref{fig:canal-aleatoire}, we again compare Vu-Paulraj's algorithm and our proposal.
Matrix ${\bf A}$ is generated according to (\ref{eq:exempleA}),
the angles being chosen at random. The transmit and receive antennas correlations are exponential with parameter
$0 < \rho_{T} < 1$ and $0 < \rho_{R} < 1$ respectively. In the experiments, $r=t=4$, while various values of
$\rho_T$, $\rho_R$ and of the Rice factor $K$ have been considered. As in the previous experiment,
the maximum number of iterations for both algorithms is 10, while the number of trials generated to
evaluate the average mutual informations and their derivatives is equal to 30.000. Our approach again provides the same
results than Vu-Paulraj's algorithm, except for low SNRs for $K=1, \rho_T=0.5, \rho_R=0.8$ where our method
gives better results: at these points, the Vu-Paulraj's algorithm seems not to have converge
at the 10th iteration.

\begin{figure}
\label{table}
\begin{center}
\begin{tabular}{c|c|c|c|}
\cline{2-4}
& $n=N=2$ & $n=N=4$ & $n=N=8$ \\
\hline
\multicolumn{1}{|c|}{Vu-Paulraj} & $0.75$ & $8.2$& $138$\\
\hline
\multicolumn{1}{|c|}{New algorithm} & $10^{-2}$ & $3.10^{-2}$ & $7.10^{-2}$ \\
\hline
\end{tabular}
\end{center}
\caption{Average time per iteration in seconds}
\end{figure}

%++++++++++++++++++++++++++++++++++++++++++++++++++++++++++++++++++++

%++++++++++++++++++++++++++++++++++++++++++++++++++++++++++++++++++++

%++++++++++++++++++++++++++++++++++++++++++++++++++++++++++++++++++++

\section{Conclusions} \label{sec:conclusion}

%++++++++++++++++++++++++++++++++++++++++++++++++++++++++++++++++++++

%++++++++++++++++++++++++++++++++++++++++++++++++++++++++++++++++++++

%++++++++++++++++++++++++++++++++++++++++++++++++++++++++++++++++++++

In this paper, an explicit approximation for the ergodic mutual
information for Rician MIMO channels with transmit and receive antenna
correlation is provided. This approximation is based on the asymptotic
Random Matrix Theory. The accuracy of the approximation has been
studied both analytically and numerically. It has been shown to be
very accurate even for small MIMO systems: The relative error is less
than $5 \% $ for a $2 \times 2$ MIMO channel and less $1 \ \%$ for an
$8 \times 8$ MIMO channel.

The derived expression for the EMI has been exploited 
to derive an efficient optimization algorithm providing the optimum covariance matrix.

%++++++++++++++++++++++++++++++++++++++++++++++++++++++++++++++++++++

%++++++++++++++++++++++++++++++++++++++++++++++++++++++++++++++++++++

%++++++++++++++++++++++++++++++++++++++++++++++++++++++++++++++++++++

\appendices

\section{Proof of the existence and uniqueness of the system (\ref{eq:canonique-preparatoire}).}
\label{proof-existence-unicite}
We consider functions $g(\kappa, \tilde{\kappa})$ and  $\tilde{g}(\kappa, \tilde{\kappa})$ defined by 
\begin{equation}
\label{eq:defgtildeg}
\begin{array}{l}
g(\kappa, \tilde{{\kappa}})  
= \displaystyle{ \frac{1}{\kappa} \, \frac{1}{t} \textrm{Tr} \left[ {\bf D} \left( \sigma^{2}({\bf I}_r+{\bf D} \tilde{\kappa}) + {\bf B} ({\bf I}_t+\tilde{{\bf D}}\kappa)^{-1} {\bf B}^{H} \right)^{-1} \right] }
\\
\\
\tilde{g}(\kappa, \tilde{{\kappa}}) = \displaystyle{ \frac{1}{\tilde{\kappa}} \, \frac{1}{t} \textrm{Tr} \left[ \tilde{{\bf D}} \left( \sigma^{2}({\bf I}_t+\tilde{{\bf D}} \kappa) + {\bf B}^{H} 
({\bf I}_r+ {\bf D} \tilde{\kappa})^{-1} {\bf B} \right)^{-1} \right]}
\end{array} .
\end{equation}
For each $\tilde{\kappa} > 0$ fixed, function $\kappa \rightarrow g(\kappa, \tilde{\kappa})$ is clearly strictly 
decreasing, converges toward $+\infty$ if $\kappa \rightarrow 0$ and converges to $0$ if 
$\kappa \rightarrow +\infty$. Therefore, there exists a unique $\kappa > 0$ satisfying $g(\kappa, \tilde{\kappa}) = 1$. 
As this solution depends on $\tilde{\kappa}$, it is denoted $h(\tilde{\kappa})$ in the following. We claim that
\begin{itemize}
\item (i) Function $\tilde{\kappa} \rightarrow h(\tilde{\kappa})$ is strictly decreasing,
\item (ii) Function  $\tilde{\kappa} \rightarrow \tilde{\kappa} h(\tilde{\kappa})$ is strictly increasing.
\end{itemize}
In fact, consider $\tilde{\kappa}_2 > \tilde{\kappa}_1$. It is easily checked that 
for each $\kappa > 0$, $g(\kappa,  \tilde{\kappa}_1) > g(\kappa,  \tilde{\kappa}_2)$. Hence, 
the solution $h(\tilde{\kappa}_1)$ and  $h(\tilde{\kappa}_2)$ of the equations 
$g(\kappa,  \tilde{\kappa}_1) = 1$ and $g(\kappa,  \tilde{\kappa}_2) = 1$ satisfy 
$h(\tilde{\kappa}_1) > h(\tilde{\kappa}_2)$. This establishes (i). To prove (ii), we use the
obvious relation $g(h(\tilde{\kappa}_1), \tilde{\kappa}_1) - g(h(\tilde{\kappa}_2), \tilde{\kappa}_2) = 0$. 
We denote by $({\bf U}_i)_{i=1,2}$ the matrices 
\[
{\bf U}_i = \sigma^{2} \left( h(\tilde{\kappa}_i) {\bf I} + \tilde{\kappa_i} h(\tilde{\kappa}_i) {\bf D} \right)  + 
{\bf B} \left( \frac{{\bf I}}{h(\tilde{\kappa}_i)}  +  \tilde{{\bf D}} \right)^{-1} {\bf B}^{H} 
\]
It is clear that $g(h(\tilde{\kappa}_i), \tilde{\kappa}_i) = \frac{1}{t} \mathrm{Tr} {\bf D} {\bf U}_i^{-1}$. 
We express $g(h(\tilde{\kappa}_1), \tilde{\kappa}_1) - g(h(\tilde{\kappa}_2), \tilde{\kappa}_2)$ 
as 
\[
g(h(\tilde{\kappa}_1), \tilde{\kappa}_1) - g(h(\tilde{\kappa}_2), \tilde{\kappa}_2) = \frac{1}{t} \mathrm{Tr} 
{\bf D} ({\bf U}_1^{-1} - {\bf U}_2^{-1})
\]
and use the identity 
\begin{equation}
\label{eq:difference-inverses}
{\bf U}_1^{-1} - {\bf U}_2^{-1} = {\bf U}_1^{-1} \left( {\bf U}_2 - {\bf U}_1 \right) {\bf U}_2^{-1}\ .
\end{equation}
Using the form of matrices $({\bf U}_i)_{i=1,2}$, we eventually obtain that 
\[
g(h(\tilde{\kappa}_1), \tilde{\kappa}_1) - g(h(\tilde{\kappa}_2),
\tilde{\kappa}_2) = u (h(\tilde{\kappa}_2) -h(\tilde{\kappa}_1)) + v
(\tilde{\kappa}_2 h(\tilde{\kappa}_2) - \tilde{\kappa}_1
h(\tilde{\kappa}_1))\ ,
\]
where $u$ and $v$ are the strictly positive terms defined by 
\[
u= \frac{1}{t} \mathrm{Tr} {\bf D}  {\bf U}_1^{-1} \left( \sigma^{2} {\bf I} +  {\bf B} ({\bf I} + 
h(\tilde{\kappa}_2) \tilde{{\bf D}})^{-1}  ({\bf I} + h(\tilde{\kappa}_1) \tilde{{\bf D}})^{-1} {\bf B}^{H} 
\right) {\bf U}_2^{-1}  
\]
and 
\[
v = \frac{1}{t} \mathrm{Tr} {\bf D} {\bf U}_1^{-1}  {\bf D} {\bf U}_2^{-1}\ .
\]

As $u (h(\tilde{\kappa}_2) -h(\tilde{\kappa}_1)) + v  (\tilde{\kappa}_2 h(\tilde{\kappa}_2) - \tilde{\kappa}_1 h(\tilde{\kappa}_1)) = 0$, 
$(h(\tilde{\kappa}_2) -h(\tilde{\kappa}_1)) < 0$ implies that $\tilde{\kappa}_2 h(\tilde{\kappa}_2) - \tilde{\kappa}_1 h(\tilde{\kappa}_1) > 0$. 
Hence, $\tilde{\kappa} h(\tilde{\kappa})$ is a strictly increasing function as expected. 

From this, it follows that function $\tilde{\kappa} \rightarrow
\tilde{g}(h(\tilde{\kappa}), \tilde{\kappa})$ is strictly decreasing.
This function converges to $+\infty$ if $\tilde{\kappa} \rightarrow 0$
and to $0$ if $\tilde{\kappa} \rightarrow +\infty$. Therefore, the
equation
\[
\tilde{\kappa} \rightarrow \tilde{g}(h(\tilde{\kappa}), \tilde{\kappa}) = 1
\] 
has a unique strictly positive solution $\tilde{\beta}$. If $\beta = h(\tilde{\beta})$, it is clear that 
$g(\beta, \tilde{\beta}) = 1$ and $\tilde{g}(\beta, \tilde{\beta}) = 1$. Therefore, we have shown 
that $(\beta, \tilde{\beta})$ is the unique solution of (\ref{eq:canonique-preparatoire}) satisfying $\beta > 0$ and 
$\tilde{\beta} > 0$.

\section{Proof of Theorem \ref{theo:convergence-rate}}
\label{proof-theo-convergence-rate}
This section is organized as follows. We first recall in subsection
\ref{subsec:tools} some useful mathematical tools. In subsection
\ref{subsec:immediat}, we establish
(\ref{eq:vitesse-convergence-normalise}). In \ref{subsec:sketch}, we
prove (\ref{eq:vitesse-convergence-biais}) and
(\ref{eq:vitesse-rapide}).   \\

We shall use the following notations. If $u$ is a random variable, the
zero mean random variable $u - \mathbb{E}(u)$ is denoted by
$\overcirc{u}$. If $z = x + \mathbf{i} y$ is a complex number, the differential
operators $\frac{\partial}{\partial z}$ and $\frac{\partial}{\partial
  \bar{z}}$ are defined respectively by $\frac 12 \left( \frac{\partial}{\partial
  x} - \mathbf{i} \frac{\partial}{\partial y}\right)$ and 
  $\frac 12 \left( \frac{\partial}{\partial x} + \mathbf{i} \frac{\partial}{\partial y}\right)$.
Finally, if $\bs{\Sigma}, {\bf B}, {\bf Y}$ are given matrices, we denote respectively by 
$\bs{{\bf \xi}}_j, {\bf b}_j, {\bf y}_j$ their columns.
\subsection{Mathematical tools.}
\label{subsec:tools}
\subsubsection{The Poincaré-Nash inequality} (see e.g. \cite{Che-82}, \cite{Hou-Per-Sur-98}).
Let $ {\bf x} = [ x_1, \ldots, x_M ]^T$ be a complex Gaussian random
vector whose law is given by $\EE [ {\bf x} ] = {\bf 0}$,
$\EE [ {\bf x} {\bf x}^T ] = {\bf 0}$, and
$\EE [ {\bf x} {\bf x}^* ] = {\bs \Xi}$. Let 
$\Phi = \Phi(x_1, \ldots, x_M, \bar{x}_1, \ldots, \bar{x}_M)$ be a $C^1$ complex function polynomially bounded together
with its partial derivatives. Then the following inequality holds true: 
$$
\mathrm{Var}(\Phi({\bf x})) \leq 
\EE \left[ \nabla_z \Phi({\bf x})^T \ {\bs \Xi} \ 
\overline{\nabla_z \Phi({\bf x})} \right] + \EE \left[ \left(\nabla_{\bar{z}} \Phi({\bf x}) \right)^H \ {\bs \Xi} \ 
\nabla_{\bar{z}} \Phi({\bf x}) \right] \ ,
$$ 
where $\nabla_z \Phi
= [ \partial\Phi / \partial z_1, \ldots, 
\partial\Phi / \partial z_M ]^T$ and  $\nabla_{\bar{z}} \Phi
= [ \partial\Phi / \partial \bar{z}_1, \ldots, 
\partial\Phi / \partial \bar{z}_M ]^T$. 
 
Let ${\bf Y}$ be the $r× t$ matrix ${\bf Y}=\frac 1{\sqrt{t}}
{\bf D}^{\frac 12} {\bf X} \tilde{\bf D}^{\frac 12}$, where ${\bf X}$
has i.i.d. ${\mathcal C}N(0,1)$ entries and consider the stacked 
$rt × 1$ vector
${\bf x} = [ Y_{11}, \ldots, Y_{rt} ]^T$. In this case, Poincaré-Nash inequality
writes:
\begin{equation}
\label{eq:np-our-model} 
\mathrm{Var} \left( \Phi({\bf Y}) \right) \leq 
\frac{1}{t} \sum_{i=1}^r \sum_{j=1}^t d_i \tilde d_j 
\EE \left[ \left| \frac{\partial \Phi({\bf Y})}
{\partial Y_{i,j}} \right|^2 +  \left| \frac{\partial \Phi({\bf Y})}
{\partial \overline{Y}_{i,j}} \right|^2 \right]\ . 
\end{equation}

\subsubsection{The differentiation formula for functions of Gaussian 
random vectors} 

With ${\bf x}$ and $\Phi$ given as above, we have the following 
\begin{equation}
\label{eq:integ-parts}
\EE \left[ x_p \Phi({\bf x}) \right] = 
\sum_{m=1}^M \left[ {\bs \Xi} \right]_{pm} 
\EE \left[ \frac{\partial \Phi({\bf x})}{\partial \bar{x}_{m}} \right]\ .
\end{equation}

This formula relies on an integration by parts, and is thus referred to as the Integration by parts 
formula for Gaussian vectors. It is widely used in Mathematical Physics (\cite{Gli-Jaf-87}) and has been used 
in Random Matrix Theory in \cite{Kho-Pas-93} and \cite{Pas-Kho-Vas-95}. 

If ${\bf x}$ coincides with the $rt × 1$ vector ${\bf x} = [ Y_{11}, \ldots, Y_{rt} ]^T$, 
relation (\ref{eq:integ-parts}) becomes 
\begin{equation}
\label{eq:integ-parts-our-model}
\EE \left[ Y_{pq}  \Phi({\bf Y}) \right] = 
\frac{d_p \tilde{d}_q  }{t}
\EE \left[ \frac{\partial \Phi({\bf Y})}{\partial \overline{Y}_{pq}} \right] \ . 
\end{equation} 
Replacing matrix ${\bf Y}$ by matrix $\bar{{\bf Y}}$ also provides
\begin{equation}
\label{eq:integ-conj-parts-our-model}
\EE \left[ \overline{Y}_{pq}  \Phi({\bf Y}) \right] = 
\frac{d_p \tilde{d}_q}{t}  
\EE \left[ \frac{\partial \Phi({\bf Y})}{\partial Y_{pq}} \right] \ . 
\end{equation} 
\subsubsection{Some useful differentiation formulas}
The following partial derivatives $\frac{\partial ({\bf
    S}_{pq})}{\partial Y_{ij}}$ and $\frac{\partial {\bf
    S}_{pq}}{\partial \overline{Y}_{ij}}$ for each $p,q \in
\{1, \ldots, r \}$ and $1 \leq i \leq r, 1 \leq j \leq t$ will
be of use in the sequel. Straightforward computations yield:
\begin{equation}
\label{eq:derivees-S}
\left\{ 
\begin{array}{lcr} 
\frac{\partial S_{pq}}{\partial Y_{ij}} & = & - S_{p,i} \left( {\bs \xi}_j^{H} {\bf S} \right)_{q}  \\
\frac{\partial S_{pq}}{\partial \bar{Y_{ij}}} & = & - S_{i,q} \left( {\bf S} {\bs \xi} \right)_{p} 
\end{array}\right.\ .
\end{equation}

\subsection{Proof of (\ref{eq:vitesse-convergence-normalise})}
\label{subsec:immediat}
We just prove that the variance of $\frac{1}{t} \mbox{Tr}({\bf M} {\bf
  S})$ is a $O({t^{-2}})$ term. For this, we note that the
random variable $\frac{1}{t} \mbox{Tr}({\bf M} {\bf S})$ can be
interpreted as a function $\Phi({\bf Y})$ of the entries of matrix
${\bf Y}$, and use the Poincaré-Nash inequality
(\ref{eq:np-our-model}) to $\Phi({\bf Y})$. Function $\Phi({\bf Y})$ is equal
to
\[
\Phi({\bf Y}) = \frac{1}{t} \sum_{p,q} M_{q,p} S_{p,q}\ .
\]
Therefore, the partial derivative of $\Phi({\bf Y})$ with respect to $Y_{ij}$ is given by 
$
\frac{\partial \Phi({\bf Y}) }{\partial Y_{ij}} =  \frac{1}{t} \sum_{p,q} {\bf M}_{q,p} \frac{\partial S_{pq}}{\partial Y_{ij}}
$
which, by (\ref{eq:derivees-S}), coincides with
\[
\frac{\partial \Phi({\bf Y}) }{\partial Y_{ij}} = - \frac{1}{t}
\sum_{p,q} {\bf M}_{q,p} {\bf S}_{p,i} ({\bs \xi}_j^{H} {\bf S})_q =
-\frac{1}{t} \left( {\bs \xi}_j^{H} {\bf S} {\bf M} {\bf S} \right)_i\ .
\]
As $d_i \leq d_{\max}$ and $\tilde{d}_j \leq \tilde{d}_{\max}$, it is clear that  
\[
\sum_{i=1}^{r} \sum_{j=1}^{t} d_i \tilde{d}_j \mathbb{E} \left|
  \frac{\partial \Phi({\bf Y}) }{\partial Y_{ij}} \right|^{2} \leq
d_{\max} \tilde{d}_{\max} \sum_{i=1}^{r} \sum_{j=1}^{t} \mathbb{E}
\left| \frac{\partial \Phi({\bf Y}) }{\partial Y_{ij}} \right|^{2}\ .
\]
It is easily seen that 
\[
\sum_{i=1}^{r} \mathbb{E} \left| \frac{\partial \Phi({\bf Y})
  }{\partial Y_{ij}} \right|^{2} = \frac{1}{t^{2}} \mathbb{E} \left(
  {\bs \xi}_j^{H} {\bf S} {\bf M} {\bf S}^{2} {\bf M}^{H} {\bf S} {\bs
    \xi}_j^{H} \right)\ .
\]
As $\| {\bf S} \| \leq \frac{1}{\sigma^{2}}$ and $\sup_t \|{\bf M} \| < \infty$,  ${\bs \xi}_j^{H} {\bf S} {\bf M} {\bf S}^{2} {\bf M}^{H} {\bf S} {\bs \xi}_j^{H}$ 
is less than $ \frac{1}{\sigma^{8}} \sup_t \|{\bf M} \|^{2} \; \|{\bs \xi}_j \|^{2}$. Moreover, 
$\mathbb{E} \|{\bs \xi}_j \|^{2}$ coincides with $\|{\bf b}_j \|^{2} + \frac{1}{t} \tilde{d}_j \sum_{i=1}^{r} d_i$, 
which is itself less than $b_{\max}^{2} + d_{\max} \tilde{d}_{\max} \frac{r}{t}$, a uniformly bounded term.    
Therefore, 
$
\sum_{i=1}^{r} \mathbb{E} \left| \frac{\partial \Phi({\bf Y}) }{\partial Y_{ij}} \right|^{2} 
$
is a $O({t^{-2}})$ term. This proves that
\[
\frac{1}{t} \sum_{i=1}^{r} \sum_{j=1}^{t} d_i \tilde{d}_j \mathbb{E}
\left| \frac{\partial \Phi({\bf Y}) }{\partial Y_{ij}} \right|^{2} =
O\left(\frac{1}{t^{2}}\right)\ .
\]
It can be shown similarly that 
$
{t}^{-1} \sum_{i=1}^{r} \sum_{j=1}^{t} d_i \tilde{d}_j \mathbb{E}
\left| \frac{\partial \Phi({\bf Y}) }{\partial \overline{Y_{ij}}}
\right|^{2} = O\left({t^{-2}}\right)\ .
$
The conclusion follows from Poincaré-Nash inequality (\ref{eq:np-our-model}). 

\subsection{Proof of (\ref{eq:vitesse-convergence-biais}) and (\ref{eq:vitesse-rapide}).}
\label{subsec:sketch}
As we shall see, proofs of (\ref{eq:vitesse-convergence-biais}) and
(\ref{eq:vitesse-rapide}) are demanding.  We first introduce the
following notations: Define scalar parameters $\eta(\sigma^{2}),
\alpha(\sigma^{2}), \tilde{\alpha}(\sigma^{2})$ as
\begin{equation}
\label{eq:def-eta-alpha}
\begin{array}{ccc}
\eta(\sigma^{2}) & = & \frac{1}{t} \mbox{Tr} \left({\bf D} {\bf S}(\sigma^{2}) \right) \\
\alpha(\sigma^{2}) & = & \EE \left[ \frac{1}{t} \mbox{Tr} \left({\bf D}{\bf S}(\sigma^{2}) \right) \right]\\
\tilde{\alpha}(\sigma^{2}) & = &  \EE \left[ \frac{1}{t} \mbox{Tr} \left(\tilde{{\bf D}} \tilde{{\bf S}}(\sigma^{2}) \right) \right]
\end{array}
\end{equation}
and matrices ${\bf R}(\sigma^{2}), \tilde{{\bf R}}(\sigma^{2})$ as
\begin{equation}
\label{eq:def-R-tildeR}
\begin{array}{ccc}
  {\bf R}(\sigma^{2}) &=  & \left[ \sigma^{2} \left( {\bf I} + \tilde{\alpha} {\bf D} \right) + {\bf B} \left( {\bf I} + \alpha \tilde{{\bf D}} \right)^{-1} {\bf B}^{H} \right]^{-1} \\
  \tilde{{\bf R}}(\sigma^{2}) &=  & \left[ \sigma^{2} \left( {\bf I} + \alpha \tilde{\bf D} \right) + {\bf B}^H \left( {\bf I} + \tilde{\alpha} {{\bf D}} \right)^{-1} {\bf B} \right]^{-1}
\end{array}\ .
\end{equation} 
We note that, as $\alpha(\sigma^{2}) \geq 0$ and $\tilde{\alpha}(\sigma^{2}) \geq 0$, 
then
\begin{equation}
\label{eq:bornesR}
0 < {\bf R}(\sigma^{2}) < \frac{{\bf I}_r}{\sigma^{2}}, \; 0 < \tilde{{\bf R}}(\sigma^{2}) < \frac{{\bf I}_t}{\sigma^{2}}
\end{equation}
It is difficult to study directly the term $\frac{1}{r} \mathrm{Tr}
{\bf M} (\mathbb{E}({\bf S}) - {\bf T})$. In some sense, matrix ${\bf
  R}$ can be seen as an intermediate quantity between $\mathbb{E}({\bf S})$ and
${\bf T}$. Thus the proof consists into two steps: 1) for each uniformly bounded 
matrix ${\bf M}$, we first prove that
$\frac{1}{r}\mathrm{Tr} {\bf M} (\mathbb{E}({\bf S}) - {\bf R})$ and  
$\frac{1}{r} \mathrm{Tr} {\bf M} ({\bf R} - {\bf T})$ converge to $0$ 
as $t \to \infty$; 2) we then refine the previous result and establish in fact that 
$\frac{1}{r}\mathrm{Tr} {\bf M} (\mathbb{E}({\bf S}) - {\bf R})$ and  
$\frac{1}{r} \mathrm{Tr} {\bf M} ({\bf R} - {\bf T})$ are $O({t^{-2}})$ terms. 
This, of course, imply (\ref{eq:vitesse-convergence-biais}). Eq. (\ref{eq:vitesse-rapide}) eventually follows 
from Eq. (\ref{eq:vitesse-convergence-biais}), the integral representation
\begin{equation}
\label{eq:expre-ecartEMI}
\overline{J}(\sigma^{2}) - J(\sigma^{2}) = \int_{\sigma^{2}}^{\infty} \mathrm{Tr} \left( \mathrm{E}({\bf S}(\omega)) - {\bf T}(\omega) \right) \, d \, \omega,
\end{equation}
which follows from \eqref{eq:lien-trace-log} and \eqref{eq:defJbarre}, 
as well as a dominated convergence argument that is omitted. 

\subsubsection{First step: Convergence of $\frac{1}{r}\mathrm{Tr} {\bf M} (\mathbb{E}({\bf S}) - {\bf R})$ and  
$\frac{1}{r} \mathrm{Tr} {\bf M} ({\bf R} - {\bf T})$ to zero}
The first step consists in showing the following Proposition.
\begin{proposition}
\label{prop:E(S)-R}
For each deterministic $r × r$ matrix ${\bf M}$, uniformly
bounded (for the spectral norm) as $t \to \infty$, we have:
\begin{equation}
\label{eq:E(S)-R-zero}
\lim_{t \to +\infty} \frac{1}{t} \mathrm{Tr}\left[ {\bf M} \left( \mathbb{E}({\bf S}) - {\bf R} \right) \right] = 0
\end{equation}
\begin{equation}
\label{eq:R-T-zero}
\lim_{t \to +\infty} \frac{1}{t} \mathrm{Tr}\left[ {\bf M} \left( {\bf R}) - {\bf T} \right) \right] = 0
\end{equation}
\end{proposition}

\begin{proof}
We first prove (\ref{eq:E(S)-R-zero}). For this, we state the following  useful Lemma. 
\begin{lemma}
\label{le:eval-variances-utiles}
Let ${\bf P}, {\bf P}_1$ and ${\bf P}_2$ be deterministic $r × t,
t × t, t × r$ matrices respectively, uniformly bounded with
respect to the spectral norm as $t\to \infty$. Consider the
following functions of ${\bf Y}$.
$$
\bs{\Phi}({\bf Y})= \frac{1}{t} \mbox{Tr} \left[ {\bf S} {\bf P} {\bs{\Sigma}}^{H} \right],\quad
\bs{\Psi}({\bf Y}) = \frac{1}{t} \mbox{Tr} \left[ {\bf S}  {\bs{\Sigma}} {\bf P}_1 {\bs{\Sigma}}^{H} {\bf P}_2 \right],\quad
\bs{\Psi}^{'}({\bf Y}) = \frac{1}{t} \mbox{Tr} \left[ {\bf S}  {\bs{\Sigma}} {\bf P}_1 {\bf{Y}}^{H} {\bf P}_2 \right] \ .
$$
Then, the following estimates hold true:
$$
\mbox{Var}(\bs{\Phi}) = O\left(\frac{1}{t^{2}}\right),\quad 
\mbox{Var}(\bs{\Psi}) = O\left(\frac{1}{t^{2}}\right),\quad 
\mbox{Var}(\bs{\Psi}^{'}) = O\left(\frac{1}{t^{2}}\right)\ .
$$
\end{lemma}
The proof, based on the Poincaré-Nash inequality (\ref{eq:np-our-model}), is omitted.

In order to use the Integration by parts formula 
(\ref{eq:integ-parts-our-model}), notice that
\begin{equation}
\label{eq:resolvent}
\sigma^{2} {\bf S}(\sigma^{2}) +  {\bf S}(\sigma^{2}) {\bs \Sigma}  {\bs \Sigma}^{H} = {\bf I}\ .
\end{equation}           
Taking the mathematical expectation, we have for each $p,q \in \{1, \ldots, r \}$:
\begin{equation}
\label{eq:resolvent_(p,q)}
\sigma^{2} \mathbb{E}(S_{pq}) +  \mathbb{E} \left[( {\bf S} {\bs \Sigma}  {\bs \Sigma}^{H})_{pq} \right] = \delta(p-q)\ .
\end{equation}
A convenient use of the Integration by parts formula allows to express  $\mathbb{E} \left[ ({\bf S} {\bs \Sigma}  {\bs \Sigma}^{H})_{pq} \right]$
in terms of the entries of $\mathbb{E}({\bf S})$. To see this, note that 
\[
\mathbb{E} \left[ ({\bf S} {\bs \Sigma} {\bs \Sigma}^{H})_{pq} \right]
= \sum_{j=1}^{t} \sum_{i=1}^{r} \mathbb{E} (S_{pi} \Sigma_{ij}
\overline{\Sigma_{qj}})\ .
\]
For each $i$,  $\mathbb{E} (S_{pi} \Sigma_{ij} \overline{\Sigma_{qj}})$ can be written as
\[
\mathbb{E} (S_{pi} \Sigma_{ij} \overline{\Sigma_{qj}}) = \mathbb{E}
(S_{pi}) B_{ij} \overline{B_{qj}} + \mathbb{E} \left( S_{pi}
  \overline{Y_{qj}} \right) B_{ij} + \mathbb{E} \left( S_{pi} Y_{ij}
  \overline{\Sigma_{qj}} \right)\ .
\]
Using (\ref{eq:integ-parts-our-model}) with function $\Phi({\bf Y}) =  S_{pi}  \overline{\Sigma_{qj}}$ and (\ref{eq:integ-conj-parts-our-model})
with $\Phi({\bf Y}) =  S_{pi}$, and summing over index $i$ yields:
\begin{equation}
\label{eq:point-de-depart}
\mathbb{E} \left[ ({\bf S} {\bs \xi}_j)_p  \overline{\Sigma_{q,j}} \right] = \frac{d_q \tilde{d}_j}{t} \mathbb{E}(S_{pq}) - \tilde{d}_j \EE \left[ \eta  ({\bf S} {\bs \xi}_j)_p  \overline{\Sigma_{q,j}} \right] -  \frac{d_q \tilde{d}_j}{t} \EE \left[ S_{pq} {\bs \xi}_j^{H} {\bf S} {\bf b}_j \right] + \EE \left[ ({\bf S} {\bf b}_j)_p \right] \overline{B_{qj}}\ .
\end{equation}
Eq. (\ref{eq:vitesse-convergence-normalise}) for ${\bf M} = {\bf D}$
implies that $\mathrm{Var}(\eta) = O({t^{-2}})$, or
equivalently that $\mathbb{E} (\overcirc{\eta}^{2}) =
O({t^{-2}})$.
We now complete proof of (\ref{eq:E(S)-R-zero}). We take Eq.
(\ref{eq:point-de-depart}) as a starting point, and write $\eta$ as
$\eta = \mathbb{E}(\eta) + \overcirc{\eta} = \alpha +
\overcirc{\eta}$.  Therefore,
\[
\mathbb{E}\left[\eta \, ({\bf S} \bs{\xi}_j)_p \, \overline{{\Sigma}_{q,j}} \right] = \alpha \mathbb{E}\left[ ({\bf S} \bs{\xi}_j)_p \, \overline{{\Sigma}_{q,j}} \right] + \mathbb{E}\left[\overcirc{\eta} \, ({\bf S} \bs{\xi}_j)_p \, \overline{{\Sigma}_{q,j}} \right]\ .
\]
Plugging this relation into (\ref{eq:point-de-depart}), and solving
w.r.t.  $\mathbb{E}\left[ ({\bf S} \bs{\xi}_j)_p \,
  \overline{{\Sigma}_q,j} \right]$ yields
\begin{eqnarray*}
\mathbb{E} \left[ ({\bf S} {\bs \xi}_j)_p  \overline{\Sigma_{q,j}} \right] &  = & \frac{1}{t} \frac{d_q \tilde{d}_j}{1 + \alpha \tilde{d}_j} \EE ( S_{pq}) + 
\frac{1}{1 + \alpha \tilde{d}_j}  \EE \left[ ({\bf S} {\bf b}_j)_p \right] \overline{B_{qj}}  \\
    &  &   - \, \frac{1}{t} \frac{d_q \tilde{d}_j}{1 + \alpha \tilde{d}_j}  \EE \left[ S_{pq} {\bs \xi}_j^{H} {\bf S} {\bf b}_j \right] - \frac{\tilde{d}_j}{1 + \alpha \tilde{d}_j} \mathbb{E}\left[\overcirc{\eta} \, ({\bf S} \bs{\xi}_j)_p \, \overline{{\Sigma}_{q,j}} \right]\ .
\end{eqnarray*}
Writing  ${\bs \xi}_j = {\bf b}_j + {\bf y}_j$, and 
summing over $j$ provides the following  expression of $\mathbb{E} \left[ ({\bf S} {\bs \Sigma}  {\bs \Sigma}^{H})_{pq} \right]$:
\begin{multline}
\label{eq:intermediaire-1-exact}
\mathbb{E} \left[ ({\bf S} {\bs \Sigma} {\bs \Sigma}^{H})_{pq} \right]
\quad =\quad  d_q \, \frac{1}{t} \mathrm{Tr} \left[ \tilde{{\bf D}} ({\bf I} +
  \alpha \tilde{{\bf D}})^{-1} \right]
\mathbb{E}(S_{pq}) \\
+ \EE \left[ \left( {\bf S} {\bf B} ({\bf I} + \alpha \tilde{{\bf
        D}})^{-1} {\bf B}^{H} \right)_{pq} \right] - d_q \, \EE \left[
  S_{pq} \frac{1}{t} \mathrm{Tr} \left( {\bf S} {\bf B} \tilde{{\bf
        D}} ({\bf I} + \alpha \tilde{{\bf D}})^{-1} {\bf B}^{H}
  \right) \right]
\\
- d_q  \, \EE \left[ S_{pq} \frac{1}{t} \mathrm{Tr} \left( {\bf S} {\bf B} \tilde{{\bf D}} ({\bf I} + \alpha \tilde{{\bf D}})^{-1} {\bf Y}^{H} \right) \right] 
- \mathbb{E} \left[ \overcirc{\eta} \left( {\bf S} {\bs \Sigma}
    \tilde{{\bf D}} ({\bf I} + \alpha \tilde{{\bf D}})^{-1} {\bs
      \Sigma}^{H} \right)_{p,q} \right]\ .
\end{multline}
The resolvent identity (\ref{eq:resolvent}) thus implies that 
\begin{multline}
\label{eq:intermediaire-2-exact}
\delta(p-q) \quad = \quad \sigma^{2} \EE (S_{pq})
+   \frac{d_q}{t} \mathrm{Tr} \left[ \tilde{{\bf D}} ({\bf I} + \alpha \tilde{{\bf D}})^{-1} \right] \mathbb{E}(S_{pq})\\
+ \EE \left[ \left( {\bf S} {\bf B} ({\bf I} + \alpha \tilde{{\bf
        D}})^{-1} {\bf B}^{H} \right)_{pq} \right]
- d_q \, \EE \left[ S_{pq} \frac{1}{t} \mathrm{Tr} \left( {\bf S} {\bf B} \tilde{{\bf D}} ({\bf I} + \alpha \tilde{{\bf D}})^{-1} {\bf B}^{H} \right) \right] \\
- d_q  \, \EE \left[ S_{pq} \frac{1}{t} \mathrm{Tr} \left( {\bf S} {\bf B} \tilde{{\bf D}} ({\bf I} + \alpha \tilde{{\bf D}})^{-1} {\bf Y}^{H} \right) \right] 
- \mathbb{E} \left[ \overcirc{\eta} \left( {\bf S} {\bs \Sigma}
    \tilde{{\bf D}} ({\bf I} + \alpha \tilde{{\bf D}})^{-1} {\bs
      \Sigma}^{H} \right)_{p,q} \right]\ .
\end{multline}

In order to simplify the notations, we define $\rho_1$ and $\rho_2$ by 
$$
\rho_1 \ = \  \frac{1}{t} \mathrm{Tr} \left( {\bf S} {\bf B} \tilde{{\bf D}} ({\bf I} + \alpha \tilde{{\bf D}})^{-1} {\bf B}^{H} \right)  \quad\textrm{and}\quad
\rho_2 \ =\     \frac{1}{t} \mathrm{Tr} \left( {\bf S} {\bf B} \tilde{{\bf D}} ({\bf I} + \alpha \tilde{{\bf D}})^{-1} {\bf Y}^{H} \right)\ .
$$
For $i=1,2$, we write $\EE  (S_{pq} \rho_i)$ as
\[
\EE  (S_{pq} \rho_i) = \EE (S_{pq}) \, \EE (\rho_i) + \EE \left( \overcirc{S_{pq}} \; \overcirc{\rho_i} \right)\ .
\]
Thus, (\ref{eq:intermediaire-2-exact}) can be written as 
\begin{multline}
\label{eq:intermediaire-3-exact}
\delta(p-q) \quad =  \quad  \sigma^{2} \EE (S_{pq}) +  d_q \, \frac{1}{t} \mathrm{Tr} \left[ \tilde{{\bf D}} ({\bf I} + \alpha \tilde{{\bf D}})^{-1} \right] \mathbb{E}(S_{pq}) \\
 +   \left( \EE({\bf S}) {\bf B}  ({\bf I} + \alpha \tilde{{\bf D}})^{-1} {\bf B}^{H} \right)_{pq}  
 - d_q \, \EE (S_{pq})  \frac{1}{t} \mathrm{Tr} \left( \EE({\bf S}) {\bf B} \tilde{{\bf D}} ({\bf I} + \alpha \tilde{{\bf D}})^{-1} {\bf B}^{H} \right)  \\
 - d_q  \, \EE (S_{pq}) \EE \left[\frac{1}{t} \mathrm{Tr} \left( {\bf S} {\bf B} \tilde{{\bf D}} ({\bf I} + \alpha \tilde{{\bf D}})^{-1} {\bf Y}^{H} \right) \right] 
  -d_q \, \EE \left(\overcirc{S_{pq}} \; \overcirc{\rho_1} \right) 
- d_q \, \EE \left( \overcirc{S_{pq}} \; \overcirc{\rho_2} \right) \\
 - \mathbb{E} \left[ \overcirc{\eta} \left( {\bf S} {\bs \Sigma} \tilde{{\bf D}} ({\bf I} + \alpha \tilde{{\bf D}})^{-1} {\bs \Sigma}^{H} \right)_{p,q} \right] \ .
\end{multline}
We now establish the following lemma.
\begin{lemma}
\begin{eqnarray}
\EE \rho_2 &=& \EE \left[\frac{1}{t} \mathrm{Tr} \left( {\bf S} {\bf B} \tilde{{\bf D}} ({\bf I} + \alpha \tilde{{\bf D}})^{-1} {\bf Y}^{H} \right) \right] \nonumber \\
&=& - \alpha \, 
\frac{1}{t} \mathrm{Tr} \left( \mathbb{E}({\bf S}) {\bf B} \tilde{{\bf D}}^{2} ({\bf I} + \alpha \tilde{{\bf D}})^{-2} {\bf B}^{H} \right)
- \EE \left( \overcirc{\eta} \; \overcirc{\rho_3} \right) \label{eq:utile-exact}\ ,
\end{eqnarray}
 where $\rho_3$ is defined by
\[
\rho_3 = \frac{1}{t} \mathrm{Tr} \left( {\bf S} {\bf B} \tilde{{\bf D}}^{2} ({\bf I} + \alpha \tilde{{\bf D}})^{-2} {\bs \Sigma}^{H} \right)\ .
\]
\end{lemma}
\begin{proof} We express $\EE (\rho_2)$ as 
\begin{equation}
\label{eq:utile-exact-1}
\begin{array}{ccc}
\EE (\rho_2) & = & \frac{1}{t} \sum_{j=1}^{t} \frac{\tilde{d}_j}{1 + \alpha \tilde{d}_j} \; \EE ({\bf y}_j^{H} {\bf S} {\bf b}_j) \\
\    &  = &  \frac{1}{t} \sum_{j=1}^{t} \frac{\tilde{d}_j}{1 + \alpha \tilde{d}_j} \, \sum_{i=1}^{r} \EE \left(({\bf S} {\bf b}_j)_i\overline{Y_{ij}}\right)
\end{array}
\end{equation}
and evaluate $ \EE \left(({\bf S} {\bf b}_j)_i\overline{Y_{ij}}\right)$ using formula (\ref{eq:integ-conj-parts-our-model}) for 
$\Phi({\bf Y}) = ({\bf S} {\bf b}_j)_i$. This gives
\[
 \EE \left(({\bf S} {\bf b}_j)_i\overline{Y_{ij}}\right) = \frac{1}{t} d_i \tilde{d}_j \sum_{k=1}^{r} \EE \left( \frac{\partial S_{ik}}{\partial Y_{ij}} \right) B_{kj}\ .
\]
By (\ref{eq:derivees-S}), 
\[
\EE \left( \frac{\partial S_{ik}}{\partial Y_{ij}} \right) =  - \EE \left(S_{ii} ({\bf b}_j^{H} {\bf S})_{k} \right)
- \EE  \left(S_{ii} ({\bf y}_j^{H} {\bf S})_{k} \right)\ .
\]
Therefore, 
\[
\EE \left( {\bf y}_j^{H} {\bf S} {\bf b}_j \right) = -\tilde{d}_j \EE \left( \eta \, {\bf b}_j^{H} {\bf S} {\bf b}_j \right) -  \tilde{d}_j \EE \left( \eta \, {\bf y}_j^{H} {\bf S} {\bf b}_j \right)\ .
\]
Writing again $\eta = \EE (\eta) + \overcirc{\eta} = \alpha + \overcirc{\eta}$, we get that 
\begin{equation}
\begin{array}{ccc}
\EE \left( {\bf y}_j^{H} {\bf S} {\bf b}_j \right) &  =  & -\alpha \tilde{d}_j \EE \left( {\bf b}_j^{H} {\bf S} {\bf b}_j \right) - \alpha  \tilde{d}_j \EE \left( {\bf y}_j^{H} {\bf S} {\bf b}_j \right) \\
  &  &  -\tilde{d}_j \EE \left( \overcirc{\eta} \, {\bf b}_j^{H} {\bf S} {\bf b}_j \right) -  \tilde{d}_j \EE \left( \overcirc{\eta} \, {\bf y}_j^{H} {\bf S} {\bf b}_j \right)\ .
\end{array}
\end{equation}
Solving this equation w.r.t. $\EE \left( {\bf y}_j^{H} {\bf S} {\bf b}_j \right)$ yields
\begin{equation}
\EE \left( {\bf y}_j^{H} {\bf S} {\bf b}_j \right) =   -\frac{\alpha \tilde{d}_j}{1 + \alpha \tilde{d}_j}  \EE \left( {\bf b}_j^{H} {\bf S} {\bf b}_j \right) 
 -\frac{\tilde{d}_j}{1 + \alpha \tilde{d}_j}  \EE \left( \overcirc{\eta} \, {\bf b}_j^{H} {\bf S} {\bf b}_j \right) -  \frac{\tilde{d}_j}{1+\alpha \tilde{d}_j}  \EE \left( \overcirc{\eta} \, {\bf y}_j^{H} {\bf S} {\bf b}_j \right)
\end{equation}
or equivalently
\begin{equation}
\label{eq:utile-exact-2}
\EE \left( {\bf y}_j^{H} {\bf S} {\bf b}_j \right) =   -\frac{\alpha \tilde{d}_j}{1 + \alpha \tilde{d}_j}  \EE \left( {\bf b}_j^{H} {\bf S} {\bf b}_j \right) 
 -\frac{\tilde{d}_j}{1 + \alpha \tilde{d}_j}  \EE \left( \overcirc{\eta} \, {\bs \xi}_j^{H} {\bf S} {\bf b}_j \right)\ .
\end{equation}
Eq. (\ref{eq:utile-exact}) immediately follows from (\ref{eq:utile-exact-1}), (\ref{eq:utile-exact-2}), 
and the relation $\mathbb{E}(  \overcirc{\eta} \, \rho_3) = \mathbb{E}(  \overcirc{\eta} \, \overcirc{\rho_3})$. 
\end{proof}

Plugging (\ref{eq:utile-exact}) into (\ref{eq:intermediaire-3-exact}) yields
\begin{multline}
\label{eq:intermediaire-4-exact}
\delta(p-q) + \Delta_{pq}  \\
=    \EE (S_{pq}) \left[ \sigma^{2} + d_q \left( \frac{1}{t} \mathrm{Tr} 
\tilde{{\bf D}} ({\bf I} + \alpha \tilde{{\bf D}})^{-1} - \EE ( \rho_1) + \alpha  \frac{1}{t} \mathrm{Tr} 
\EE({\bf S}) {\bf B} \tilde{{\bf D}}^{2} ({\bf I} + \alpha \tilde{{\bf D}})^{-2} {\bf B}^{H} \right) \right] \\
   + \left[ \EE({\bf S}) {\bf B} ({\bf I} + \alpha \tilde{{\bf D}})^{-1} {\bf B}^{H} \right]_{pq}
\end{multline}
where ${\bs \Delta}$ is the $r × r$ matrix defined by
\[
\Delta_{pq} = \EE \left[ \overcirc{\eta} \left( {\bf S} {\bs \Sigma} \tilde{{\bf D}} ({\bf I} + \alpha \tilde{{\bf D}})^{-1} {\bs \Sigma}^{H}\right)_{pq} \right] + d_q \EE \left( \overcirc{S}_{pq}(\overcirc{\rho_1} + \overcirc{\rho_2}) \right) - d_q \EE(S_{pq}) \, \EE \left(  \overcirc{\eta} \, \overcirc{\rho_3} \right) 
\]
for each $p,q$ or equivalently by
\[
{\bs \Delta} =  \EE \left[ \overcirc{\eta} \left( {\bf S} {\bs \Sigma} \tilde{{\bf D}} ({\bf I} + \alpha \tilde{{\bf D}})^{-1} {\bs \Sigma}^{H}\right) \right] + \EE \left((\overcirc{\rho_1} + \overcirc{\rho_2}) \, \overcirc{{\bf S}} \right)
\, {\bf D} - \EE \left(  \overcirc{\eta} \, \overcirc{\rho_3} \right) \, \EE({\bf S}) \, {\bf D}\ .
\]
Using the relation $
\alpha \tilde{{\bf D}} ({\bf I} + \alpha \tilde{{\bf D}})^{-1} = {\bf I} -  ({\bf I} + \alpha \tilde{{\bf D}})^{-1}$,
we obtain that 
\begin{eqnarray}
\lefteqn{ 
\alpha  \frac{1}{t} \mathrm{Tr} \left( 
\EE({\bf S}) {\bf B} \tilde{{\bf D}}^{2} ({\bf I} + \alpha \tilde{{\bf D}})^{-2} {\bf B}^{H} \right) }\nonumber \\
& = &  \frac{1}{t} \mathrm{Tr} \left(
\EE({\bf S}) {\bf B} \tilde{{\bf D}} ({\bf I} + \alpha \tilde{{\bf D}})^{-1} {\bf B}^{H} \right)
- \frac{1}{t} \mathrm{Tr} \left(
\EE({\bf S}) {\bf B} \tilde{{\bf D}} ({\bf I} + \alpha \tilde{{\bf D}})^{-2} {\bf B}^{H} \right) \nonumber \\
 &  =  & \EE (\rho_1) -  \frac{1}{t} \mathrm{Tr} \left( 
\EE({\bf S}) {\bf B} \tilde{{\bf D}} ({\bf I} + \alpha \tilde{{\bf D}})^{-2} {\bf B}^{H} \right)\ .
\end{eqnarray}
Therefore, the term
\[
 \frac{1}{t} \mathrm{Tr} 
\tilde{{\bf D}} ({\bf I} + \alpha \tilde{{\bf D}})^{-1} - \EE ( \rho_1) + \alpha  \frac{1}{t} \mathrm{Tr} \left(
\EE({\bf S}) {\bf B} \tilde{{\bf D}}^{2} ({\bf I} + \alpha \tilde{{\bf D}})^{-2} {\bf B}^{H} \right)
\]
is equal to
\begin{multline*}
 \frac{1}{t} \mathrm{Tr} \tilde{{\bf D}} ({\bf I} + \alpha \tilde{{\bf D}})^{-1}  -  \frac{1}{t} \mathrm{Tr} \left( 
\EE({\bf S}) {\bf B} \tilde{{\bf D}} ({\bf I} + \alpha \tilde{{\bf D}})^{-2} {\bf B}^{H} \right) \\
=  \frac{1}{t} \mathrm{Tr} \left[ \tilde{{\bf D}} ({\bf I} + \alpha \tilde{{\bf D}})^{-1} \left( {\bf I} - 
{\bf B}^{H} \EE ({\bf S}) {\bf B}  ({\bf I} + \alpha \tilde{{\bf D}})^{-1} \right) \right]
\end{multline*}
which, in turn, coincides with $\sigma^{2} \, \tilde{\tau}$, where $\tilde{\tau}$ is defined by
 \begin{equation}
 \label{eq:def-tildetau}
 \tilde{\tau}(\sigma^{2}) = \frac{1}{t} \mbox{Tr} \left[ \tilde{{\bf D}} \left( \sigma^{2}({\bf I} + \alpha \tilde{{\bf D}}) \right)^{-1}
 \left( {\bf I} - {\bf B}^{H} \, \mathbb{E}({\bf S}(\sigma^{2})) \, {\bf B} ({\bf I} + \alpha \tilde{{\bf D}})^{-1} \right) \right] 
 \end{equation}
Eq. (\ref{eq:intermediaire-4-exact}) 
is thus equivalent to 
\begin{equation}
\label{eq:intermediaire-5-exact}
\left( \EE({\bf S}) \left[ \sigma^{2}({\bf I} + \tilde{\tau} {\bf D}) + {\bf B} ({\bf I} + \alpha \tilde{{\bf D}})^{-1} {\bf B}^{H} \right] \right) = {\bf I} + {\bs \Delta}
\end{equation}
or equivalently to 
\[
\left( \EE \left[ {\bf S} \left( \sigma^{2}({\bf I} + \tilde{\alpha} {\bf D}) + {\bf B} ({\bf I} + \alpha \tilde{{\bf D}})^{-1} {\bf B}^{H} \right) \right] \right) =  {\bf I} + \sigma^{2} (\tilde{\alpha} - \tilde{\tau}) \mathbb{E}({\bf S}) {\bf D} + \bs{\Delta}
\]
or to 
\begin{equation}
\label{eq:second-approximation-E(S)}
\EE({\bf S}) = {\bf R} +  \sigma^{2} (\tilde{\alpha} - \tilde{\tau}) \mathbb{E}({\bf S}) {\bf D} {\bf R} + \bs{\Delta} {\bf R}\ .
\end{equation}
We now verify that if ${\bf M}$ 
is a deterministic, uniformly bounded matrix for the spectral norm as $t \to \infty$, then
$
{t}^{-1} \mathrm{Tr} {\bs \Delta} {\bf R} {\bf M} = O\left({t^{-2}}\right)\ .
$
For this, we write $\frac{1}{t} \, \mathrm{Tr} {\bs \Delta} {\bf R} {\bf M}$ as
$\frac{1}{t} \, \mathrm{Tr} {\bs \Delta} {\bf R} {\bf M} = T_1 + T_2 - T_3$ 
where 
\[
\begin{array}{ccc}
T_1  & = &  \EE \left[ \overcirc{\eta}  \; \frac{1}{t} \mathrm{Tr} \left( {\bf S} {\bs \Sigma} \tilde{{\bf D}} ({\bf I} + \alpha \tilde{{\bf D}})^{-1} {\bs \Sigma}^{H} {\bf R} {\bf M}\right) \right]  \ ,\\
T_2 & = &  \EE \left((\overcirc{\rho_1} + \overcirc{\rho_2}) \,  \frac{1}{t} \, \mathrm{Tr} (\overcirc{{\bf S}} {\bf D} {\bf R} {\bf M}) \right)\ , \\
T_3 & = &  \EE \left(  \overcirc{\eta} \, \overcirc{\rho_3} \right) \,  \frac{1}{t} \, \mathrm{Tr} \left( \EE({\bf S}) {\bf D} {\bf R} {\bf M} \right)\ .
\end{array}
\]
We denote by $\rho_4$ the term 
\[
\rho_4 =  \frac{1}{t} \mathrm{Tr} \left( {\bf S} {\bs \Sigma} \tilde{{\bf D}} ({\bf I} + \alpha \tilde{{\bf D}})^{-1} {\bs \Sigma}^{H} {\bf R} {\bf M}\right)
\]
and notice that $T_1 = \EE ( \overcirc{\eta} \; \overcirc{\rho_4})$.
Eq. (\ref{eq:vitesse-convergence-normalise}) implies that $\EE(
\overcirc{\eta}^{2})$ and $\EE \left[ \frac{1}{t} \, \mathrm{Tr}
  \left( \overcirc{{\bf S}} {\bf D} {\bf R} {\bf M}) \right)
\right]^{2}$ are $O({t^{-2}})$ terms. Moreover, matrix ${\bf R}$ 
is uniformly bounded for the spectral norm as $t \to \infty$ (see (\ref{eq:bornesR}). Lemma
\ref{le:eval-variances-utiles} immediately shows that for each $i=1,
2, 3$, $\EE ( \overcirc{\rho_i}^{2})$ is a $O({t^{-2}})$
term.  The Cauchy-Schwarz inequality eventually provides $\frac{1}{t}
\mathrm{Tr} {\bs \Delta} {\bf R} {\bf M} = O({t^{-2}})$. 

In order to establish (\ref{eq:E(S)-R-zero}), it remains to show that $\tilde{\alpha} - \tilde{\tau} \to 0$. 
For this, we remark that exchanging the roles of 
matrices $\bs{\Sigma}$ and $\bs{\Sigma}^{H}$ leads to the following relation
\begin{equation}
\label{eq:second-approximation-E(tildeS)}
\EE(\tilde{{\bf S}}) = \tilde{{\bf R}} +  \sigma^{2} (\alpha - \tau) \mathbb{E}(\tilde{{\bf S}} \tilde{{\bf D}} \tilde{{\bf R}})
+ \tilde{\bs{\Delta}} \tilde{{\bf R}}
\end{equation}
where $\tau(\sigma^{2})$ is defined by
\begin{equation}
\label{eq:def-tau}
\tau(\sigma^{2}) = \frac{1}{t} \mbox{Tr} \left[ {\bf D} \left( \sigma^{2}({\bf I} + \tilde{\alpha} {\bf D}) \right)^{-1}
\left( {\bf I} - {\bf B} \, \mathbb{E}(\tilde{{\bf S}}(\sigma^{2})) \, {\bf B}^{H} ({\bf I} + \tilde{\alpha} {\bf D})^{-1} \right) \right] 
\end{equation}
and where $\tilde{\bs{\Delta}}$, the analogue of $\bs{\Delta}$, satisfies
\begin{equation}
\label{eq:prop-tildeDelta}
\frac{1}{t} \mathrm{Tr}(\tilde{\bs{\Delta}} \tilde{{\bf M}}) = O\left(\frac{1}{t^{2}}\right)
\end{equation}
for every matrix $\tilde{{\bf M}}$ uniformly bounded for the spectral norm.

Equations (\ref{eq:second-approximation-E(S)}) and (\ref{eq:second-approximation-E(tildeS)}) allow to evaluate $\tilde{\alpha}$ and $\tilde{\tau}$.  
More precisely, writing $\tilde{\alpha} = \frac{1}{t} \mathrm{Tr}(\tilde{{\bf D}} \mathbb{E}(\tilde{{\bf S}}))$ and using the expression 
(\ref{eq:second-approximation-E(tildeS)}) of $\mathbb{E}(\tilde{{\bf S}})$, we obtain that 
\begin{equation}
\label{eq:exprebis-tildealpha}
\tilde{\alpha} = \frac{1}{t} \mathrm{Tr}(\tilde{{\bf D}} \tilde{{\bf R}}) + \sigma^{2}(\alpha - \tau)  \frac{1}{t} \mathrm{Tr}(\tilde{{\bf D}} 
\mathbb{E}(\tilde{{\bf S}}) \tilde{{\bf D}} \tilde{{\bf R}}) + \frac{1}{t} \mathrm{Tr} (\tilde{{\bf D}} \tilde{\bs{\Delta}} \tilde{{\bf R}})\ .
\end{equation}
Similarly, replacing $\mathbb{E}({\bf S})$ by  (\ref{eq:second-approximation-E(S)}) into the expression (\ref{eq:def-tildetau}) of $\tilde{\tau}$, we get that 
\begin{equation}
\label{eq:exprebis-tildetau}
\begin{array}{ccc}
\tilde{\tau} & = & \frac{1}{t} \mathrm{Tr} \left[ \tilde{{\bf D}} (\sigma^{2}({\bf I} + \alpha \tilde{{\bf D}})^{-1}({\bf I} 
- {\bf B}^{H} {\bf R} {\bf B} ({\bf I} + \alpha \tilde{{\bf D}})^{-1} \right] \\
     &   & -(\tilde{\alpha} - \tilde{\tau})  \frac{1}{t} \mathrm{Tr} \left[ \tilde{{\bf D}} ({\bf I} + \alpha \tilde{{\bf D}})^{-1}
 {\bf B}^{H} \mathbb{E}({\bf S}) {\bf D} {\bf R} {\bf B} ({\bf I} + \alpha \tilde{{\bf D}})^{-1} \right] \\
    &    & -  \frac{1}{t} \mathrm{Tr} \left[ \tilde{{\bf D}} (\sigma^{2}({\bf I} + \alpha \tilde{{\bf D}})^{-1} {\bf B}^{H} 
\bs{\Delta} {\bf R} {\bf B} ({\bf I} + \alpha \tilde{{\bf D}})^{-1} \right] \ .
\end{array}
\end{equation}
Using standard algebra, it is easy to check that the first term of the righthandside of (\ref{eq:exprebis-tildetau}) 
coincides with $\frac{1}{t} \mathrm{Tr} (\tilde{{\bf D}} \tilde{{\bf R}}$). Substracting (\ref{eq:exprebis-tildetau}) from
(\ref{eq:exprebis-tildealpha}), we get that 
\begin{equation}
\label{eq:equation-systeme-2}
(\alpha - \tau) \tilde{u}_0 + (\tilde{\alpha} - \tilde{\tau}) \tilde{v}_0 = \tilde{\epsilon}
\end{equation}
where 
\begin{equation}
\label{eq:uv0tilde}
\begin{array}{ccc}
\tilde{u}_0 & = &   \sigma^{2} \frac{1}{t} \mathrm{Tr}(\tilde{{\bf D}} 
\mathbb{E}(\tilde{{\bf S}}) \tilde{{\bf D}} \tilde{{\bf R}}) \\
\tilde{v}_0 & = & 1 -  \frac{1}{t} \mathrm{Tr} \left[ \tilde{{\bf D}} ({\bf I} + \alpha \tilde{{\bf D}})^{-1}
 {\bf B}^{H} \mathbb{E}({\bf S}) {\bf D} {\bf R} {\bf B} ({\bf I} + \alpha \tilde{{\bf D}})^{-1} \right] \\
\tilde{\epsilon} & = & \frac{1}{t} \mathrm{Tr} (\tilde{{\bf D}} \tilde{\bs{\Delta}} \tilde{{\bf R}}) +  \frac{1}{t} \mathrm{Tr} \left[ \tilde{{\bf D}} (\sigma^{2}({\bf I} + \alpha \tilde{{\bf D}})^{-1} {\bf B}^{H} 
\bs{\Delta} {\bf R} {\bf B} ({\bf I} + \alpha \tilde{{\bf D}})^{-1} \right]\ .
\end{array}
\end{equation}
Using the properties of $\bs{\Delta}$ and $\tilde{\bs{\Delta}}$, we get that $\tilde{\epsilon} = 0({t^{-2}})$. 

Similar calculations allow to evaluate $\alpha$ and $\tau$, and to obtain 
\begin{equation}
\label{eq:equation-systeme-1}
(\alpha - \tau) u_0 + (\tilde{\alpha} - \tilde{\tau}) v_0 = \epsilon
\end{equation}
where 
\begin{equation}
\label{eq:uv0}
\begin{array}{ccc}
u_0 & = &  1 -  \frac{1}{t} \mathrm{Tr} \left[ {\bf D} ({\bf I} + \tilde{\alpha} {\bf D})^{-1}
 {\bf B} \mathbb{E}(\tilde{{\bf S}}) \tilde{{\bf D}} \tilde{{\bf R}} {\bf B}^{H} ({\bf I} + \tilde{\alpha} {\bf D})^{-1} \right] \\
v_0 & = &  \sigma^{2} \frac{1}{t} \mathrm{Tr}({\bf D} 
\mathbb{E}({\bf S}) {\bf D} {\bf R}) 
\end{array}
\end{equation}
and where $\epsilon = O({t^{-2}})$. (\ref{eq:equation-systeme-1}, \ref{eq:equation-systeme-2}) can be written as
\begin{equation}
\label{eq:system_alpha-tau}
\left( \begin{array}{cc} u_0 & v_0 \\ \tilde{u}_0 & \tilde{v}_0 \end{array} \right) \left( \begin{array}{c} \alpha - \tau \\ \tilde{\alpha} - \tilde{\tau} \end{array} \right) = 
\left( \begin{array}{c} \epsilon \\ \tilde{\epsilon} \end{array} \right) \ .
\end{equation}
If the determinant $u_0 \tilde{v}_0 - \tilde{u}_0 v_0$ of the $2 × 2$ matrix governing the system is nonzero, 
$\alpha - \tau$ and $\tilde{\alpha} - \tilde{\tau}$ are given by: 
\begin{equation}
\label{eq:expre_alpha-tau}
\alpha - \tau = \  \frac{\tilde{v_0} \epsilon - v_0 \tilde{\epsilon} }{u_0 \tilde{v}_0 - \tilde{u}_0 v_0},\qquad
\tilde{\alpha} - \tilde{\tau} \ = \ \frac{u_0 \tilde{\epsilon} - \tilde{u}_0 \epsilon}{u_0 \tilde{v}_0 - \tilde{u}_0 v_0}\ ,
\end{equation}
As matrices ${\bf R}$ and $\mathbb{E}({\bf S})$ are less than $\frac{1}{\sigma^{2}} {\bf I}_r$ 
and matrices  $\tilde{{\bf R}}$ and $\mathbb{E}(\tilde{{\bf S}})$ are less than $\frac{1}{\sigma^{2}} {\bf I}_t$, it is easy to check that 
$u_0, v_0, \tilde{u}_0, \tilde{v}_0$ are uniformly bounded. As
$\epsilon$ and $\tilde{\epsilon}$ are $O({t^{-2}})$ terms,
$(\alpha - \tau)$ and $(\tilde{\alpha} - \tilde{\tau})$ will
converge to $0$ as long as the inverse
$(u_0 \tilde{v}_0 - \tilde{u}_0 v_0)^{-1}$ of the determinant is
uniformly bounded. For the moment, we show this property for $\sigma^{2}$ large enough. 
For this, we study the behaviour of coefficients $u_0, \tilde{u}_0, v_0, \tilde{v}_0$ for large
enough values of $\sigma^{2}$. It is easy to check that:
\begin{equation}
\label{eq:bornes_u0_v0-tildeu0-tildev0}
\begin{array}{ccc}
u_0 & \geq &1 - \frac{1}{\sigma^{4}} \; \frac{r}{t} \; d_{\max} \tilde{d}_{\max} b_{\max}^{2}\ , \\
\tilde{v}_0 & \geq &  1 - \frac{1}{\sigma^{4}}  \; d_{\max} \tilde{d}_{\max} b_{\max}^{2} \ ,\\
\tilde{u}_0 &  \leq & \frac{\tilde{d}_{\max}^{2}}{\sigma^{2}} \\
v_0 &  \leq & \frac{r}{t} \; \frac{d_{\max}^{2}}{\sigma^{2}}
\end{array}  
\end{equation}
As $\frac{t}{r} \to c$, it is clear that there exists
$\sigma^{2}_0$ and an integer $t_0$ for which $u_0 \geq 1/2,
\tilde{v}_0 \geq 1/2, \tilde{u}_0 \leq 1/4, v_0 \leq 1/4$ for $t \geq
t_0$ and $\sigma^{2} \geq \sigma^{2}_0$. Therefore, $u_0 \tilde{v}_0 -
\tilde{u}_0 v_0 > \frac{3}{16}$ for $t \geq t_0$ and $\sigma^{2} \geq
\sigma^{2}_0$. Eq. (\ref{eq:expre_alpha-tau}) thus implies that if
$\sigma^{2} \geq \sigma^{2}_0$, then $\alpha - \tau$ and
$\tilde{\alpha} - \tilde{\tau}$ are of the same order of magnitude as
$\epsilon = O({t^{-2}})$, and therefore converge to 0 when $t
\to \infty$.  It remains to prove that this convergence still
holds for $0 < \sigma^{2} < \sigma^{2}_0$. For this, we shall rely on
Montel's theorem (see e.g. \cite{cartan-1978}), a tool frequently used in the context of large
random matrices. It is based on the observation that, considered as
functions of parameter $\sigma^{2}$, $\alpha(\sigma^{2}) -
\tau(\sigma^{2})$ and $\tilde{\alpha}(\sigma^{2}) -
\tilde{\tau}(\sigma^{2})$ can be extended to holomorphic functions on
$\mathbb{C} - \mathbb{R}^{-}$ by replacing $\sigma^{2}$ by a complex
number $z$. Moreover, it can be shown that these holomorphic functions
are uniformly bounded on each compact subset $K$ of $\mathbb{C} -
\mathbb{R}^{-}$, in the sense that $\sup_{t} \sup_{z \in K} \;
|\alpha(z) - \tau(z) | < \infty$ and $\sup_{t} \sup_{z \in K} \;
|\tilde{\alpha}(z) - \tilde{\tau}(z) | < \infty$.  Using Montel's
theorem, it can thus be shown that if $\alpha(\sigma^{2}) -
\tau(\sigma^{2})$ and $\tilde{\alpha}(\sigma^{2}) -
\tilde{\tau}(\sigma^{2})$ converge towards zero for each $\sigma^{2} >
\sigma^{2}_0$, then for each $z \in \mathbb{C} - \mathbb{R}^{-}$,
$\alpha(z) - \tau(z)$ and $\tilde{\alpha}(z) - \tilde{\tau}(z)$
converge as well towards 0.  This in particular implies that
$\alpha(\sigma^{2}) - \tau(\sigma^{2})$ and
$\tilde{\alpha}(\sigma^{2}) - \tilde{\tau}(\sigma^{2})$ converge
towards 0 for each $\sigma^{2} > 0$.  For more details, the reader may e.g.
refer to \cite{HLN07}. This completes the proof of (\ref{eq:E(S)-R-zero}). 

We note that Montel's theorem does not guarantee that $\alpha - \tau$
and $\tilde{\alpha} - \tilde{\tau}$ are still $O(t^{-2})$
terms for $\sigma^{2} < \sigma^{2}_0$. This is one of the purpose of
the proof of Step 2 below. \\

In order to finish the proof of Proposition \ref{prop:E(S)-R}, it remains to check 
that (\ref{eq:R-T-zero}) holds. We first observe that 
$ {\bf R} - {\bf T} = {\bf R} \left( {\bf
    T}^{-1} - {\bf R}^{-1} \right) {\bf T} $. Using the expressions
of ${\bf R}^{-1}$ and ${\bf T}^{-1}$, multiplying by ${\bf M}$, and
taking the trace yields:
\begin{eqnarray}
\label{eq:expreR-T}
\frac{1}{t} \mbox{Tr} \left[ {\bf M} \left( {\bf R} - {\bf T} \right) \right] & = & (\tilde{\beta} - \tilde{\alpha}) \; \sigma^{2} \frac{1}{t} \mbox{Tr} ({\bf M} {\bf R} {\bf D} {\bf T}) + \nonumber \\
  &  & (\alpha - \beta) \; \frac{1}{t} \mbox{Tr} \left[ {\bf M} {\bf R} {\bf B} ({\bf I} + \beta \tilde{{\bf D}})^{-1}  \tilde{{\bf D}} ({\bf I} + \beta \tilde{{\bf D}})^{-1} {\bf B}^{H} {\bf T} \right]\ .
\end{eqnarray}
As the terms $  \frac{\sigma^{2}}{t} \mbox{Tr} ({\bf M} {\bf R} {\bf
  D} {\bf T})$ and $\frac{1}{t} \mbox{Tr} \left[ {\bf M} {\bf R} {\bf
    B} ({\bf I} + \beta \tilde{{\bf D}})^{-1} \tilde{{\bf D}} ({\bf I}
  + \beta \tilde{{\bf D}})^{-1} {\bf B}^{H} {\bf T} \right]$ are
uniformly bounded, it is sufficient to establish that $ (\alpha -
\beta)$ and $(\tilde{\alpha} - \tilde{\beta})$ converge towards $0$. 
For this, we note that (\ref{eq:E(S)-R-zero}) implies that 
\begin{equation}
\label{eq:premier-rapprochement}
\alpha \ =\ \frac{1}{t} \mbox{Tr} \left({\bf D} {\bf R} \right) + \epsilon^{'},\qquad 
\tilde{\alpha} \ = \ \frac{1}{t} \mbox{Tr} \left(\tilde{{\bf D}} \tilde{{\bf R}} \right) + \tilde{\epsilon}^{'}\ ,
\end{equation}
where $\epsilon^{'}$ and $\tilde{\epsilon}^{'}$ converge towards 0. We express  
$(\alpha - \beta) = \frac{1}{t} \mbox{Tr}{\bf D} ({\bf R} - {\bf T}) + \epsilon$. Using  $ {\bf R} - {\bf T} = {\bf R} \left( {\bf
    T}^{-1} - {\bf R}^{-1} \right) {\bf T} $, multiplying by ${\bf D}$ from both sides, and taking the trace yields
\begin{equation}
\label{eq:alpha-beta}
(\alpha - \beta) \; \left( 1 -  \frac{1}{t} \mbox{Tr} \left[ {\bf D} {\bf R} {\bf B} ({\bf I} + \beta \tilde{{\bf D}})^{-1}  \tilde{{\bf D}} ({\bf I} + \beta \tilde{{\bf D}})^{-1} {\bf B}^{H} {\bf T} \right] \right) + (\tilde{\alpha} - \tilde{\beta}) \;  \sigma^{2} \frac{1}{t} \mbox{Tr} ({\bf D} {\bf R} {\bf D} {\bf T}) = \epsilon^{'}.
\end{equation}
Similarly, we obtain that
\begin{equation}
\label{eq:tildealpha-tildebeta}
(\alpha - \beta) \;  \sigma^{2} \frac{1}{t} \mbox{Tr} (\tilde{{\bf D}} \tilde{{\bf R}} \tilde{{\bf D}} \tilde{{\bf T}}) + 
 (\tilde{\alpha} - \tilde{\beta}) \; \left( 1 -  \frac{1}{t} \mbox{Tr} \left[ \tilde{{\bf D}} \tilde{{\bf R}} {\bf B}^{H} ({\bf I} + \tilde{\beta} {\bf D})^{-1}  {\bf D} ({\bf I} + \tilde{\beta} {\bf D})^{-1} {\bf D} \tilde{{\bf T}} \right] \right) =  \tilde{\epsilon}^{'}. 
\end{equation}
Equations (\ref{eq:alpha-beta}) and (\ref{eq:tildealpha-tildebeta}) can be interpreted as a linear systems w.r.t. $(\alpha - \beta)$ and 
$(\tilde{\alpha} - \tilde{\beta})$. Using the same approach as in the proof of (\ref{eq:E(S)-R-zero}), we prove 
that  $(\alpha - \beta)$ and $(\tilde{\alpha} - \tilde{\beta})$ converge towards 0. This establishes (\ref{eq:R-T-zero})
and completes the proof of Proposition (\ref{prop:E(S)-R}). 
\end{proof}

\subsubsection{Second step: $\frac{1}{r}\mathrm{Tr} {\bf M}
  (\mathbb{E}({\bf S}) - {\bf R})$ and $\frac{1}{r} \mathrm{Tr} {\bf
    M} ({\bf R} - {\bf T})$ are $O(t^{-2})$ terms} This section is devoted to the proof of the following proposition.
\begin{proposition}
\label{prop:rate-E(S)-R}
For each deterministic $r × r$ matrix ${\bf M}$, uniformly
bounded (for the spectral norm) as $t \to \infty$, we have:
\begin{equation}
\label{eq:E(S)-R}
 \frac{1}{t} \mathrm{Tr}\left[ {\bf M} \left( \mathbb{E}({\bf S}) - {\bf R} \right) \right] = O({t^{-2}})
\end{equation}
\begin{equation}
\label{eq:R-T}
\frac{1}{t} \mathrm{Tr}\left[ {\bf M} \left( {\bf R}) - {\bf T} \right) \right] = O({t^{-2}})
\end{equation}
\end{proposition}
\begin{proof}
We first establish (\ref{eq:E(S)-R}). For this, we prove that the inverse of the  determinant 
$u_0 \tilde{v}_0 - \tilde{u}_0 v_0$ of linear system (\ref{eq:system_alpha-tau})
is uniformly bounded for each $\sigma^{2} > 0$.  In order to state the corresponding result, we
define $(u,v,\tilde{u}, \tilde{v})$ by
\begin{equation}
\label{eq:defuv}
\begin{array}{ccc}
u & = &   1 -  \frac{1}{t} \mathrm{Tr}( \tilde{{\bf D}} \tilde{{\bf T}} {\bf B}^{H} 
({\bf I} + \tilde{\beta} {\bf D})^{-1}  {\bf D} ({\bf I} + \tilde{\beta} {\bf D})^{-1} {\bf B} \tilde{{\bf T}})\\
\tilde{v} & = &    1 -  \frac{1}{t} \mathrm{Tr}( {\bf D} {\bf T} {\bf B} 
({\bf I} + \beta \tilde{{\bf D}})^{-1} \tilde{{\bf D}} ({\bf I} + \beta \tilde{{\bf D}})^{-1} {\bf B}^{H} {\bf T}) \\
v & = & \sigma^{2} \frac{1}{t} \mathrm{Tr}( {\bf D} {\bf T} {\bf D} {\bf T}) \\
\tilde{u} & = &  \sigma^{2} \frac{1}{t} \mathrm{Tr}( \tilde{{\bf D}} \tilde{{\bf T}} \tilde{{\bf D}} \tilde{{\bf T}})
\end{array}\ .
\end{equation}
The expressions of $(u,v,\tilde{u}, \tilde{v})$ nearly coincide with the expressions of coefficients 
$(u_0,v_0,\tilde{u}_0, \tilde{v}_0)$, the only difference being that, in the definition of $(u,v,\tilde{u}, \tilde{v})$, matrices ($\mathbb{E}({\bf S}),{\bf R}$) 
are both replaced by matrix ${\bf T}$, matrices  ($\mathbb{E}(\tilde{{\bf S}}),\tilde{{\bf R}}$) 
are both replaced by matrix $\tilde{{\bf T}}$ and scalars $(\alpha, \tilde{\alpha})$ are
replaced by scalars  $(\beta, \tilde{\beta})$. (\ref{eq:E(S)-R-zero}) and 
(\ref{eq:R-T-zero}) immediately imply that 
$(u_0,v_0,\tilde{u}_0, \tilde{v}_0)$ can be written as
\begin{equation}
\label{eq:comparaison-uv}
u_0  =  u + \epsilon_{u},\quad \tilde{v}_0  =  \tilde{v} + \tilde{\epsilon}_v,\quad v_0   =  v + \epsilon_v,\quad 
\tilde{u}_0  =  \tilde{u} + \tilde{\epsilon}_u \ ,
\end{equation}
where $\epsilon_u, \tilde{\epsilon}_v, \tilde{\epsilon}_u, \epsilon_v$ converge to $0$ when $t \to \infty$. 
The behaviour of $u \tilde{v} - \tilde{u} v$ is provided in the following Lemma, whose proof is given in paragraph 
\ref{subsub:proof-le-determinant}. 
\begin{lemma}
\label{le:determinant}
Coefficients $(u,v, \tilde{u}, \tilde{v})$ satisfy: 
(i) $u = \tilde{v}$,
(ii) $0 < u < 1$ and $\inf_{t} u > 0$,
(iii) $ 0 < u \tilde{v} - \tilde{u} v < 1$ and $\sup_{t} \frac{1}{u \tilde{v} - \tilde{u} v} < \infty$.
\end{lemma}
(\ref{eq:comparaison-uv}) and Lemma \ref{le:determinant} immediately imply that it exists $t_0$ such that 
$ 0 < u_0 \tilde{v}_0 - \tilde{u}_0 v_0 \leq 1$ for each $t \geq t_0$ and 
\begin{equation}
\label{eq:borne_inverse-det}
\sup_{t \geq t_0} \frac{1}{u_0 \tilde{v}_0 - \tilde{u}_0 v_0} < \infty\ .
\end{equation}
This eventually shows $\alpha - \tau$ and $\tilde{\alpha} -
\tilde{\tau}$ are of the same order of magnitude than $\epsilon$ and
$\tilde{\epsilon}$, i.e. are $O({t^{-2}})$ terms. 

In order to prove (\ref{eq:R-T}), we first remark that, by (\ref{eq:E(S)-R}), $\epsilon^{'}$ and $\tilde{\epsilon}^{'}$ defined by
(\ref{eq:premier-rapprochement}) are $O({t^{-2}})$ terms. It is thus sufficient to establish that the inverse of the determinant
of the linear system associated to equations (\ref{eq:alpha-beta}) and (\ref{eq:tildealpha-tildebeta}) 
is uniformly bounded. Eq. (\ref{eq:R-T-zero}) implies that the behaviour of this determinant is equivalent to the study 
of $u \tilde{v} - \tilde{u} v$. Eq. (\ref{eq:R-T}) thus follows from Lemma \ref{le:determinant}. This completes the proof 
of Proposition \ref{prop:rate-E(S)-R}. 

\end{proof}

\subsubsection{Proof of Lemma \ref{le:determinant}.} 
\label{subsub:proof-le-determinant}
In order to establish item (i), we notice that a direct application of the matrix inversion
Lemma yields:
\begin{equation}
\label{eq:lemme-inversion}
\tilde{{\bf T}} {\bf B}^{H} ({\bf I} + \tilde{\beta} {\bf D})^{-1} =  ({\bf I} + \beta \tilde{{\bf D}})^{-1} {\bf B}^{H} {\bf T}\ .
\end{equation}
The equality $u = \tilde{v}$ immediately follows from (\ref{eq:lemme-inversion}). 

The proofs of (ii) and (iii) are based on the observation that
function $\sigma^{2} \to \sigma^{2} \beta(\sigma^{2})$ is increasing
while function $\sigma^{2} \to \tilde{\beta}(\sigma^{2})$ is
decreasing. This claim is a consequence of Eq.
(\ref{eq:represention-beta}) that we recall below:
$$
\beta(\sigma^{2}) \ = \ \int_{\mathrm{R}^{+}} \frac{d \mu_b(\lambda) }{\lambda + \sigma^{2}} ,\quad 
\tilde{\beta}(\sigma^{2}) \ = \ \int_{\mathrm{R}^{+}} \frac{d \tilde{\mu}_b(\lambda) }{\lambda + \sigma^{2}} \ ,
$$
where 
$\mu_b(\mathbb{R}^{+}) = \frac{1}{t} \mathrm{Tr}({\bf D})$ and 
$\tilde{\mu}_b(\mathbb{R}^{+}) = \frac{1}{t} \mathrm{Tr}(\tilde{{\bf
    D}})$. Note that $\tilde{\beta}$ is decreasing because $\sigma^{2}
\mapsto \frac{1}{\lambda + \sigma^{2}}$ is decreasing and $\sigma^{2}
\beta(\sigma^{2})$ is increasing because $\sigma^{2} \mapsto
\frac{\sigma^{2}}{\lambda + \sigma^{2}}$ is increasing. Denote by $'$
the differentiation operator w.r.t. $\sigma^{2}$. Then, $(\sigma^{2}
\beta)^{'} > 0$ and $\tilde{\beta}^{'} < 0$ for each $\sigma^{2}$.  We
now differentiate relations (\ref{eq:exprebeta}) w.r.t. $\sigma^{2}$.
After some algebra, we obtain:
\begin{equation}
\label{eq:derivation}
\begin{array}{ccc}
\tilde{v} \; (\sigma^{2} \beta)^{'} + \sigma^{2} v \; \tilde{\beta}^{'} & = & \frac{1}{t} \mathrm{Tr}( {\bf D} {\bf T} {\bf B} 
({\bf I} + \beta \tilde{{\bf D}})^{-1} ({\bf I} + \beta \tilde{{\bf D}})^{-1} {\bf B}^{H} {\bf T}) \\
\frac{\tilde{u}}{\sigma^{2}} \; (\sigma^{2} \beta)^{'} + u \tilde{\beta}^{'} & = &  -\frac{1}{t} \mathrm{Tr} \tilde{{\bf T}} \tilde{{\bf D}} \tilde{{\bf T}}
\end{array}\ .
\end{equation}
As $\tilde{\beta}^{'} < 0$, the first equation of
(\ref{eq:derivation}) implies that $\tilde{v} \; (\sigma^{2} \beta)^{'} > 0$.
As $ (\sigma^{2} \beta)^{'} > 0$, this yields $\tilde{v} > 0$. As $\tilde{v} < 1$
clearly holds, the first part of (ii) is proved.

We now prove that $\inf_{t} \tilde{v} > 0$. The first equation of (\ref{eq:derivation})
yields: 
\begin{equation}
\label{eq:minoration1-u}
\tilde{v} > - \sigma^{2} v  \tilde{\beta}^{'} \frac{1}{ (\sigma^{2} \beta)^{'}}\ .
\end{equation}
In the following, we show that $\inf_t \frac{1}{ (\sigma^{2} \beta)^{'}} > 0$, $\inf_t |\tilde{\beta}^{'}| > 0$
and that $\inf_t \, v > 0$. 

By representation (\ref{eq:represention-beta}), 
$$
-  \tilde{\beta}^{'} \ = \  \int_{\mathrm{R}^{+}} \frac{d \tilde{\mu}_b(\lambda)}{(\lambda + \sigma^{2})^{2}}  
\quad\textrm{and}\quad
(\sigma^{2} \beta(\sigma^{2}))^{'} \ 
= \ \int_{\mathrm{R}^{+}} \frac{\lambda d \mu_b(\lambda) }{(\lambda + \sigma^{2})^{2}}\ .
$$
As $\frac{\lambda}{(\lambda + \sigma^{2})^{2}} \leq
\frac{1}{\sigma^{2}}$ for $\lambda \geq 0$, $(\sigma^{2} \beta)^{'}
\leq \frac{1}{\sigma^{2}} \mu_b(\mathbb{R}^{+}) = \frac{1}{t}
\mathrm{Tr} {\bf D}$. Therefore, the term $\frac{1}{ (\sigma^{2}
  \beta)^{'}}$ is lowerbounded by $\sigma^{2} (\frac{1}{t} \mathrm{Tr}
{\bf D})^{-1}$.  As $\frac{1}{t} \mathrm{Tr} {\bf D} \leq \frac{r}{t}
d_{\max}$, we have $\inf_{t} \frac{1}{ (\sigma^{2} \beta)^{'}} > 0$.
 
We now establish that $\inf_{t} | \tilde{\beta}^{'}| > 0$. We first use Jensen's inequality: As measure 
$ (\frac{1}{t} \mathrm{Tr} \tilde{{\bf D}})^{-1} \;d \tilde{\mu}_b(\lambda)$ is a probability distribution:
\[
\left[ \int_{\mathrm{R}^{+}} \frac{1}{\lambda + \sigma^{2}}  \left(\frac{1}{t} \mathrm{Tr} \tilde{{\bf D}}\right)^{-1} \;d \tilde{\mu}_b(\lambda) \right]^{2} 
\leq \int_{\mathrm{R}^{+}} \frac{1}{(\lambda + \sigma^{2})^{2}}  
\left(\frac{1}{t} \mathrm{Tr} \tilde{{\bf D}}\right)^{-1} \;d \tilde{\mu}_b(\lambda)\ .
\]
In other words, $| \tilde{\beta}^{'} | =  \int_{\mathrm{R}^{+}} \frac{1}{(\lambda + \sigma^{2})^{2}} d \tilde{\mu}_b(\lambda)$ 
satisfies
\[
| \tilde{\beta}^{'} | \geq \frac{1}{\frac{1}{t} \mathrm{Tr} \tilde{{\bf D}}}  
\left[ \int_{\mathrm{R}^{+}} \frac{1}{\lambda + \sigma^{2}} d \tilde{\mu}_b(\lambda) \right]^{2} = \frac{1}{\frac{1}{t} \mathrm{Tr} \tilde{{\bf D}}} \tilde{\beta}^{2}\ .
\]
As mentioned above, $(\frac{1}{t} \mathrm{Tr} \tilde{{\bf D}})^{-1}$ is lower-bounded by $(d_{\max})^{-1}$. Therefore, it remains to 
establish that $\inf_{t} \tilde{\beta}^{2} > 0$, or equivalently that $\inf_{t} \tilde{\beta} > 0$. For this, we assume that
$\inf_{t} \tilde{\beta}_t(\sigma^{2}) = 0$ (we indicate that $\tilde{\beta}$ depends both on $\sigma^{2}$ and $t$). 
Therefore, there exists an increasing sequence of integers $(t_k)_{k \geq 0}$ for which 
$
\lim_{k \to \infty} \tilde{\beta}_{t_k}(\sigma^{2}) = 0
$
i.e. 
$
\lim_{k \to \infty} \int_{\mathrm{R}^{+}} \frac{1}{\lambda + \sigma^{2}} \; d \tilde{\mu}_b^{(t_k)}(\lambda) = 0\ ,
$
where $\tilde{\mu}_b^{(t_k)}$ is the positive measure associated with $\tilde{\beta}_{t_k}(\sigma^{2})$. As $\tilde{{\bf D}}$ is 
uniformly bounded, the sequence $(\tilde{\mu}_b^{(t_k)})_{k \geq 0}$ is tight. One can therefore extract from  $(\tilde{\mu}_b^{(t_k)})_{k \geq 0}$
a subsequence $(\tilde{\mu}_b^{(t^{'}_l)})_{l \geq 0}$ that converges weakly to a certain measure $\tilde{\mu}_b^{*}$ which of course satisfies 
\[
\int_{\mathrm{R}^{+}} \frac{1}{\lambda + \sigma^{2}} \; d \tilde{\mu}_b^{*}(\lambda) = 0 \ .
\]
This implies that $\tilde{\mu}_b^{*} = 0$, and thus $\tilde{\mu}_b^{*}(\mathrm{R}^{+}) = 0$, while the convergence of $(\tilde{\mu}_b^{(t^{'}_l)})_{l \geq 0}$ gives
\[
\tilde{\mu}_b^{*}(\mathrm{R}^{+}) = \lim_{l \to \infty}  \tilde{\mu}_b^{(t^{'}_l)}(\mathrm{R}^{+}) = \lim_{l \to \infty} \frac{1}{t^{'}_l}
\mathrm{Tr} \tilde{{\bf D}}_{t^{'}_l} > 0
\]
by assumption (\ref{ass-borneinf}). Therefore, the assumption $\inf_{t} \tilde{\beta}_t(\sigma^{2}) = 0$ leads to a contradiction. 
Thus,  $\inf_{t} \tilde{\beta}_t(\sigma^{2}) > 0$ and $\inf_t |\tilde{\beta}^{'}| > 0$ is proved.  

We finally establish that $v$ is lower-bounded, i.e. that $\inf_{t} \frac{1}{t} \mathrm{Tr} {\bf D} {\bf T} {\bf D} {\bf T} > 0$. 
For any Hermitian positive matrix ${\bf M}$, 
\[
\frac{1}{t} \mathrm{Tr}({\bf M}^{2}) \geq \left[ \frac{1}{t}  \mathrm{Tr}({\bf M}) \right]^{2}\ .
\]
We use this inequality for ${\bf M} =  {\bf T}^{1/2} {\bf D} {\bf T}^{1/2}$. This leads to 
\[
\frac{1}{t} \mathrm{Tr} {\bf D} {\bf T} {\bf D} {\bf T} = \frac{1}{t} \mathrm{Tr} {\bf M}^{2} >  \left[ \frac{1}{t}  \mathrm{Tr}({\bf M}) \right]^{2} = 
\left[ \frac{1}{t}  \mathrm{Tr}({\bf D} {\bf T}) \right]^{2} = \beta^{2}\ .
\]
Therefore, $\inf_{t} \frac{1}{t} \mathrm{Tr} {\bf D} {\bf T} {\bf D} {\bf T} \geq \inf_{t} \beta^{2}$. Using the same approach as
above, we can prove that $\inf_{t} \beta^{2} > 0$. Proof of (ii) is completed. \\

In order to establish (iii), we use the first equation of (\ref{eq:derivation}) to express $(\sigma^{2} \beta)^{'}$ in terms of $\tilde{\beta}^{'}$, and plug 
this relation into the second equation of (\ref{eq:derivation}). This gives:
\begin{equation}
\label{eq:determinant-1}
\left( u - \frac{1}{\tilde{v}} \tilde{u} v \right) \; \tilde{\beta}^{'} =  -\frac{1}{t} \mathrm{Tr} \tilde{{\bf T}} \tilde{{\bf D}} \tilde{{\bf T}} - \frac{\tilde{u}}{\sigma^{2} \tilde{v}}  \frac{1}{t} \mathrm{Tr}( {\bf D} {\bf T} {\bf B} 
({\bf I} + \beta \tilde{{\bf D}})^{-1} ({\bf I} + \beta \tilde{{\bf D}})^{-1} {\bf B}^{H} {\bf T})\ .
\end{equation}
The righthand side of (\ref{eq:determinant-1}) is negative as well as $\tilde{\beta}^{'}$. Therefore, 
$u - \frac{1}{\tilde{v}} \tilde{u} v > 0$. As $\tilde{v}$ is positive, $u \tilde{v} - \tilde{u} v $ is also positive. Moreover, 
$u$ et $\tilde{v}$ are strictly less than 1. As $\tilde{u}$ and $v$ are both strictly positive,  $u \tilde{v} - \tilde{u} v $ is strictly less than 1.  
To complete the proof of (iii), we notice that by (\ref{eq:determinant-1}), 
\[
\frac{1}{u  \tilde{v} - \tilde{u}  v} \leq \frac{| \tilde{\beta}^{'} |}{\tilde{v} \frac{1}{t} \mathrm{Tr} \tilde{{\bf T}} \tilde{{\bf D}} \tilde{{\bf T}}}.
\]
$| \tilde{\beta}^{'} | $ clearly satisfies $| \tilde{\beta}^{'} | \leq
\frac{1}{\sigma^{4}} \frac{1}{t} \mathrm{Tr} \tilde{{\bf D}}$ and is
thus upper bounded by $\frac{\tilde{d}_{\max}}{\sigma^{4}}$. (ii)
implies that $\sup_{t} \frac{1}{\tilde{v}} < + \infty$. It remains to verify
that $\inf_{t} \frac{1}{t} \mathrm{Tr} \tilde{{\bf T}} \tilde{{\bf D}}
\tilde{{\bf T}} > 0$. Denote by $x = \frac{1}{t} \mathrm{Tr} \tilde{{\bf
    T}} \tilde{{\bf D}} \tilde{{\bf T}}$.
\[
x  = \frac{1}{t} \sum_{i=1}^{t} \tilde{d}_i \sum_{j=1}^{t} |\tilde{T}_{i,j}|^{2}\ .
\] 
In order to use Jensen's inequality, we consider $\tilde{\kappa}_i =
\frac{\tilde{d}_i}{\frac{1}{t} \mathrm{Tr} \tilde{{\bf D}}}$, and
notice that $\frac{1}{t} \sum_{i=1}^{t} \tilde{\kappa}_i = 1$. $x$ can
be written as
\[
x = \frac{1}{t} \mathrm{Tr} \tilde{{\bf D}} \; \frac{1}{t} \sum_{i=1}^{t} \tilde{\kappa}_i \left[ (\sum_{j=1}^{t} |\tilde{T}_{i,j}|^{2})^{1/2} \right]^{2}\ .
\]
By Jensen's inequality 
\[
\frac{1}{t} \sum_{i=1}^{t} \tilde{\kappa}_i \left[ (\sum_{j=1}^{t} |T_{i,j}|^{2})^{1/2} \right]^{2} \geq \left[ \frac{1}{t} \sum_{i=1}^{t} 
\tilde{\kappa_i}  (\sum_{j=1}^{t} |\tilde{T}_{i,j}|^{2})^{1/2} \right]^{2} \ .
\]
Moreover, 
\[
\left[ \frac{1}{t} \sum_{i=1}^{t} \tilde{\kappa}_i  (\sum_{j=1}^{t} |\tilde{T}_{i,j}|^{2})^{1/2} \right]^{2}  \geq \left[ \frac{1}{t} \sum_{i=1}^{t} \tilde{\kappa}_i  \ \tilde{T}_{i,i} \right]^{2} = \left[ \left( \frac{1}{t} \mathrm{Tr} \tilde{{\bf D}} \right)^{-1} \, \tilde{\beta} \right]^{2}\ .
\]
Finally, 
\[
x = \frac{1}{t} \mathrm{Tr} \tilde{{\bf T}} \tilde{{\bf D}} \tilde{{\bf T}} \geq \left( \frac{1}{t} \mathrm{Tr} \tilde{{\bf D}} \right)^{-1} 
\tilde{\beta}^{2}\ .
\]
Since $\inf_{t} \tilde{\beta}^{2} > 0$, we have $\inf_{t} \frac{1}{t} \mathrm{Tr} \tilde{{\bf T}} \tilde{{\bf D}} \tilde{{\bf T}} > 0$ and the proof of (iii) is completed.

\section{Strict concavity of $\bar{I}({\bf Q})$: Remaining proofs}\label{app:constric}
\subsection{Proof of Lemma \ref{unif-strict-con}}
Remark that $\phi_m$ is strictly concave due to
  \eqref{strict-con}. Remark also that $\overline{\phi}$ is concave as a
  pointwise limit of the $\phi_m$'s. Now in order to prove the strict
  concavity of $\overline{\phi}$, assume that there exists a subinterval,
  say $(a,b) \subset [0,1]$ with $a<b$ where $\overline{\phi}$ fails to be strictly
  concave:
$$
\forall \lambda\in [0,1],\quad \overline{\phi}(\lambda a +(1-\lambda)b)= \lambda\overline{\phi}(a) +(1-\lambda) \overline{\phi}(b)\ .
$$
Otherwise stated, 
$$
\forall x\in (a,b),\quad \overline{\phi}(x) =\frac{\overline{\phi}(b) - \overline{\phi}(a)}{b-a}(x-a) +\overline{\phi}(a).
$$
Let $x\in (a,b)$ and $h>0$ be small enough so that $x-h$ and $x+h$
belong to $(a,b)$; recall the following inequality, valid for
differentiable concave functions:
$$
\frac{\phi_m(x) -\phi_m(x-h)}{h} \ge \phi_m'(x) \ge \frac{\phi_m(x+h) -\phi_m(x)}{h}\ .
$$
Letting $m\rightarrow \infty$, we obtain:
$$
\frac{\overline{\phi}(x) -\overline{\phi}(x-h)}{h}\  \ge\  \limsup_{m\rightarrow \infty} \phi_m'(x)
\ \ge\  \liminf_{m\rightarrow \infty} \phi_m'(x)\  \ge\  \frac{\overline{\phi}(x+h) -\overline{\phi}(x)}{h}\ .
$$
In particular, for all $x\in (a,b)$, $\lim_{m\rightarrow \infty} \phi_m'(x) =\frac{\overline{\phi}(b) -\overline{\phi}(a)}{b-a}\ .$
Now let $[x,x+h]\in (a,b)$. Fatou's lemma together with \eqref{strict-con} yield:
\begin{multline*}
0\ <\ \kappa h\ \le\ \int_x^{x+h} \liminf_{m\rightarrow \infty} \phi_m''(u) \,du\\  
\le\  \liminf_{m\rightarrow \infty} \int_x^{x+h} \phi_m''(u) \,du\  
=\lim_{m\rightarrow \infty} \left( \phi_m'(x+h) -\phi_m'(x)\right) \ =\ 0\ .
\end{multline*}
This yields a contradiction, therefore $\overline{\phi}$ must be strictly convex on $[0,1]$.

\subsection{Proof of (\ref{unif-strict-phi}).}
We define $\check {\bf M}$ as the $tm \times tm$ matrix given by
\[
\check {\bf M} = \check {\bf H}^{H} \left( {\bf I} + \frac{\check {\bf H} \check {\bf Q} \check {\bf H}^{H}}{\sigma^{2}} \right)^{-1} 
\frac{ \check {\bf H}}{\sigma^{2}}\ .
\]
We have:
$$
\phi_m''(\lambda) = -\frac 1m \mathbb{E}\, \mathrm{Tr} \left[ \check {\bf M}  (\check {\bf Q}_1 -\check {\bf Q}_2) \check {\bf M}  (\check {\bf Q}_1 -\check {\bf Q}_2)
\right]
$$ 
or equivalently
$$
\phi_m''(\lambda) = -\frac 1m \mathbb{E}\, \mathrm{Tr} \left[ \left( {\bf I} + \frac{\check {\bf H} \check {\bf Q} \check {\bf H}^{H}}{\sigma^{2}} \right)^{-1}
\frac{ \check {\bf H}}{\sigma^{2}}  (\check {\bf Q}_1 -\check {\bf Q}_2) \check {\bf M}  (\check {\bf Q}_1 -\check {\bf Q}_2) \check {\bf H}^{H} 
\right]
$$
Recall that $\mathrm{Tr} ({\bf AB}) \ge \lambda_{\min}({\bf
  A})\mathrm{Tr}({\bf B})$ for ${\bf A}$, ${\bf B}$ Hermitian and
nonnegative matrices. In particular:
\begin{multline*}
\mathrm{Tr} \left[ \left( {\bf I} + \frac{\check {\bf H} \check {\bf Q} \check {\bf H}^{H}}{\sigma^{2}} \right)^{-1}
\frac{ \check {\bf H}}{\sigma^{2}}  (\check {\bf Q}_1 -\check {\bf Q}_2) \check {\bf M}  (\check {\bf Q}_1 -\check {\bf Q}_2) \check {\bf H}^{H} 
\right] \\
\geq \lambda_{\min} \left({\bf I} + \frac{\check {\bf H} \check {\bf Q} \check {\bf H}^{H}}{\sigma^{2}}\right)^{-1} 
\mathrm{Tr} \left[\frac{ \check {\bf H}}{\sigma^{2}}  (\check {\bf Q}_1 -\check {\bf Q}_2) \check {\bf M}  (\check {\bf Q}_1 -\check {\bf Q}_2) \check {\bf H}^{H} \right] \ .
\end{multline*}
Similarly, we obtain that
\begin{multline*}
  \mathrm{Tr} \left[\frac{ \check {\bf H}}{\sigma^{2}}  (\check {\bf Q}_1 -\check {\bf Q}_2) \check {\bf M}  (\check {\bf Q}_1 -\check {\bf Q}_2) \check {\bf H}^{H} \right] \\
  \geq \lambda_{\min} \left({\bf I} + \frac{\check {\bf H} \check {\bf
        Q} \check {\bf H}^{H}}{\sigma^{2}}\right)^{-1} \mathrm{Tr}
  \left[ \frac{\check {\bf H}}{\sigma^{2}} (\check {\bf Q}_1 -\check
    {\bf Q}_2) \frac{ \check {\bf H}^{H} \check {\bf H}}{\sigma^{2}}
    (\check {\bf Q}_1 -\check {\bf Q}_2) \check {\bf H}^{H} \right]\ .
\end{multline*}
This eventually implies that
\begin{multline*}
\mathrm{Tr} \left[ \left( {\bf I} + \frac{\check {\bf H} \check {\bf Q} \check {\bf H}^{H}}{\sigma^{2}} \right)^{-1}
\frac{ \check {\bf H}}{\sigma^{2}}  (\check {\bf Q}_1 -\check {\bf Q}_2) \check {\bf M}  (\check {\bf Q}_1 -\check {\bf Q}_2) \check {\bf H}^{H} 
\right] \geq \\  \lambda_{\min}^{2} \left({\bf I} + \frac{\check {\bf H} \check {\bf Q} \check {\bf H}^{H}}{\sigma^{2}}\right)^{-1}
\mathrm{Tr}\left[ \frac{\check {\bf H}^H \check {\bf H}}{\sigma^2}
      (\check {\bf Q}_1 -\check {\bf Q}_2) \frac{\check {\bf H}^H
        \check {\bf H}}{\sigma^2} (\check {\bf Q}_1 -\check{\bf Q}_2)
    \right]\ .
\end{multline*}
As
$$
\lambda^2_{\min}\left({\bf I} +\frac{\check {\bf H} \check {\bf Q} \check {\bf
      H}^{H}}{\sigma^2} \right)^{-1} \ge \frac
1{\lambda_{\max}^2 \left({\bf I} +\frac{\check {\bf H} \check {\bf Q} \check {\bf
        H}^{H}}{\sigma^2} \right)} \ge \frac 1{\left(1+ \sigma^{-2}
    \| \check {\bf Q} \|\,  \| \check {\bf H}^H \check {\bf H} \|\right)^2}\ ,
$$
we have:
$$
\phi_m''(\lambda) \le -\frac 1m \mathbb{E} \left[\left( \frac 1{\left(1+ \sigma^{-2}
    \| \check {\bf Q} \|\,  \| \check {\bf H}^H \check {\bf H} \|\right)^2} 
\right)
\times  \mathrm{Tr}\left( \frac{\check {\bf H}^H \check {\bf H}}{\sigma^2}
      (\check {\bf Q}_1 -\check {\bf Q}_2) \frac{\check {\bf H}^H
        \check {\bf H}}{\sigma^2} (\check {\bf Q}_1 -\check{\bf Q}_2)
    \right) \right]\ .
$$
Let us introduce the following notations: 
$$
\alpha_m = \frac 1{\left(1+ \sigma^{-2}
    \| \check {\bf Q} \|\,  \| \check {\bf H}^H \check {\bf H} \|\right)^2}, \qquad 
\beta_m = \frac 1m \mathrm{Tr}\left[ \frac{\check {\bf H}^H \check {\bf H}}{\sigma^2}
      (\check {\bf Q}_1 -\check {\bf Q}_2) \frac{\check {\bf H}^H
        \check {\bf H}}{\sigma^2} (\check {\bf Q}_1 -\check{\bf Q}_2)
    \right]\ . 
$$
The following properties whose proofs are postponed to Appendix \ref{proof:various}
hold true:
\begin{proposition}\label{prop:various}
\begin{itemize}
\item[(i)] $\lim_{m\rightarrow \infty} \mathrm{var}(\beta_m) =0\ ,$
\item[(ii)] For all $m\ge 1$, 
$
\mathbb{E}(\beta_m) = \mathbb{E}(\beta_1) =\mathbb{E} \mathrm{Tr}
\left[ \frac{{\bf H}^H {\bf H}}{\sigma^2}
      ({\bf Q}_1 -{\bf Q}_2) \frac{{\bf H}^H
        {\bf H}}{\sigma^2} ( {\bf Q}_1 -{\bf Q}_2)
    \right] > 0\ ,
$
\item[(iii)] There exists $\delta>0$ such that for all $\lambda \in [0,1]$,
$
\liminf_{m\rightarrow \infty} \mathbb{E} (\alpha_m) \ge \delta >0\ .
$
\end{itemize}
\end{proposition}

We are now in position to establish \eqref{unif-strict-phi}. By Proposition \ref{prop:various}-(i), we have
$$
|\mathbb{E} (\alpha_m \beta_m) - \mathbb{E} (\alpha_m) \mathbb{E} (\beta_m) | \le \sqrt{\mathrm{var}(\beta_m)}
\sqrt{\mathbb{E}(\alpha_m^2)} \le \sqrt{\mathrm{var}(\beta_m)} \xrightarrow[m\rightarrow \infty]{} 0\ .
$$
By Proposition \ref{prop:various}-(ii),(iii), we have:
$$
\liminf_{m\rightarrow \infty} \mathbb{E} (\alpha_m \beta_m) =\liminf_{m\rightarrow \infty} \mathbb{E} (\alpha_m)
\mathbb{E}(\beta_m) = \mathbb{E} (\beta_1) \liminf_{m\rightarrow \infty} \mathbb{E} (\alpha_m)
\ge \delta \mathbb{E}(\beta_1)>0\ .
$$ 
The bound \eqref{unif-strict-phi} is now established for $\kappa = - \delta \mathrm{E}(\beta_1)$. Applying Lemma
\ref{unif-strict-con} to $\phi_m(\lambda)$, we conclude that
$\lambda \mapsto \bar{\phi}(\lambda)$ is
strictly concave for every ${\bf Q}_1, {\bf Q}_2$ in ${\cal C}_1$
(${\bf Q}_1 \neq {\bf Q}_2$), and so is ${\bf Q}\mapsto \bar{I}({\bf
  Q})$ by Proposition \ref{prop:prop-concavite}.

\subsection{Proof of Proposition \ref{prop:various}}\label{proof:various}
%Recall the following notations:
%$$
%\alpha_m = \frac 1{\left(1+ \sigma^{-2}
%    \| \tilde {\bf Q} \|\,  \| \tilde {\bf H}^H \tilde {\bf H} \|\right)^2}, \qquad 
%\beta_m = \frac 1m \mathrm{Tr}\left[ \frac{\tilde {\bf H}^H \tilde {\bf H}}{\sigma^2}
%      (\tilde {\bf Q}_1 -\tilde {\bf Q}_2) \frac{\tilde {\bf H}^H
%        \tilde {\bf H}}{\sigma^2} (\tilde {\bf Q}_1 -\tilde{\bf Q}_2)
%    \right]\ . 
%$$

\begin{proof}[Proof of (i)] In order to prove that $\lim_m
  \mathrm{var}(\beta_m)=0$, we shall rely on Poincar\'e-Nash
  inequality. We shall use the following decomposition\footnote{Note
    that the notations introduced hereafter slightly differ from those
    introduced in Section \ref{sec:virtual} but this should not
    disturb the reader.}:
$$
\frac{{\bf C}^{\frac 12}}{\sqrt{K+1}} =  {\bf U D}^{\frac 12} {\bf U}^H;\qquad \tilde{\bf C}^{\frac 12} = \tilde {\bf U}
\tilde{\bf D}^{\frac 12} \tilde {\bf U}^H.
$$
In particular, ${\bf H}$ writes
\begin{eqnarray*}
{\bf U}^H {\bf H} \tilde{\bf U} & = & \sqrt{\frac K{K+1}} {\bf U}^H {\bf A} \tilde {\bf U} 
+ {\bf D}^{\frac 12} \frac{ {\bf U}^H {\bf W} \tilde{\bf U}}{\sqrt{t}} \tilde{\bf D}^{\frac 12}\\
&\stackrel{\triangle}= & {\bf B } +  {\bf D}^{\frac 12} \frac{ {\bf X}}{\sqrt{t}} \tilde{\bf D}^{\frac 12}\quad 
\stackrel{\triangle}= \quad {\bf B} +{\bf Y} \\
& \stackrel{\triangle}=& {\bf \Sigma}\ ,
\end{eqnarray*}
where ${\bf X}$ is a $r\times t$ matrix with i.i.d. ${\mathcal C}N(0,1)$ entries.
Consider now the following matrices:
$$
\check{\bf B} = {\bf I}_m \otimes {\bf B} ,\quad 
{\bf \Gamma}= {\bf I}_m \otimes {\bf D},\quad \tilde {\bf \Gamma}= {\bf I}_m \otimes \tilde {\bf D},\quad
{\bf V}= {\bf I}_m \otimes {\bf U},\quad \tilde {\bf V}= {\bf I}_m \otimes \tilde {\bf U}.
$$
Similarly, $\check{\bf H}$ writes:
$$
{\bf V}^H\check{\bf H} \tilde{\bf V} 
=  \check {\bf B } 
+  {\bf \Gamma}^{\frac 12} 
\frac{ \check{\bf X}}{\sqrt{m t}} \tilde{\bf \Gamma}^{\frac 12}\quad 
\stackrel{\triangle}= \quad \check{\bf B} +\check{\bf Y}\ \quad 
\stackrel{\triangle}= \quad \check{\bf \Sigma}\ ,
$$
where $\check{\bf X}$ is a $mr\times mt$ matrix with i.i.d. ${\mathcal
  C}N(0,1)$ entries.  Denote by ${\bf \Theta} = \tilde{\bf U}^H
({\bf Q}_1 - {\bf Q}_2) \tilde {\bf U}$ and by $\check{\bf \Theta} = \tilde{\bf V}^H
(\check{\bf Q}_1 - \check{\bf Q}_2) \tilde {\bf V} (= {\bf I}_m \otimes {\bf \Theta}) $. The quantity
$\beta_m$ writes then: 
$ \beta_m= \frac 1{\sigma^4 m} \mathrm{Tr}\,
\check{\bf \Theta} \check{\bf \Sigma}^H\check{\bf \Sigma} \check{\bf
  \Theta} \check{\bf \Sigma}^H\check{\bf \Sigma} $. 
Considering
$\beta_m$ as a function of the entries of $\check{\bf
  X}=(\check{X}_{ij})$, i.e. $\beta_m=\phi(\check{\bf X})$, standard
computations yield
$$
\frac{\partial \phi(\check{\bf X})}{\partial \check{X}_{ij}} =
\frac{2}{m} \left( \check{\bf \Theta} \check{\bf \Sigma}^H\check{\bf
    \Sigma} \check{\bf \Theta} \check{\bf \Sigma}^H\right)_{ji}\ .
$$
Poincar\'e-Nash inequality yields then 
\begin{eqnarray*}
\mathrm{var}\,(\beta_m) &\le & \frac 1{mt} \sum_{i,j} {\bf \Gamma}_i \tilde{\bf \Gamma}_j
\mathbb{E} \left| \frac{\partial \phi(\check{\bf X})}{\partial \check{X}_{ij}} \right|^2\\
& = & \frac 1{mt} \sum_{i,j} {\bf \Gamma}_i \tilde{\bf \Gamma}_j \frac{4}{m^2 t^2}
\mathbb{E} \left|\left( \check{\bf \Theta} \check{\bf \Sigma}^H\check{\bf
    \Sigma} \check{\bf \Theta} \check{\bf \Sigma}^H\right)_{ji} \right|^2 \\
&\le & \frac{4 d_{\max} \tilde d_{\max} }{m^3 t^3} 
\mathbb{E} 
\mathrm{Tr}\, \left( \check{\bf \Theta} \check{\bf \Sigma}^H\check{\bf
    \Sigma} \check{\bf \Theta} \check{\bf \Sigma}^H \check{\bf \Sigma}\check{\bf \Theta}^H \check{\bf \Sigma}^H
\check{\bf \Sigma}\check{\bf \Theta}^H
\right)\\
&\le & \frac{4 d_{\max} \tilde d_{\max} }{m^2 t^2} \,  \| \check{\bf \Theta}^H \check{\bf \Theta} \| 
\mathbb{E} \left( \frac{1}{mt} \mathrm{Tr}  \check{\bf \Sigma}^H \check{\bf \Sigma} \check{\bf \Theta} 
 \check{\bf \Sigma}^H \check{\bf \Sigma} \check{\bf \Theta}^{H}  \check{\bf \Sigma}^H \check{\bf \Sigma} \right)\ .
\end{eqnarray*}
Moreover, Schwarz inequality yields
\[
\frac{1}{mt} \mathrm{Tr}  \check{\bf \Sigma}^H \check{\bf \Sigma} \check{\bf \Theta}^H 
(\check{\bf \Sigma}^H \check{\bf \Sigma})^{2} \check{\bf \Theta} \leq 
\left[ \frac{1}{mt} \mathrm{Tr}( \check{\bf \Sigma}^H \check{\bf \Sigma})^{2}\right]^{1/2}
\left[ \frac{1}{mt} \mathrm{Tr}\left(  \check{\bf \Theta}^H  (\check{\bf \Sigma}^H \check{\bf \Sigma})^{2} 
\check{\bf \Theta} \check{\bf \Theta}^H  (\check{\bf \Sigma}^H  \check{\bf \Sigma})^{2} \check{\bf \Theta} \right)
\right]^{1/2}
\]
so that
\[
\frac{1}{mt} \mathrm{Tr}  \check{\bf \Sigma}^H \check{\bf \Sigma} \check{\bf \Theta}^H 
\check{\bf \Sigma}^H \check{\bf \Sigma} \check{\bf \Theta} \leq 
 \| \check{\bf \Theta}^H \check{\bf \Theta} \| \left[ \frac{1}{mt} \mathrm{Tr}( \check{\bf \Sigma}^H \check{\bf \Sigma})^{2}\right]^{1/2} \left[ \frac{1}{mt} \mathrm{Tr}( \check{\bf \Sigma}^H \check{\bf \Sigma})^{4}\right]^{1/2}\ .
\]
Schwarz inequality yields then
\[
\mathbb{E} \left(\frac{1}{mt} \mathrm{Tr}  \check{\bf \Sigma}^H \check{\bf \Sigma} \check{\bf \Theta}^H 
\check{\bf \Sigma}^H \check{\bf \Sigma} \check{\bf \Theta} \right) \leq 
 \| \check{\bf \Theta}^H \check{\bf \Theta} \| \left[ \mathbb{E} \left( \frac{1}{mt} \mathrm{Tr}  (\check{\bf \Sigma}^H \check{\bf \Sigma})^{2} \right) \right]^{1/2} \left[ \mathbb{E} \left( \frac{1}{mt} \mathrm{Tr}  (\check{\bf \Sigma}^H \check{\bf \Sigma})^{4} \right) \right]^{1/2} \ .  
\]
It is tedious, but straightforward, to check that 
\[
\sup_{m} \mathbb{E} \left( \frac{1}{mt} \mathrm{Tr}  (\check{\bf \Sigma}^H \check{\bf \Sigma})^{2} \right) < +\infty
\]
and 
\[
\sup_{m} \mathbb{E} \left( \frac{1}{mt} \mathrm{Tr}  (\check{\bf \Sigma}^H \check{\bf \Sigma})^{4} \right) < +\infty
\] 
which, in turn, imply that $\mathrm{var}(\beta_m) = O(\frac{1}{m^{2}})$. 
\end{proof}

\begin{proof}[Proof of (ii)]
Write $\mathbb{E}\,\beta_m$ as
\begin{eqnarray*}
  \mathbb{E}\, \beta_m & = & \frac 1{\sigma^4 m} \mathbb{E}\,\mathrm{Tr}\,
  \check{\bf \Sigma}^H\check{\bf \Sigma} \check{\bf
    \Theta} \check{\bf \Sigma}^H\check{\bf \Sigma} \check{\bf \Theta} \\
  &=& \frac 1{\sigma^4 m} \mathbb{E}\,\mathrm{Tr}\,
  \left( \check{\bf B}^H\check{\bf B}
    +\check{\bf B}^H\check{\bf Y}
    +\check{\bf Y}^H\check{\bf B}
    +\check{\bf Y}^H\check{\bf Y}\right)
  \check{\bf \Theta} 
  \left( \check{\bf B}^H\check{\bf B}
    +\check{\bf B}^H\check{\bf Y}
    +\check{\bf Y}^H\check{\bf B}
    +\check{\bf Y}^H\check{\bf Y}\right)
  \check{\bf \Theta}\\
  &\stackrel{(a)}=& \frac 1{\sigma^4 m}\,\mathrm{Tr}\, \check{\bf B}^H\check{\bf B} \check{\bf \Theta} \check{\bf B}^H\check{\bf B} \check{\bf \Theta} 
  \quad +\quad  \frac 1{\sigma^4 m}\,\mathbb{E}\,\mathrm{Tr}\, \check{\bf B}^H\check{\bf B} \check{\bf \Theta} \check{\bf Y}^H\check{\bf Y} \check{\bf \Theta}\\
  && + \frac 1{\sigma^4 m}\,\mathbb{E}\,\mathrm{Tr}\, \check{\bf B}^H\check{\bf Y} \check{\bf \Theta} \check{\bf Y}^H\check{\bf B} \check{\bf \Theta}
\quad +\quad  \frac 1{\sigma^4 m}\,\mathbb{E}\,\mathrm{Tr}\, \check{\bf Y}^H\check{\bf B} \check{\bf \Theta} \check{\bf B}^H\check{\bf Y} \check{\bf \Theta}\\
&& + \frac 1{\sigma^4 m}\,\mathbb{E}\,\mathrm{Tr}\, \check{\bf Y}^H\check{\bf Y} \check{\bf \Theta} \check{\bf B}^H\check{\bf B} \check{\bf \Theta}\quad 
 + \quad \frac 1{\sigma^4 m}\,\mathbb{E}\,\mathrm{Tr}\, \check{\bf Y}^H\check{\bf Y} \check{\bf \Theta} \check{\bf Y}^H\check{\bf Y} \check{\bf \Theta}\ ,\\
\end{eqnarray*}
where $(a)$ follows from the fact that the terms where $\check{\bf Y}$
appears one or three times are readily zero, and so are the terms like
$\mathbb{E}\,\mathrm{Tr}\, \check{\bf B}^H\check{\bf Y} \check{\bf
  \Theta} \check{\bf B}^H\check{\bf Y} \check{\bf \Theta}$. Therefore, it remains to compute the following four terms:
\begin{eqnarray*}
  T_1&\stackrel{\triangle}=& \frac 1m \mathrm{Tr}\, \check{\bf B}^H\check{\bf B} \check{\bf \Theta} \check{\bf B}^H\check{\bf B} \check{\bf \Theta}\ , \\
  T_2&\stackrel{\triangle}=& \frac 1m \mathbb{E}\,\mathrm{Tr}\, \check{\bf B}^H\check{\bf B} \check{\bf \Theta} \check{\bf Y}^H\check{\bf Y} \check{\bf \Theta}\ ,\\
  T_3&\stackrel{\triangle}=&\frac 1m \mathbb{E}\,\mathrm{Tr}\, \check{\bf B}^H\check{\bf Y} \check{\bf \Theta} \check{\bf Y}^H\check{\bf B} \check{\bf \Theta}\ ,\\
  T_4&\stackrel{\triangle}=&\frac 1m \mathbb{E}\,\mathrm{Tr}\, \check{\bf Y}^H\check{\bf Y} \check{\bf \Theta} \check{\bf Y}^H\check{\bf Y} \check{\bf \Theta}\ .
\end{eqnarray*}
Due to the block nature of the matrices involved, $T_1= \mathrm{Tr}\,
{\bf B}^H{\bf B} {\bf \Theta} {\bf B}^H{\bf B} {\bf \Theta}$; in
particular, $T_1$ does not depend on $m$. Let us now compute $T_2$. We
have $T_2 = m^{-1} \mathrm{Tr}\, \check{\bf B}^H\check{\bf B}
\check{\bf \Theta} \mathbb{E} \left( \check{\bf Y}^H\check{\bf
    Y}\right) \check{\bf \Theta}$ and $\mathbb{E} \left( \check{\bf
    Y}^H\check{\bf Y}\right) = (mt)^{-1} \tilde{\bf \Gamma}^{\frac 12}
\mathbb{E} \left(\check{\bf X} {\bf \Gamma} \check{\bf X}\right) 
\tilde{\bf \Gamma}^{\frac 12} = (mt)^{-1} \mathrm{Tr}({\bf \Gamma}) \tilde{\bf \Gamma}$. Therefore, $T_2$ writes:
$$
T_2= \frac 1m \mathrm{Tr}\left({\bf \Gamma} \right)
\frac 1{mt}  \mathrm{Tr}\, \left(\check{\bf B}^H \check{\bf B}\check{\bf \Theta} 
\tilde{\bf \Gamma} \check{\bf \Theta}\right) = \mathrm{Tr}\left({\bf D} \right)
\frac 1t \mathrm{Tr}\, \left({\bf B}^H {\bf B}{\bf \Theta} 
\tilde{\bf D} {\bf \Theta}\right),
$$
and this quantity does not depend on $m$. We now turn to the term
$T_3$. We have $T_3= m^{-1} \mathrm{Tr}\, \check{\bf
  B}^H\mathbb{E}\,\left( \check{\bf Y} \check{\bf \Theta} \check{\bf
    Y}^H\right)\check{\bf B} \check{\bf \Theta}$. The same
computations as before yield $\mathbb{E}\,\left( \check{\bf Y}
  \check{\bf \Theta} \check{\bf Y}^H\right) = (mt)^{-1}
\mathrm{Tr}\left( \tilde{\bf \Gamma}^{\frac 12} \check{\bf \Theta}
  \tilde{\bf \Gamma}^{\frac 12} \right) {\bf \Gamma}$. Therefore $T_3$ writes:
$$
T_3 =\frac 1m \mathrm{Tr} \left( \tilde{\bf \Gamma}^{\frac 12} \check{\bf \Theta}
  \tilde{\bf \Gamma}^{\frac 12} \right) \frac 1{mt} \mathrm{Tr}\, 
\left( \check{\bf B}^H {\bf \Gamma} \check{\bf B}\check{\bf \Theta} \right)
= \mathrm{Tr} \left( \tilde{\bf D}^{\frac 12} {\bf \Theta}
  \tilde{\bf D}^{\frac 12} \right) \frac 1{t} \mathrm{Tr}\, 
\left( {\bf B}^H {\bf D} {\bf B}{\bf \Theta}\right)\ , 
$$
which does not depend on $m$. It remains to compute $T_4
 =  \frac 1m \mathrm{Tr}\left[ \mathbb{E}\left( 
\check{\bf Y}^H\check{\bf Y} \check{\bf \Theta} \check{\bf Y}^H\check{\bf Y}
\right) \check{\bf \Theta} \right]$. 
$$
\mathbb{E} \left( \check{\bf Y}^H\check{\bf Y} \check{\bf \Theta} \check{\bf Y}^H\check{\bf Y}
\right)\quad  =\quad  \frac 1{(mt)^2} \tilde{\bf \Gamma}^{\frac 12}
\mathbb{E} \left(\check{\bf X} {\bf \Gamma} \check{\bf X}
\tilde{\bf \Gamma}^{\frac 12} \check{\bf \Theta} \tilde{\bf \Gamma}^{\frac 12}
\check{\bf X} {\bf \Gamma} \check{\bf X}\right) 
\tilde{\bf \Gamma}^{\frac 12}\ .
$$
Computing the individual terms of matrix $\mathbb{E} \left(\check{\bf X} {\bf \Gamma} \check{\bf X}
\tilde{\bf \Gamma}^{\frac 12} \check{\bf \Theta} \tilde{\bf \Gamma}^{\frac 12}
\check{\bf X} {\bf \Gamma} \check{\bf X}\right)$ yields (denote by ${\bf G}= 
\tilde{\bf \Gamma}^{\frac 12} \check{\bf \Theta} \tilde{\bf \Gamma}^{\frac 12}$ for the sake of simplicity):
\begin{eqnarray*}
\left[ \mathbb{E} \left(\check{\bf X} {\bf \Gamma} \check{\bf X}
{\bf G}
\check{\bf X} {\bf \Gamma} \check{\bf X}\right) \right]_{k\ell} & = & 
\sum_{i_1, j_1, j_2, i_2} \mathbb{E} 
\left( \overline{\check{\bf X}}_{i_1, k} \check{\bf X}_{i_1, j_1} \overline{\check{\bf X}}_{i_2, j_2}
\check{\bf X}_{i_2, \ell}\right) {\bf \Gamma}_{i_1, i_1} {\bf G}_{j_1, j_2} {\bf \Gamma}_{i_2, i_2} \\
&= & \left( \mathrm{Tr}\, {\bf \Gamma}\right)^2 {\bf G}_{k\ell} + \mathrm{Tr} 
\left({\bf \Gamma} ^2\right) \, \mathrm{Tr}\, {\bf G}\ \delta_{k\ell}\ , 
\end{eqnarray*}
where $\delta_{k \ell}$ stands for the Kronecker symbol (i.e.
$\delta_{k \ell}=1$ if $k=\ell$, and 0 otherwise). This yields
$$
\mathbb{E} \left( \check{\bf Y}^H\check{\bf Y} \check{\bf \Theta} \check{\bf Y}^H\check{\bf Y}
\right)\quad = \quad
\frac 1{(mt)^2} \left( \mathrm{Tr}\, {\bf \Gamma}\right)^2 \tilde{\bf \Gamma} \check{\bf \Theta} \tilde{\bf \Gamma}
+ \frac 1{(mt)^2} \mathrm{Tr} 
\left({\bf \Gamma} ^2\right) \, \mathrm{Tr}\left( \check{\bf \Theta} \tilde{\bf \Gamma} \right) \tilde{\bf \Gamma}
$$
and 
\begin{eqnarray*}
T_4& = & \frac 1{t^2} \left( \frac{\mathrm{Tr}{\bf \Gamma}
  }{m}\right)^2 \frac 1m \mathrm{Tr} \left( \tilde{\bf \Gamma}
  \check{\bf \Theta} \tilde{\bf \Gamma}\check{\bf \Theta}\right)
+ \frac 1{t^2} \frac 1{m} \mathrm{Tr}\left({\bf \Gamma} ^2\right)
\frac 1{m^2} \left( \mathrm{Tr} \check{\bf \Theta} \tilde{\bf \Gamma} \right)^2 \\
& =& \frac 1{t^2} \left( \mathrm{Tr}{\bf D}\right)^2 \mathrm{Tr} \left( \tilde{\bf D}
  {\bf \Theta} \tilde{\bf D}{\bf \Theta}\right)
+ \frac 1{t^2}  \mathrm{Tr}\left({\bf D}^2\right)
 \left( \mathrm{Tr} {\bf \Theta} \tilde{\bf D} \right)^2 \ ,
\end{eqnarray*}
which does not depend on $m$. This shows that $\mathbb{E} \beta_m$ does not depend on  
$m$, and thus coincides with $\mathbb{E} \beta_1$. In order to complete the proof of (ii), it remains
to verify that $\mathbb{E} \beta_1 > 0$, or equivalenlty that $\mathbb{E} \beta_1$ is not equal to 0. 
If $\mathbb{E} \beta_1$ was indeed equal to 0, then, matrix
$$
\left({\bf H}^{H} {\bf H} \right)^{1/2} \left( {\bf Q}_1 - {\bf Q}_2 \right) \left({\bf H}^{H} {\bf H} \right)^{1/2}
$$
or equivalently matrix 
$$
{\bf H}^{H} {\bf H} ({\bf Q}_1 - {\bf Q}_2)
$$ 
would be equal to zero almost everywhere. As ${\bf Q}_1 \neq {\bf Q}_2$, it would exist a deterministic 
non zero vector ${\bf x}$ such that ${\bf x}^{H} {\bf H}^{H} {\bf H} {\bf x} = 0$ almost everywhere, 
i.e. ${\bf H} {\bf x} = 0$, or equivalently
\begin{equation}
\label{eq:impossible}
{\bf W} \tilde{{\bf C}}^{1/2} {\bf x} = - \sqrt{K t} {\bf C}^{-1/2} {\bf A} {\bf x}\ .
\end{equation}
As matrix $\tilde{{\bf C}}^{1/2}$ is positive definite, vector  $\tilde{{\bf C}}^{1/2} {\bf x}$ is non zero. 
Relation (\ref{eq:impossible}) leads to a contradiction because the joint distribution 
of the entries of ${\bf W}$ is absolutely continuous. This shows that $\mathbb{E} \beta_1 > 0$. The proof of
(ii) is complete. 
\end{proof}

\begin{proof}[Proof of (iii)]
In order to control $\alpha_m = \frac 1{\left(1+ \sigma^{-2}
    \| \tilde {\bf Q} \|\,  \| \tilde {\bf H}^H \tilde {\bf H} \|\right)^2}$, first notice that 
$\| \tilde {\bf Q} \|= \| {\bf Q}\|$. Now $ \| \tilde {\bf H}^H \tilde {\bf H} \| = \| \tilde {\bf H} \|^2$
and 
$$
\| \tilde {\bf H} \| \le \sqrt{\frac K{K+1}} \|\check{\bf A} \| +\frac1{\sqrt{K+1}} \| {\bf \Delta}^{\frac 12} \| \,
\| \tilde {\bf \Delta}^{\frac 12} \|\, \left\|\frac{\check{\bf W}}{\sqrt{mt}} \right\|\ .  
$$
Now $\|\check{\bf A} \|= \| {\bf A}\|$, $\| {\bf \Delta}^{\frac 12} \| = \|{\bf C}^{\frac 12} \|$ and 
$\| \tilde{\bf \Delta}^{\frac 12} \| = \|\tilde{\bf C}^{\frac 12} \|$. 
The behaviour of the spectral norm of $(mt)^{-\frac 12}\check{\bf W}$
is well-known (see for instance \cite{Sil85, BaiSil98}): $\left\|(mt)^{-\frac
    12} \check{\bf W}\right\| \rightarrow_{m\rightarrow \infty} 1+
\sqrt{1/c}$ almost surely. Therefore, Fatou's lemma yields the
desired result: $ \liminf_m \mathbb{E} \alpha_m \ge \delta> 0 \ , $
and (iii) is proved.
\end{proof}

\section{Proof of Proposition \ref{properties-ibar}, item (i).}

  By \eqref{eq:defG} and \eqref{eq:expreVter},
  $(\kappa, \tilde \kappa, \tQ) \mapsto V(\kappa, \tilde \kappa, \tQ)$
  is differentiable from $\mathbb{R}^+ \times \mathbb{R}^+ \times
  {\cal C}_1$ to $\mathbb{R}$. In order to prove that
  $\bar{I}(\tQ)=V(\delta(\tQ), \tilde \delta(\tQ), \tQ)$ is
  differentiable, it is sufficient to prove the differentiability of
  $\delta, \tilde \delta:{\cal C}_1 \rightarrow \mathbb{R}$. Recall
  that $\delta$ and $\tilde \delta$ are solution of system
  \eqref{eq:canonique} associated with matrix $\tQ$. In order to apply
  the implicit function theorem, which will immediatly yield the
  differentiablity of $\delta$ and $\tilde \delta$ with respect to
  $\tQ$, we must check that:
\begin{enumerate}
\item The function 
$$
(\delta, \tilde \delta, \tQ) \mapsto \Upsilon (\delta, \tilde \delta, \tQ)=
\left( 
\begin{array}{c} 
\delta - f(\delta, \tilde \delta, \tQ)\\ 
\tilde \delta - \tilde f(\delta, \tilde \delta, \tQ)
\end{array}
\right)
$$
is differentiable.
\item The partial jacobian 
$$
D_{(\delta, \tilde \delta)} \Upsilon (\delta, \tilde \delta, \tQ)=
\left( 
\begin{array}{cc}
1 - \frac{\partial f}{\partial \delta} (\delta, \tilde \delta, \tQ) 
& -\frac {\partial f}{\partial \tilde \delta} (\delta, \tilde \delta, \tQ)\\
- \frac{\partial \tilde f}{\partial \delta} (\delta, \tilde \delta, \tQ) 
& 1 -\frac {\partial \tilde f}{\partial \tilde \delta} (\delta, \tilde \delta, \tQ)
\end{array}\right)
$$
is invertible for every $\tQ \in {\cal C}_1$.
\end{enumerate}

In order to check the differentiability of $\Upsilon$, recall the following matrix equality
\begin{equation}\label{eq:magic}
({\bf I}+{\bf UV})^{-1}{\bf U} = {\bf U}({\bf I}+{\bf VU})^{-1}
\end{equation}
which follows from elementary matrix manipulations (cf. \cite[Section 0.7.4]{HorJoh94}).
Applying this equality to ${\bf U}=\tQ^{\frac 12}$ and $V=\delta \tilde {\bf C} \tQ^{\frac 12}$, we obtain:
$$
{\bf A} \tQ^{\frac 12} \left( {\bf I} + \delta \tQ^{\frac 12} \tilde
  {\bf C} \tQ^{\frac 12} \right)^{-1} \tQ^{\frac 12} {\bf A}^H = {\bf A} \tQ
\left( {\bf I} + \delta  \tilde
  {\bf C} \tQ \right)^{-1}{\bf A}^H
$$
which yields
$$
f(\delta, \tilde \delta, {\bf Q}) =\frac 1t \mathrm{Tr} \bigg\{ {\bf
  C} \Big[ \sigma^2\,\big({\bf
  I}_r+  \frac{\tilde\delta}{K+1} \, {\bf C}\,\big)
+ \frac{K}{K+1} {\bf A}{\bf Q}\pa{{\bf I}_t+
  \frac{\delta}{K+1} \,{\bf \tilde{C}}{\bf
    Q}}^{-1}{\bf A}^H \Big]^{-1}\bigg\}\ .
$$
Clearly, $f$ is differentiable with respect to the three variables
$\delta, \tilde \delta$ and $\tQ$. Similar computations yield
$$
\tilde{f}(\delta, \tilde{\delta}, \mb{Q}) =  \f{1}{t}\, \mathrm{Tr}
\bigg\{ {\bf Q}{\bf \tilde{C}}  \Big[\sigma^2\,\big({\bf I}_t+\frac{\delta}{K+1}\tilde{{\bf
C}}{\bf Q} \,\big)\\
+ \frac{K}{K+1} {\bf A}^H\pa{{\bf
I}_r+\frac{\tilde\delta}{K+1}\,{\bf C}}^{-1}{\bf A}{\bf Q}\Big]^{-1}\bigg\}\ ,
$$
and the same conclusion holds for $\tilde f$. Therefore, $(\delta, \tilde \delta, \tQ) \mapsto
\Upsilon(\delta, \tilde \delta, \tQ)$ is differentiable and 1) is proved. 
 
In order to study the jacobian $D_{(\delta, \tilde \delta)} \Upsilon$,
let us compute first $\frac{\partial f}{\partial \delta}$.
\begin{eqnarray*}
\frac{\partial f}{\partial \delta}(\delta, \tilde \delta, {\bf Q}) &=&
\frac 1t \mathrm{Tr}\ {\bf C} {\bf T}_K {\bf A} \tQ \left( {\bf I} +\frac{\delta}{K+1} \tilde {\bf C} {\bf Q}\right)^{-1}
\frac {\tilde {\bf C} {\bf Q}}{K+1} \left( {\bf I} +\frac{\delta}{K+1} \tilde {\bf C} {\bf Q}\right)^{-1} {\bf A}^H
{\bf T}_K \frac{K}{K+1}\ ,\\
&= & \frac 1t \mathrm{Tr}\  {\bf C} {\bf T}_K {\bf A} \tQ^{\frac 12} \left( {\bf I} +\frac{\delta}{K+1} \tQ^{\frac 12}
\tilde {\bf C} {\bf Q}^{\frac 12} \right)^{-1}
\frac {\tQ^{\frac 12} \tilde {\bf C} {\bf Q}^{\frac 12} }{K+1} \\
&\phantom{=}&\quad \times  \frac {\tQ^{\frac 12} \tilde {\bf C} {\bf Q}^{\frac 12} }{K+1} \left( {\bf I} +\frac{\delta}{K+1} \tQ^{\frac 12}
\tilde {\bf C} {\bf Q}^{\frac 12} \right)^{-1}
\tQ^{\frac 12} {\bf A}^H {\bf T}_K \frac{K}{K+1}\ ,\\
&\stackrel{(a)}{=} & \frac{1}{t} \mathrm{Tr}\ ( {\bf D} {\bf T} {\bf B} 
({\bf I} + \beta \tilde{{\bf D}})^{-1} \tilde{{\bf D}} ({\bf I} + \beta \tilde{{\bf D}})^{-1} {\bf B}^{H} {\bf T})
\end{eqnarray*}
where $(a)$ follows from the virtual channel equivalences
\eqref{eq:eigen}, \eqref{eq:definition-B} together with
\eqref{eq:lien-delta-beta} and \eqref{equiv-tk}. Finally, we end up
with the following:
$$
1-\frac{\partial f}{\partial \delta}(\delta, \tilde \delta, {\bf Q}) = 1-
\frac{1}{t} \mathrm{Tr}( {\bf D} {\bf T} {\bf B} ({\bf I} + \beta
\tilde{{\bf D}})^{-1} \tilde{{\bf D}} ({\bf I} + \beta \tilde{{\bf
    D}})^{-1} {\bf B}^{H} {\bf T})\ .
$$
Similar computations yield
\begin{eqnarray*}
1-\frac{\partial \tilde f}{\partial \tilde \delta}(\delta, \tilde \delta, {\bf Q}) &=& 1-
\frac{1}{t} \mathrm{Tr}( \tilde {\bf D} \tilde {\bf T} {\bf B}^H ({\bf I} + \tilde \beta
{{\bf D}})^{-1} {{\bf D}} ({\bf I} + \tilde \beta {{\bf
    D}})^{-1} {\bf B} \tilde {\bf T})\ ,\\
- \frac{\partial f}{\partial \tilde \delta} (\delta, \tilde \delta, {\bf Q}) & = & 
\frac{\sigma^2}t \mathrm{Tr}\ ({\bf DTDT})\ ,\\
- \frac{\partial \tilde f}{\partial \delta} (\delta, \tilde \delta, {\bf Q}) & = & 
\frac{\sigma^2}t \mathrm{Tr}\ ( \tilde {\bf D}\tilde {\bf T}\tilde {\bf D}\tilde {\bf T})\ .
\end{eqnarray*}
The invertibility of the jacobian $D_{(\delta, \tilde
  \delta)}\Upsilon$ follows then from Lemma \ref{le:determinant} in
Appendix \ref{subsec:sketch} and 2) is proved. In particular, we can
assert that ${\cal C}_1 \ni \tQ \mapsto \delta(\tQ)$ and ${\cal C}_1
\ni \tQ \mapsto \tilde \delta(\tQ)$ are differentiable due to the
Implicit function theorem. Item (i) is proved.

\section{Proof of Proposition \ref{prop:convergence-algo}}

First note that the sequence $(\mb{Q}_{k})$ belongs to the compact set
${\cal C}_1$. Therefore, in order to show that the sequence converges,
it is sufficient to establish that the limits of all convergent
subsequences coincide.  We thus consider a convergent subsequence
extracted from $(\mb{Q}_{k})_{k \geq 0}$, say $(\mb{Q}_{\psi(k)})_{k
  \geq 0}$, where for each $k$, $\psi(k)$ is an integer, and denote by
$\mb{Q}_{*}^{\psi}$ its limit. If we prove that
\begin{equation}
\label{eq:cqfd} <\nabla \bar{I}(\mb{Q}_{*}^{\psi}), \mb{Q} -
\mb{Q}_{*}^{\psi}> \; \leq 0
\end{equation}
for each $\mb{Q} \in {\cal C}_1$, Proposition
\ref{properties-ibar}-(ii) will imply that $\mb{Q}_{*}^{\psi}$
coincides with the argmax $\tO$ of $\bar{I}$ over ${\cal C}_1$.
This will prove that the limit of every convergent subsequence
converges towards $\tO$, which in turn will show that the whole
sequence $(\mb{Q}_{k})_{k \geq 0}$ converges to $\tO$.

In order to prove (\ref{eq:cqfd}), consider the iteration $\psi(k)$ of
the algorithm. The matrix $\mb{Q}_{\psi(k)}$ maximizes the function
$\mb{Q} \mapsto V(\delta_{\psi(k)}, \tilde{\delta}_{\psi(k)},
\mb{Q})$. As this function is strictly concave andd differentiable,
Proposition 
\ref{prop:caracterisation} implies that
\begin{equation}
\label{eq:etape1} <\nabla_{\tQ} V (\delta_{\psi(k)}, \tilde{\delta}_{\psi(k)},
\mb{Q}_{\psi(k)}), \mb{Q} -  \mb{Q}_{\psi(k)}> \; \leq 0
\end{equation}
for every $\mb{Q} \in {\cal C}_1$ (recall that $\nabla_{\tQ}$
represents the derivative of $V(\kappa, \tilde \kappa, \tQ)$ with
respect to $V$'s third component).  We now consider the pair of
solutions $(\delta_{\psi(k)+1}, \tilde{\delta}_{\psi(k)+1})$ of the
system (\ref{eq:canonique}) associated with matrix $\mb{Q}_{\psi(k)}$.

Due to the continuity of $\delta(\tQ)$ and $\tilde \delta(\tQ)$, the
convergence of the subsequence $\mb{Q}_{\psi(k)}$ implies the
convergence of the subsequences $(\delta_{\psi(k)+1},
\tilde{\delta}_{\psi(k)+1})$ towards a limit $(\delta_{*}^{\psi},
\tilde{\delta}_{*}^{\psi})$.  The pair $(\delta_{*}^{\psi},
\tilde{\delta}_{*}^{\psi})$ is the solution of system
(\ref{eq:canonique}) associated with $\mb{Q}_{*}^{\psi}$ i.e.
$\delta_{*}^{\psi} = \delta(\mb{Q}_{*}^{\psi})$ and
$\tilde{\delta}_{*}^{\psi} = \tilde{\delta}(\mb{Q}_{*}^{\psi})$; in particular:
$$
\frac{\partial V}{\partial \kappa} (\delta_*^\psi, \tilde{\delta}_*^\psi, {\bf Q}_*^\psi) =  
\frac{\partial V}{\partial \tilde{\kappa}} (\delta_*^\psi, \tilde{\delta}_*^\psi, {\bf Q}_*^\psi)=0\ 
$$
(see for instance \eqref{eq:crucial}). Using the same computation as in the proof of Proposition
\ref{prop:waterfilling}, we obtain
\begin{equation}
\label{eq:egalite-bis} \langle \nabla \bar{I}(\mb{Q}_{*}^{\psi}), \mb{Q} -
\mb{Q}_{*}^{\psi}\rangle  = \langle \nabla V\left( \delta_{*}^{\psi},
\tilde{\delta}_{*}^{\psi}, \mb{Q}_{*}^{\psi} \right),  \mb{Q} -
\mb{Q}_{*}^{\psi} \rangle
\end{equation}
for every $\mb{Q} \in {\cal C}_1$. Now condition (\ref{eq:condition})
implies that the subsequence $(\delta_{\psi(k)},
\tilde{\delta}_{\psi(k)})$ also converges toward
$(\delta_{*}^{\psi}, \tilde{\delta}_{*}^{\psi})$. As a consequence,
$$
\lim_{k \rightarrow +\infty} \langle \nabla V(\delta_{\psi(k)},
\tilde{\delta}_{\psi(k)}, \mb{Q}_{\psi(k)}), \mb{Q} -
\mb{Q}_{\psi(k)}\rangle = \langle \nabla V(\delta_{*}^{\psi},
\tilde{\delta}_{*}^{\psi}, \mb{Q}_{*}^{\psi}), \mb{Q} -
\mb{Q}_{*}^{\psi}\rangle \ .
$$
Inequality (\ref{eq:etape1}) thus implies that
$
\langle  \nabla V(\delta_{*}^{\psi},
\tilde{\delta}_{*}^{\psi},\mb{Q}_{*}^{\psi}), \mb{Q} -
\mb{Q}_{*}^{\psi}\rangle   \leq 0 
$ and relation (\ref{eq:egalite-bis}) allows us to conclude the proof.

\section{End of proof of Proposition \ref{approx-capacity}}\label{proof-approx-capacity}
Proof of Proposition \ref{approx-capacity} relies on properties of $\tO$ established in Proposition 
\ref{properties-ibar}--(iii). Denote by 
$$
A = \max\left( \sup_t \| {\bf A} \| , 
\sup_t \| \tilde{\bf C} \| , \sup_t \| {\bf C} \| \right) < \infty\quad \textrm{and} \quad
a = \min\left( \inf_t \lambda_{\mathrm{min}} ( \tilde{\bf C} ), 
\inf_t \lambda_{\mathrm{min}} ( {\bf C} )\right) > 0\ . 
$$
\paragraph*{Proof of (i)} 
Recall that by Proposition
\ref{properties-ibar}--(iii), $\tO$ maximizes 
$\log\det({\bf I} + {\bf QG}(\delta_*, \tilde{\delta}_*))$. This implies that 
the eigenvalues $( \lambda_j(\tO) )$ are the solutions of the waterfilling equation
$$
\forall \ j=1,\ldots, t, \quad \lambda_j(\tO) = \max\left( \gamma - 
\frac{1}{\lambda_j({\bf G})}, 0 \right) 
$$
where $\gamma$ is tuned in such a way that $\sum_j \lambda_j(\tO) =
t$. It is clear from this equation that $\| \tO \| \leq \gamma$. If 
$\gamma \leq \lambda_{\min}({\bf
  G})^{-1}$ then $\|\tO \| \le \lambda_{\min}({\bf
  G})^{-1}$. If $\gamma \geq \lambda_{\min}({\bf G})^{-1}$ then 
$\gamma \ge \lambda_{j}({\bf G})^{-1}$ and we have:
$$
t = \sum_j \lambda_j(\tO) = \gamma t - \sum_j \frac{1}{\lambda_j({\bf G})}\ ,
$$ 
hence 
$$
\gamma = 1 + \frac 1t \sum_j \frac{1}{\lambda_j({\bf G})} \leq 1 +\frac{1}{\lambda_{\min}({\bf G})}\ .
$$
In both cases, we have 
\begin{equation}\label{majo-gamma}
\| \tO\| \le  1 +\frac{1}{\lambda_{\min}({\bf G})}\ .
\end{equation}
It remains to prove 
\begin{equation}
\label{eq-lambda-min(G)}
\forall \ {\bf Q} \in {\cal C}_1, \quad 
\inf_t \lambda_{\mathrm{min}}\left( {\bf G}(\delta({\bf Q}), \tilde{\delta}({\bf Q})) \right)
> 0 
\end{equation}
and we are done. To this end, we first show that $\inf_t \delta({\bf Q}) > 0$ for all 
${\bf Q} \in {\cal C}_1$. From Equations \eqref{eq:defTK} and \eqref{eq:expredelta}, we
have:
\begin{eqnarray}
\delta({\bf Q}) &=& \frac 1t \tr {\bf C} {\bf T}_K(\sigma^2) \nonumber \\
&\geq& \lambda_{\min}({\bf C}) \frac 1t \tr {\bf T}_K(\sigma^2) \nonumber \\
&\stackrel{(a)}{\geq}& \lambda_{\min}({\bf C}) 
\left[
\frac 1t \tr\left(
\sigma^2 {\bf I}_r + \frac{\sigma^2}{K+1} \tilde{\delta} {\bf C} \right. \right. \nonumber \\ 
& & \left. \left. + \frac{K}{K+1} {\bf AQ}^{1/2} 
\left( {\bf I}_t + \frac{\delta}{K+1} {\bf Q}^{1/2} \tilde{\bf C} {\bf Q}^{1/2} \right)^{-1} 
{\bf Q}^{1/2} {\bf A}^H \right) \right]^{-1} \nonumber \\
&\stackrel{(b)}{\geq}& 
 \lambda_{\min}({\bf C}) 
\left(
\frac 1t \tr\left(
\sigma^2 {\bf I}_r + \frac{\sigma^2}{K+1} \tilde{\delta} {\bf C}  
 + \frac{K}{K+1} {\bf AQA}^{H} \right) \right)^{-1} 
\label{eq-borneinf-delta} 
\end{eqnarray} 
where $(a)$ follows from Jensen's Inequality and $(b)$ is due to the facts that 
$\| ({\bf I}_t + {\bf Y} )^{-1} \| \leq 1$ and 
$\tr( {\bf XY}) \leq \| {\bf X} \| \tr({\bf Y})$ when ${\bf Y}$ is a nonnegative matrix. 
We now find an upper bound for $\tilde\delta$. From \eqref{equiv-tk} and (\ref{eq:borne-T}), we
have $\| \tilde {\bf T}_K(\sigma^2) \| \leq 1/\sigma^2$. Using \eqref{eq:expredelta} we
then have
$$
\tilde{\delta} \leq \| \tilde {\bf T}_K \| \frac 1t \tr \tilde{\bf C} {\bf Q} 
\leq \| \tilde {\bf T}_K \| \| \tilde{\bf C} \| \frac 1t \tr {\bf Q} 
\leq \frac{A}{\sigma^2} \ 
$$
(recall that $\frac 1t \tr \,{\bf Q} = 1$). Getting back to
\eqref{eq-borneinf-delta}, we easily obtain   
$$
\frac 1t \tr\left( \sigma^2 {\bf I}_r + \frac{\sigma^2}{K+1}
  \tilde{\delta} {\bf C} + \frac{K}{K+1} {\bf AQA}^{H} \right) \leq
\frac rt \left( \sigma^2 + \frac{A}{K+1} \right) + \frac{A^2 K}{K+1}
\le C_0 \quad \forall (t,r),\ \frac tr \rightarrow c\ 
$$
where $C_0$ is a certain constant term. Hence we have $\delta({\bf Q}) \geq a C_0^{-1} $. By inspecting the expression 
\eqref{eq:defG} of ${\bf G}(\delta, \tilde\delta)$, we then obtain 
$$
\lambda_{\min}({\bf G}) \geq \frac{a C_0^{-1}}{K+1} \lambda_{\min}(\tilde{\bf C}) 
\geq \frac{a^2 C_0^{-1}}{K+1} = C_1 > 0 
$$
and \eqref{eq-lambda-min(G)} is proven. It remains to plug this estimate into \eqref{majo-gamma}
and (i) is proved.
 
\paragraph*{Proof of (ii)} 
We begin by restricting the maximization of $I(\mb{Q})$ to the set 
${\cal C}^{\mathrm{d}}_1 = \{ {\bf Q} \ : \ 
{\bf Q} = \diag(q_1, \ldots, q_t) \geq {\bf 0}, 
\tr({\bf Q}) = t \}$ of the diagonal matrices within 
${\cal C}_1$, and show that
${\bf Q}^{\mathrm{d}}_* = \arg\max_{{\bf Q} \in {\cal C}^{\mathrm{d}}_1} 
I({\bf Q})$ satisfies $\sup_t \| \mb{Q}_*^{\mathrm{d}} \| < \infty$ where the bound
is a function of $(a,A,\sigma^2, c,K)$ only. 
The set ${\cal C}^{\mathrm{d}}_1$ is clearly convex and the solution ${\bf Q}^{\mathrm{d}}_*$
is given by the Lagrange Karush-Kuhn-Tucker (KKT) conditions
\begin{equation} 
\label{eq-KKT-I(Q)}
\frac{\partial I(\mb{Q})}{\partial q_j} = 
\frac{\partial}{\partial q_j} 
\EE\left[ {\cal I}({\bf Q}) \right] = 
\eta - \beta_j 
\end{equation} 
where ${\cal I}({\bf Q}) = \log\det\left( {\bf I}_r + 
\frac{1}{\sigma^2} {\bf HQH}^H \right)$ and the Lagrange multipliers 
$\eta$ and the $\beta_i$ are associated with the power constraint 
and with the positivity constraints respectively. 
More specifically, $\eta$ is the unique real positive number for which 
$\sum_{j=1}^t q_j = t$, and the $\beta_j$ satisfy 
$\beta_j = 0$ if $q_j > 0$ and $\beta_j \geq 0$ if $q_j = 0$. 
We have
$$
\frac{\partial {\cal I}({\bf Q})}{\partial q_j} 
= 
\frac{1}{\sigma^2} {\bf h}_j^H 
\left( {\bf I}_r + \frac{1}{\sigma^2} {\bf HQH}^H \right)^{-1} {\bf h}_j 
$$
where ${\bf h}_j$ the $j^{\mathrm{th}}$ column of ${\bf H}$. By consequence, 
$
\EE \left[ \partial {\cal I}({\bf Q}) / \partial q_j \right] 
\leq
\frac{1}{\sigma^2}  \EE \left[ \| {\bf h}_j \|^2 \right] 
$. As ${\bf h}_j$ is a Gaussian vector, the righthand side of this inequality is defined
and therefore, by the Dominated Convergence Theorem, 
we can exchange $\partial/\partial q_j$ with $\EE$ in Equation 
\eqref{eq-KKT-I(Q)} and write
\begin{equation}
\label{eq-derivee-KKT-I(Q)} 
\frac{\partial I(\mb{Q})}{\partial q_j} = 
\frac{1}{\sigma^2} \EE\left[ {\bf h}_j^H 
\left( {\bf I}_r + \frac{1}{\sigma^2} {\bf HQH}^H \right)^{-1} 
{\bf h}_j \right] 
\end{equation} 
Let us denote by ${\bf H}_j$ the $r \times (t-1)$ matrix that remains after extracting 
${\bf h}_j$ from ${\bf H}$.  Similarly, we denote by ${\bf Q}_j$ the $(t-1) \times
(t-1)$ diagonal matrix that remains after deleting row and column $j$ from ${\bf Q}$.
Writing ${\bf R}_j = 
\left( {\bf I}_r + \frac{1}{\sigma^2} {\bf H}_j {\bf Q}_j {\bf H}_j^H \right)^{-1}$, we have
by the Matrix Inversion Lemma (\cite[\S 0.7.4]{HorJoh94})
$$
\left( {\bf I}_r + \frac{1}{\sigma^2} {\bf HQH}^H \right)^{-1} = 
{\bf R}_j - 
\frac{q_j}{\sigma^2 + q_j {\bf h}_j^H {\bf R}_j {\bf h}_j} 
{\bf R}_j {\bf h}_j {\bf h}_j^H {\bf R}_j \ .
$$
By plugging this expression into the righthand side of Equation 
\eqref{eq-derivee-KKT-I(Q)}, the Lagrange-KKT conditions become
\begin{equation} 
\label{eq-kkt-I(Q)-fraction} 
\EE\left[ \frac{X_j}{\sigma^2 + q_j X_j} \right] = \eta - \beta_j 
\end{equation} 
where $X_j = {\bf h}_j^H {\bf R}_j {\bf h}_j$. 
A consequence of this last equation is that $q_j \leq 1 / \eta$ for every $j$. 
Indeed, assume that $q_j > 1 / \eta$ for some
$j$. Then $\sigma^2 + q_j X_j >  X_j / \eta$ hence
$\EE\left[ \frac{X_j}{\sigma^2 + q_j X_j} \right] < \eta$, therefore $\beta_j > 0$
\eqref{eq-kkt-I(Q)-fraction}, which implies that
$q_j = 0$, a contradiction. As a result, in order to prove that 
$\sup_t \| \mb{Q}_*^{\mathrm{d}} \| < \infty$, it will be enough to prove
that $\sup_t 1/ \eta < \infty$. To this end, we shall prove that there exists 
a constant $C > 0$ such that 
\begin{equation}
\label{eq-proba-Xj} 
\max_{j=1,\ldots,t} \PP\left( X_j \leq C \right) 
\xrightarrow[t\to\infty]{ } 0 \ .
\end{equation} 
Indeed, let us admit \eqref{eq-proba-Xj} temporarily. We have
\begin{eqnarray*}
\EE\left[ \frac{X_j}{\sigma^2 + q_j X_j} \right] 
- \frac{C}{\sigma^2 + q_j C} 
&=&
\EE\left[ \frac{X_j}{\sigma^2 + q_j X_j} {\bf 1}_{X_j > C} \right] 
- \frac{C}{\sigma^2 + q_j C} 
+ \EE\left[ \frac{X_j}{\sigma^2 + q_j X_j} {\bf 1}_{X_j \leq C} \right] \\
&\geq& 
\frac{C}{\sigma^2 + q_j C} \PP(X_j > C) 
- \frac{C}{\sigma^2 + q_j C}
\\ 
&=& 
\varepsilon_{j} 
\end{eqnarray*} 
where $\varepsilon_{j} = - \frac{C}{\sigma^2 + q_j C} \PP(X_j \leq C)$, 
and the inequality is due to the fact that the function 
$f(x) = \frac{x}{\sigma^2 + q_j x}$ is increasing. As 
\[
\max_{j=1,\ldots,t} | \varepsilon_{j} | \leq
\frac{C}{\sigma^2} \max_{j=1,\ldots,t} \PP(X_j \leq C)  
\xrightarrow[t\to\infty]{ } 0
\]
by \eqref{eq-proba-Xj}, we have 
\[
\liminf_t \min_j 
\left( \EE\left[ \frac{X_j}{\sigma^2 + q_j X_j} \right] 
-
\frac{C}{\sigma^2 + q_j C} \right) 
\geq 0 \ .  
\]
Getting back to the Lagrange KKT condition \eqref{eq-kkt-I(Q)-fraction}
we therefore have for $t$ large enough 
$\eta - \beta_j > \frac{C/2}{\sigma^2 + q_j C/2}$ for every 
$j=1,\ldots,t$. By consequence, 
$$
\frac 1 \eta \leq \frac{1}{\eta - \beta_j} < \frac{2 \sigma^2}{C} + q_j 
$$
for large $t$.
Summing over $j$ and taking into account the power constraint 
$\sum_j q_j = t$, we obtain 
$ \frac t \eta < \frac{2 \sigma^2 t}{C} + t $, i.e. $ \frac{1}{\eta} < \frac{2 \sigma^2 }{C} + 1 $ and
\begin{equation} 
\label{eq-borne-Q} 
\sup_t \| \mb{Q}_*^{\mathrm{d}} \| <  \frac{2 \sigma^2}{C} + 1 
\end{equation} 
which is the desired result. 
To prove \eqref{eq-proba-Xj}, we make use of MMSE estimation theory. 
Recall that ${\bf H} = \sqrt{\frac{K}{K+1}} {\bf A} + 
\frac{1}{\sqrt{K+1}} \frac{1}{\sqrt{t}} {\bf C}^{1/2} {\bf W} \tilde{\bf C}^{1/2}$.  
Denoting by ${\bf a}_j$ and ${\bf z}_j$ the $j^{\mathrm{th}}$ columns of the matrices 
${\bf A}$ and ${\bf W} \tilde{\bf C}^{1/2}$ respectively, we have
$$
X_j = \left( \sqrt{\frac{K}{K+1}} {\bf a}_j^H + 
\frac{1}{\sqrt{K+1}} \frac{1}{\sqrt{t}} {\bf z}_j^H {\bf C}^{1/2} \right)
{\bf R}_j 
\left( \sqrt{\frac{K}{K+1}} {\bf a}_j + 
\frac{1}{\sqrt{K+1}} \frac{1}{\sqrt{t}} {\bf C}^{1/2} {\bf z}_j \right)\ .
$$
We decompose ${\bf z}_j$ as ${\bf z}_j = {\bf u}_{j} + {\bf u}_{j}^\perp$ where
${\bf u}_j$ is the conditional expectation ${\bf u}_j = \EE\left[ {\bf z}_j \| 
{\bf z}_1, \ldots, {\bf z}_{j-1}, {\bf z}_{j+1}, \ldots, {\bf z}_t \right]$, in
other words, ${\bf u}_j$ is the MMSE estimate of ${\bf z}_j$ drawn from the other columns
of ${\bf W} \tilde{\bf C}^{1/2}$. 
Put 
\begin{eqnarray}
S_j &=& 2 \Re \left( 
\frac{1}{\sqrt{K+1}} \frac{1}{\sqrt{t}} {{\bf u}_j^\perp}^H {\bf C}^{1/2} 
{\bf R}_j 
\left( \sqrt{\frac{K}{K+1}} {\bf a}_j + 
\frac{1}{\sqrt{K+1}} \frac{1}{\sqrt{t}} {\bf C}^{1/2} {\bf u}_j \right) \right) 
\nonumber \\
& & + 
\frac{1}{t(K+1)} 
{{\bf u}_j^\perp}^H {\bf C}^{1/2} {\bf R}_j {\bf C}^{1/2} {\bf u}_j^\perp \ .
\label{eq-expression-Sj} 
\end{eqnarray} 
Then 
\begin{eqnarray}
X_j &=& S_j + \left( \sqrt{\frac{K}{K+1}} {\bf a}_j^H + 
\frac{1}{\sqrt{K+1}} \frac{1}{\sqrt{t}} {\bf u}_j^H {\bf C}^{1/2} \right)
{\bf R}_j 
\left( \sqrt{\frac{K}{K+1}} {\bf a}_j + 
\frac{1}{\sqrt{K+1}} \frac{1}{\sqrt{t}} {\bf C}^{1/2} {\bf u}_j \right) \nonumber \\
&\geq& S_j \ .
\label{eq-X_j>S_j} 
\end{eqnarray} 
Let us study the asymptotic behaviour of $S_j$. First, we note that due  
to the fact that the joint distribution of the 
elements of ${\bf W} \tilde{\bf C}^{1/2}$ is the Gaussian distribution, ${\bf u}_{j}^\perp$
and ${\bf v}_j = [ {\bf z}_1^T, \ldots, {\bf z}_{j-1}^T, {\bf z}_{j+1}^T, 
\ldots, {\bf z}_t^T ]^T$ are independent. By consequence, 
${\bf u}_j^\perp$ and $({\bf R}_j, {\bf u}_j)$ are independent. 
Let us derive the expression of the covariance matrix ${\bf R}_{\bf u} = 
\EE[ {\bf u}_j^\perp {{\bf u}_j^\perp}^H ]$. 
From the well known formulas for MMSE estimation (\cite{scharf91}), 
we have
${\bf R}_{\bf u} = 
\EE[ {\bf z}_j {\bf z}_j^H ] - \EE[ {\bf z}_j {\bf v}_j^H ] 
\left(\EE[ {\bf v}_j {\bf v}_j^H ] \right)^{-1} \EE[ {\bf v}_j {\bf z}_j^H ]$. 
To obtain ${\bf R}_{\bf u}$, we note that the 
covariance matrix of the vector ${\bf z} = [ {\bf z}_1^T, \ldots, {\bf z}_t^T ]^T$
is $\EE[ {\bf z} {\bf z}^H ] = \tilde{\bf C}^T \otimes {\bf I}_r$ (just check that
$\EE\left[ [ {\bf W} \tilde{\bf C}^{1/2} ]_{ij} 
\overline{[ {\bf W} \tilde{\bf C}^{1/2} ]_{kl}} \right] = \delta(i-k) 
[ \tilde{\bf C} ]_{lj}$). Let us denote by $\tilde c_j$, $\tilde{\bf c}_j$ and 
$\tilde{\bf C}_j$ the scalar $\tilde c_{j} = [ \tilde{\bf C} ]_{jj}$, the $j^{\mathrm{th}}$
vector column of $\tilde{\bf C}$ without element $\tilde c_{j}$, and the $(t-1) \times 
(t-1)$ matrix that remains after extracting row and column $j$ from $\tilde{\bf C}$
respectively. With these notations we have
${\bf R}_{\bf u} = 
\left( \tilde c_j - \tilde{\bf c}_j^H \tilde{\bf C}_j^{-1} \tilde{\bf c}_j \right) 
{\bf I}_r$. Recalling that ${\bf u}_j^\perp$ and $({\bf R}_j, {\bf u}_j)$ are independent,  
% and the elements of ${\bf u}_j^\perp$ are centered, 
one may see that the first term of the righthand side of \eqref{eq-expression-Sj} is
negligible while the second is close to 
$\rho_j = 
\frac{1}{t} \frac{\tilde c_j - \tilde{\bf c}_j^H \tilde{\bf C}_j^{-1} \tilde{\bf c}_j}{K+1} 
\tr( {\bf R}_j {\bf C} )$. More rigorously, using this independence 
% of ${\bf u}_j^\perp$ and $({\bf R}_j, {\bf u}_j)$ 
in addition to $A = \max( \| {\bf A} \|, \| {\bf C} \|, \| \tilde{\bf C} \| ) < \infty$ 
and $\| {\bf R}_j \| \leq 1$, 
we can prove with the help of \cite[Lemma 2.7]{BaiSil98} or by direct 
calculation that there exists a constant $C_1$ such that 
\begin{equation}
\label{eq-tchebychev} 
\EE \left[ \left( S_j - \rho_j \right)^2 \right] \leq \frac{C_1}{t} \ .
\end{equation} 
In order to prove \eqref{eq-proba-Xj}, we will prove that the $\rho_j$
are bounded away from zero in some sense. First, we have
$$
\tilde c_j - \tilde{\bf c}_j^H \tilde{\bf C}_j^{-1} \tilde{\bf c}_j 
\stackrel{(a)}{=} 
\left[ \tilde{\bf C}^{-1} \right]_{jj}^{-1} 
\stackrel{(b)}{\geq} 
\| \tilde{\bf C}^{-1} \|^{-1} = \lambda_{\mathrm{min}}(\tilde{\bf C}) 
\geq a 
$$ 
(for $(a)$ see \cite[\S 0.7.3]{HorJoh94} and for $(b)$, use
the fact that $| [{\bf X}]_{kl} | \leq \| {\bf X} \|$ for any element
$(k,l)$ of a matrix ${\bf X}$). 
% Recall that if ${\bf A}$ and ${\bf B}$ are Hermitian nonnegative matrices, then 
% $\tr({\bf AB}) \geq \lambda_{\mathrm{min}}({\bf A}) \tr({\bf B})$. 
By consequence, 
\begin{eqnarray*} 
\rho_j &\geq& \frac{a \lambda_{\mathrm{min}}({\bf C})}{K+1} \frac 1t \tr 
\left( {\bf I}_r + \frac{1}{\sigma^2} {\bf H}_j {\bf Q}_j {\bf H}_j^H \right)^{-1} \\
&\stackrel{(a)}{\geq}& 
\frac{a \lambda_{\mathrm{min}}({\bf C})}{K+1} 
\left( \frac 1t \tr 
\left( {\bf I}_r + \frac{1}{\sigma^2} {\bf H}_j {\bf Q}_j {\bf H}_j^H \right) 
\right)^{-1} \\
&\stackrel{(b)}{\geq}& 
\frac{a^2}{K+1} 
\left(\frac rt + 
\frac{1}{\sigma^2} \left( \| {\bf A} \| + \| {\bf C} \|^{1/2} \| \tilde{\bf C} \|^{1/2}  
\| \frac{1}{\sqrt{t}} {\bf W} \| \right)^{2} \frac 1t \tr ({\bf Q}) \right)^{-1}  
\end{eqnarray*} 
where $(a)$ is Jensen Inequality and $(b)$ is due to 
$\tr( {\bf XY}) \leq \| {\bf X} \| \tr({\bf Y})$ when ${\bf Y}$ is a nonnegative matrix. 
As $\lim_t \| \frac{1}{\sqrt{t}} {\bf W} \| =  1 + \sqrt{1/c}$ with probability one 
(\cite{BaiSil98}), and furthermore, $\tr({\bf Q}) = t$, we have with probability one 
\begin{equation}
\label{eq-liminf-rho} 
\liminf_t \min_{j=1,\ldots,t} \rho_j \geq \frac{a^2}{K+1} 
\left( c^{-1} + \frac{A^2}{\sigma^2} \left( 2 + c^{-1/2} \right)^2 \right)^{-1}
= C_2 \ .
\end{equation} 
Choose the constant $C$ in the lefthand side of \eqref{eq-proba-Xj} as 
$C = C_2 / 4$. From \eqref{eq-X_j>S_j} we have 
\begin{eqnarray*}
\max_j \PP\left( X_j \leq C \right)
&\leq& \max_j \PP\left( S_j \leq C \right) \\
&=& \max_j \PP\left( S_j \leq C, | S_j - \rho_j | \geq C \right) + 
 \max_j \PP\left( S_j \leq C, | S_j - \rho_j | < C \right) \\
&\leq& 
\max_j \PP\left( | S_j - \rho_j | \geq C \right) + 
\max_j \PP\left( \rho_j \leq 2 C \right)  \\
&\stackrel{(a)}{\leq}& 
\frac{1}{C^2} \max_j \EE \left[ \left( S_j - \rho_j \right)^2 \right] + 
\max_j \PP\left( \rho_j \leq 2 C \right)  \\
&\stackrel{(b)}{\leq}& 
\frac{1}{C^2} \max_j \EE \left[ \left( S_j - \rho_j \right)^2 \right] + 
\PP\left( \min_j \rho_j \leq 2C \right) \\
&\stackrel{(c)}{=}& o(1) 
\end{eqnarray*} 
where $(a)$ is Tchebychev's Inequality, $(b)$ is due to $\max_j \PP({\cal E}_j) 
\leq \PP(\cup_j {\cal E}_j)$, and $(c)$ is due to \eqref{eq-tchebychev} and to 
\eqref{eq-liminf-rho}. \\  
We have proven \eqref{eq-proba-Xj} and hence that 
${\bf Q}^{\mathrm{d}}_* = \arg\max_{{\bf Q} \in {\cal C}^{\mathrm{d}}_1} 
I({\bf Q})$ satisfies $\sup_t \| \mb{Q}_*^{\mathrm{d}} \| < \infty$. \\ 
In order to prove that 
${\bf Q}_* = \arg\max_{{\bf Q} \in {\cal C}_1} 
I({\bf Q})$ satisfies $\sup_t \| \mb{Q}_* \| < \infty$, we begin by noticing that 
\begin{equation}
\label{eq-max-max} 
\max_{{\bf Q} \in {\cal C}_1} I({\bf Q}) = 
\max_{{\bf U} \in {\cal U}_t} 
\max_{{\bs \Lambda} \in {\cal C}_1^{\mathrm{d}}} 
\EE\left[ \log\det\left( {\bf I}_r + \frac{1}{\sigma^2} 
{\bf H} {\bf U} {\bs \Lambda} {\bf U}^H {\bf H}^H \right) 
\right] 
\end{equation} 
where ${\cal U}_t$ is the group of unitary $t \times t$ matrices. 
For a given matrix ${\bf U} \in {\cal U}_t$, the inner maximization in 
\eqref{eq-max-max} is equivalent to the problem of maximizing the mutual information over 
${\cal C}_1^{\mathrm{d}}$ when the channel matrix ${\bf H}$ is  replaced with 
${\bf H}' = {\bf H} {\bf U} = \sqrt{\frac{K}{K+1}} {\bf A}' +
\frac{1}{\sqrt{K+1}} \frac{1}{\sqrt{t}} {\bf C}^{1/2} {\bf W}^{'} \tilde{\bf C}'^{1/2}$.
Here, matrix $\tilde{{\bf C}}^{'}$ is defined by $\tilde{{\bf C}}^{'} = {\bf U}^{H} \tilde{{\bf C}} {\bf U}$, 
${\bf A}^{'} = {\bf A} {\bf U}$, ${\bf W}^{'} = {\bf W} {\bs \Theta}$ where ${\bs \Theta}$ is the unitary matrix
${\bs \Theta} = \tilde{{\bf C}}^{1/2} {\bf U} \tilde{{\bf C}}^{'-1/2}$. 
As ${\bf U} \in {\cal U}_t$, we clearly have 
$\| {\bf A}' \| = \| {\bf A} \|$, 
$\| \tilde{\bf C}' \| = \| \tilde{\bf C} \|$, and 
$\| \tilde{\bf C}'^{-1} \| = \| \tilde{\bf C}^{-1} \|$. 
By consequence, the bounds $a$ and $A$, and hence the constant $C$ in 
the left hand member of \eqref{eq-proba-Xj} (which depends only on $(a,A,\sigma^2,c,K)$)
remain unchanged when we replace ${\bf H}$ with ${\bf H}'$. By consequence, for every
${\bf U} \in {\cal U}_t$ the matrix 
${\bs \Lambda}_*({\bf U})$ that maximizes 
$\EE\left[ \log\det\left( {\bf I}_r + \frac{1}{\sigma^2}
{\bf H} {\bf U} {\bs \Lambda} {\bf U}^H {\bf H}^H \right)
\right]$ satisfies $ \| {\bs \Lambda}_*({\bf U}) \| < 2 \sigma^2 / C +1$ (see
\eqref{eq-borne-Q}) which is independent of ${\bf U}$. 
Hence $\| {\bf Q}_* \| < 2 \sigma^2 / C +1$ which terminates the proof of (ii).

\bibliographystyle{plain}

%\bibliography{biblio}

%\bibliography{math}

\end{document}